\documentclass[11pt]{amsart}
\relax 
\usepackage{latexsym}
\usepackage{amssymb}
\usepackage{amscd}
\usepackage[all]{xy}
\usepackage[mathscr]{euscript}
\usepackage{dsfont}
\usepackage{hyperref}
\usepackage{graphicx}
\usepackage{amsfonts}
\usepackage{amsmath}
\usepackage{amsthm}
\usepackage{aliascnt}
\usepackage{color}

\normalsize

\RequirePackage{xspace}

\newcommand{\thought}[1]{}
\renewcommand{\thought}[1]{ \textbf{[#1]}}

\usepackage{enumerate}        
\newenvironment{roenumerate}{\begin{enumerate}[\upshape (i)]}{\end{enumerate}}

\newcommand\nc {\newcommand}
\newcommand\rnc{\renewcommand}

\newcount\blopone
\newcount\xone
\newcount\xtwo
\newcount\ytwo

\newtheorem{theorem}{Theorem}[section]
\newtheorem{prop}[theorem]{Proposition}
\newtheorem{com}[theorem]{Comment}
\newtheorem{apl}[theorem]{Application}
\newtheorem{exercise}[theorem]{Exercise}
\newtheorem{redu}[theorem]{Reduction}
\newtheorem{refinement}[theorem]{Refinement}
\newtheorem{summary}[theorem]{Summary}
\newtheorem{importnota}[theorem]{Important Notation}
\newtheorem{prblm}[theorem]{Problem}
\newtheorem{notation}[theorem]{Notation}
\newtheorem{explanation}[theorem]{Explanation}
\newtheorem{defin}[theorem]{Definition}
\newtheorem{caution}[theorem]{Caution}
\newtheorem{remark}[theorem]{Remark}
\newtheorem{reminder}[theorem]{Reminder}
\newtheorem{illustration}[theorem]{Illustration}
\newtheorem{observation}[theorem]{Observation}
\newtheorem{lemma}[theorem]{Lemma}
\newtheorem{construction}[theorem]{Construction}
\newtheorem{discussion}[theorem]{Discussion}
\newtheorem{corollary}[theorem]{Corollary}
\newtheorem{example}[theorem]{Example}
\newtheorem{conclusion}[theorem]{Conclusion}
\newtheorem{sketch}[theorem]{Sketch}
\newtheorem{triviality}[theorem]{Triviality}
\newtheorem{proto}[theorem]{Prototype Quasifibration}
\newtheorem{cauex}[theorem]{Cautionary Example}
\newtheorem{hypo}[theorem]{Hypothesis}
\newtheorem{subth}{ }[theorem]
\newtheorem{case}{Case}[theorem]
\newtheorem{ssubth}{ }[subth]
\newtheorem{facts}[theorem]{Facts}
\newtheorem{history}[theorem]{Historical Survey}
\newtheorem{proofs}[theorem]{Discussion of the Proofs, Old and New}
\newtheorem{discl}[theorem]{Disclaimer}

\nc\tri[1]{\begin{triviality}
\label{#1}}
\nc\fac[1]{\begin{facts}
\label{#1}
\begin{em}}
\nc\app[1]{\begin{apl}
\label{#1}
\begin{em}}
\nc\skt[1]{\begin{sketch}
\label{#1}
\begin{em}}
\nc\hst[1]{\begin{history}
\label{#1}
\begin{em}}
\nc\pfs[1]{\begin{proofs}
\label{#1}
\begin{em}}
\nc\cas[1]{\begin{case}
\label{#1}
\begin{em}}
\nc\rfn[1]{\begin{refinement}
\label{#1}}
\nc\prt[1]{\begin{proto}
\label{#1}}
\nc\lem[1]{\begin{lemma}
\label{#1}}
\nc\pro[1]{\begin{prop}
\label{#1}}
\nc\thm[1]{\begin{theorem}
\label{#1}}
\nc\dis[1]{\begin{discussion}
\label{#1}
\begin{em}}
\nc\dsc[1]{\begin{discl}
\label{#1}
\begin{em}}
\nc\cor[1]{\begin{corollary}
\label{#1}}
\nc\dfn[1]{\begin{defin}
\label{#1}}
\nc\sthm[1]{\begin{subth}
\label{#1}}
\nc\exm[1]{\begin{example}
\label{#1}
\begin{em}}
\nc\obs[1]{\begin{observation}
\label{#1}
\begin{em}}
\nc\plm[1]{\begin{prblm}
\label{#1}
\begin{em}}
\nc\rmk[1]{\begin{remark}
\label{#1}
\begin{em}}
\nc\rmd[1]{\begin{reminder}
\label{#1}
\begin{em}}
\nc\ntn[1]{\begin{notation}
\label{#1}
\begin{em}}
\nc\exe[1]{\begin{exercise}
\label{#1}
\begin{em}}
\nc\xpl[1]{\begin{explanation}
\label{#1}
\begin{em}}
\nc\smr[1]{\begin{summary}
\label{#1}
\begin{em}}
\nc\cau[1]{\begin{caution}
\label{#1}
\begin{em}}
\nc\hyp[1]{\begin{hypo}
\label{#1}}
\nc\imn[1]{\begin{importnota}
\label{#1}
\begin{em}}
\nc\rdn[1]{\begin{redu}
\label{#1}
\begin{em}}
\nc\cax[1]{\begin{cauex}
\label{#1}
\begin{em}}
\nc\cmt[1]{\begin{com}
\label{#1}
\begin{em}}
\nc\con[1]{\begin{construction}
\label{#1}
\begin{em}}
\nc\ill[1]{\begin{illustration}
\label{#1}
\begin{em}}
\nc\ssthm[1]{\begin{ssubth}
\label{#1}
\begin{em}}
\nc\cnc[1]{\begin{conclusion}
\label{#1}
\begin{em}}

\nc\elem{\end{lemma}}
\nc\erdn{\end{em}\end{redu}}
\nc\erfn{\end{refinement}}
\nc\eprt{\end{proto}}
\nc\ethm{\end{theorem}}
\nc\ecor{\end{corollary}}
\nc\edfn{\end{defin}}
\nc\esthm{\end{subth}}
\nc\epro{\end{prop}}
\nc\etri{\end{triviality}}
\nc\eexm{\end{em}
\end{example}}
\nc\eobs{\end{em}
\end{observation}}
\nc\ecmt{\end{em}
\end{com}}
\nc\efac{\end{em}
\end{facts}}
\nc\eapp{\end{em}
\end{apl}}
\nc\ermk{\end{em}
\end{remark}}
\nc\ermd{\end{em}
\end{reminder}}
\nc\eill{\end{em}
\end{illustration}}
\nc\eplm{\end{em}
\end{prblm}}
\nc\ecas{\end{em}
\end{case}}
\nc\eskt{\end{em}
\end{sketch}}
\nc\ecau{\end{em}
\end{caution}}
\nc\ecax{\end{em}
\end{cauex}}
\nc\eimn{\end{em}
\end{importnota}}
\nc\entn{\end{em}
\end{notation}}
\nc\eexe{\end{em}
\end{exercise}}
\nc\expl{\end{em}
\end{explanation}}
\nc\edis{\end{em}
\end{discussion}}
\nc\edsc{\end{em}
\end{discl}}
\nc\econ{\end{em}
\end{construction}}
\nc\esmr{\end{em}
\end{summary}}
\nc\ehst{\end{em}
\end{history}}
\nc\epfs{\end{em}
\end{proofs}}
\nc\ehyp{
\end{hypo}}
\nc\ecnc{\end{em}
\end{conclusion}}
\nc\essthm{\end{em}
\end{ssubth}}

\nc\sst{\scriptstyle}
\newcommand{\comment}[1]{}
\newcommand{\ri}{\longrightarrow}
\newcommand{\sr}{\rightarrow}

\newcommand{\zz}{{\mathbb Z}}

\newcommand{\nn}{{\mathbb N}}
\newcommand{\K}{{\mathbf K}}
\newcommand{\D}{{\mathbf D}}

\nc\op{^{\hbox{\rm\tiny op}}}
\nc\mth{^{\hbox{\rm\tiny th}}}

\nc\script{\mathscr}
\nc\z{\zeta}
\nc\bc{{\mathbb{BC}}}
\nc\ct{{\script T}}
\nc\cf{{\script F}}
\nc\cg{{\script G}}
\nc\ch{{\script H}}
\nc\ck{{\script K}}
\nc\cl{{\script L}}
\nc\cv{{\script V}}
\nc\ce{{\script E}}
\nc\cm{{\script M}}
\nc\cn{{\script N}}
\nc\cs{{\script S}}
\nc\car{{\script R}}
\nc\cd{{\script D}}
\nc\cc{{\script C}}
\nc\ca{{\script A}}
\nc\ci{{\script I}}
\nc\co{{\script O}}
\nc\cu{{\script U}}
\nc\cx{{\script X}}
\nc\cy{{\script Y}}
\nc\cz{{\script Z}}
\nc\Cp{{\script P}}
\nc\cq{{\script Q}}
\nc\bd{\begin{description}}
\nc\ed{\end{description}}
\nc\ctob{{\script C}at\big(\ci^{op},\ca\big)}
\nc\clim{{\ds\mathop{\rm lim}_{\ds\longleftarrow}}\,}
\nc\climone{{\ds{\mathop{\rm lim}_{\ds\longleftarrow}}}^1\,}
\nc\climi{\clim_{\!i}\,}
\nc\climn{\clim^{\!n}\,}
\nc\colim{{\ds\mathop{\rm colim}_{\ds\la}}}
\nc\colimj{{\ds\mathop{\rm colim}_{\ds\la}}{}_{j\,}}
\nc\oa{\overline{\ca}}
\nc\s{\sigma}
\nc\ta{\tau}
\nc\os{\overline\sigma}
\nc\ot{\overline\tau}
\nc\T{\Sigma}
\nc\Tm{\Sigma^{-1}}
\nc\de[1]{{\mathop{\rm deg(#1)}}}
\nc\Ad[1]{\mathop{\rm Ad}(#1)}
\nc\ad[1]{\mathop{\rm ad}(#1)}
\nc\kth{{\it K}--theory}
\nc\loc[1]{{\text{\rm Loc}(#1)}}
\nc\coloc[1]{{\text{\rm Coloc}(#1)}}

\def\der #1 {D\left(#1\right)}
\nc\prf{\begin{proof}}
\nc\eprf{\end{proof}}
\nc\ds{\displaystyle}
\nc\Tor{\text{\rm Tor}}

\nc\cb{{\script B}}
\nc\ab{{\script A}b}

\nc\be{\begin{roenumerate}}
\nc\ee{\end{roenumerate}}

\nc\cat[1]{{\script C}at\Big({\big\{#1\big\}}\op\,\,,\,\,\ab\Big)}
\nc\csab{{\script C}at\big(\cs^{op},\ab\big)}
\nc\ctab{{\script C}at\Big({\{\ct^\alpha\}}^{op},\ab\Big)}
\nc\csex{{\script E}x\big(\cs^{op},\ab\big)}
\nc\ctex{{\script E}x\Big({\{\ct^\alpha\}}^{op},\ab\Big)}
\nc\sub{\qquad\subset\qquad}
\nc\ctr[1]{{\left.\ct\left(-,#1\right)\right|}_{\cs}}
\nc\ctrf[2]{{\left.\ct\left(#1,#2\right)\right|}_{\cs}}
\nc\Ctr[1]{{\left.\ct\left(-,#1\right)\right|}_{\ct^\alpha}}
\nc\Ctrf[2]{{\left.\ct\left(#1,#2\right)\right|}_{\ct^\alpha}}

\nc\la{\longrightarrow}
\nc\nin{\noindent}
\nc\cad[1]{\text{card}(#1)}
\nc\eq{\quad=\quad}
\nc\BA{\begin{array}{c}}
\nc\EA{\end{array}}
\nc\barr{
\[
\begin{array}{cccccccccccccccc}
}
\nc\earr{
\end{array}
\]
}
\nc\as[1]{{\langle S\rangle}^{#1}}
\nc\sh{\text{\it shift}}

\nc\yy[1]{{\left.\ct\left(-,#1\right)\right|}_{\ct^c}}
\nc\vrep[2]{{\left.\ct\left(#1,#2\right)\right|}_{\ct^\alpha}}
\nc\da{\downarrow}
\nc\Hom{{\mathop{\rm Hom}}}
\nc\HHom{{\script H}{\mathop{\rm om}}}
\nc\End{{\mathop{\rm End}}}
\nc\Ext{{\mathop{\rm Ext}}}
\nc\mr\Modtc
\nc\PExt{{\mathop{\rm PExt}}}
\nc\stm{\text{\rm stmod}(kG)}
\nc\stM{\text{\rm StMod}(kG)}
\nc\e{\varepsilon}
\nc\p{\varphi}

\nc\rs{\s^{-1}A}
\nc\br{{\{\s^{-1}A\}}}
\nc\ra\ri
\nc\y[1]{\mathbf{y}#1}
\nc\x[1]{\mathbf{z}#1}
\nc\mmod[1]{#1\text{--\rm mod}}
\nc\Mod[1]{#1\text{--\rm Mod}}
\nc\Md {\ensuremath{\mathop{\textup{Mod}}}}
\rnc\mod[1]{\ensuremath{\mathop{#1\textup{--mod}}}\xspace}
\nc\MMod[1]{\text{\rm Mod-}#1}
\nc\Modtc{\Mod{\ct^c}}
\nc\pgldim[1]{\mathop{\rm pgldim}\,#1}
\nc\tf{{\rm [TR5]}}
\nc\tfs{{\rm [TR5$^*$]}}
\nc\Fun{\text{\rm Funct}(F\op,\ab)}
\nc\sym{\text{\rm Sym}}
\nc\sgn{\text{\rm sgn}}
\nc\Pro{\text{\rm Prod}^{}_\alpha(F\op,\ab)}
\nc\Yt[1]{{\left.\Hom_\ct^{}\left(-,#1\right)\right|}_F^{}}
\nc\dl{\delta}
\nc\Proj[1]{#1\text{--\rm Proj}}
\nc\proj[1]{#1\text{--\rm proj}}
\nc\Flat[1]{#1\text{--\rm Flat}}
\nc\Inj[1]{#1\text{--\rm Inj}}
\nc\Ima{\mathrm{Im}}
\nc\Ker{\mathrm{Ker}}
\nc\ov{\overline}
\nc\wt{\widetilde}
\nc\wh{\widehat}
\nc\ph{\varphi}
\nc\tstr{{\it t}--structure}
\nc\tst[1]{\left({#1}^{\leq0},{#1}^{\geq0}\right)}
\nc\tstv[2]{\left({#1}_{#2}^{\leq0},{#1}_{#2}^{\geq0}\right)}
\nc\tsth[2]{{#1}_{#2}^{\heartsuit}}

\nc\spec[1]{{\text{\rm Spec}(#1)}}

\newcommand{\fc}{\mathfrak{C}}
\newcommand{\fl}{\mathfrak{L}}
\newcommand{\fs}{\mathfrak{S}}
\nc\EProd{\text{\rm EProd}}
\nc\ECoprod{\text{\rm ECoprod}}
\nc\Prod{\text{\rm Prod}}
\nc\ldimp{\text{\rm LDim}^{\prod}}
\nc\ldimc{\text{\rm LDim}^{\coprod}}
\nc\gen[2]{{\langle#1\rangle}^{}_{#2}}
\nc\genu[3]{{\langle#1\rangle}^{[#3]}_{#2}}
\nc\ogen[1]{\ov{\langle#1\rangle}}
\nc\ogenun[2]{\ov{\langle#1\rangle}_{#2}^{}}
\nc\ogenu[3]{\ov{\langle#1\rangle}^{[#3]}_{#2}}
\nc\ogenul[3]{\ov{\langle#1\rangle}^{(-\infty,#3]}_{#2}}
\nc\ogenuf[3]{\ov{\langle#1\rangle}^{[#3,\infty)}_{#2}}
\nc\genuf[3]{{\langle#1\rangle}^{[#3,\infty)}_{#2}}
\nc\genul[3]{{\langle#1\rangle}^{(-\infty,#3]}_{#2}}
\nc\dperf[1]{\D^{\mathrm{perf}}(#1)}
\nc\dcoh{\mathbf{D}^b_{\mathrm{coh}}}
\nc\dperfs[2]{\D_{#1}^{\mathrm{perf}}(#2)}
\nc\dcohs[1]{\mathbf{D}^b_{\mathrm{coh},#1}}

\newcommand{\Dqcs}[1]{{\mathbf D_{\text{\bf qc},#1}}}

\nc\dmcoh{\mathbf{D}^-_{\mathrm{coh}}}
\nc\dmcohs[1]{\mathbf{D}^-_{\mathrm{coh,#1}}}
\nc\dscoh{\mathbf{D}^{}_{\mathrm{coh}}}
\nc\RHHom{{\script{RH}}{\mathrm{om}}}
\nc\Coprod{\mathrm{Coprod}}
\nc\COprod{\mathrm{coprod}}
\nc\add{\mathrm{add}}
\nc\Add{\mathrm{Add}}
\nc\Smr{\mathrm{smd}}
\nc\id{\mathrm{id}}
\nc\LL{\mathbf{L}}
\nc\R{\mathbf{R}}
\nc\tsb{\ct^{\mathrm{sb}}}
\nc\Projdg[1]{\text{\rm Proj}^{}_{\mathbf{dg}}#1}

\renewcommand{\leq}{\leqslant}
\renewcommand{\geq}{\geqslant}

\nc\hoco{
\begin{picture}(40,10)
\put(20,0){\makebox(0,0)[b]{\text{\rm Hocolim}}}
\put(5,-2){\vector(1,0){30}}
\end{picture}\,\,}

\nc\holim{
\begin{picture}(40,10)
\put(20,0){\makebox(0,0)[b]{\text{\rm Holim}}}
\put(35,-2){\vector(-1,0){30}}
\end{picture}}

\begin{document}

\author{Amnon Neeman}\thanks{
  The research was supported by
the ERC Advanced Grant 101095900-TriCatApp.
}
\address{Dipartimento di Matematica ``F.\ Enriques''\\
        Universit{\`a} degli Studi di Milano\\
        Via Cesare Saldini 50\\
	20133 Milano\\
        ITALY}
\email{amnon.neeman@unimi.it}

\title[Excellent metrics]{Excellent metrics on triangulated
  categories, and the involutivity of the map taking
  $\cs$ to $\fs(\cs)\op$}

\begin{abstract}
In the article
\cite{Neeman18A} we defined
good metrics on triangulated
categories, and then studied the
construction, that began with a triangulated
category $\cs$ together with a good metric
$\{\cm_i,\,i\in\nn\}$, and out of it
cooked up another triangulated category
$\fs(\cs)$. We went on
to study examples, and
produced many
for which the construction is involutive.
By this we mean
that, if you let $\ct=\fs(\cs)\op$,
then there is
a choice of metric on $\ct$ for
which $\cs=\fs(\ct)\op$.

In this article we study this phenomenon much
more carefully, with the focus
being on understanding the metrics for which
involutivity occurs.
As it turns out there is a large class
of them, the excellent metrics on
triangulated categories.

At the end we will produce a few new examples
of excellent metrics. And our reason
for going to all this trouble is that
the results of this article will
permit us to prove new and surprising
statements about uniqueness of enhancements.
Those results will come in a sequel to this
article,
which is
joint with Canonaco and Stellari.
\end{abstract}

\subjclass[2020]{Primary 18G80}

\keywords{Derived categories, {\it t}--structures, homotopy limits, metrics, completions}

\maketitle

\tableofcontents

\setcounter{section}{-1}

\section{Introduction}
\label{S900}

The article is a sequel to~\cite{Neeman18A},
which
was written to try to better understand
Rickard~\cite[Theorem~6.4]{Rickard89b}.
Rickard's old theorem asserts, among other
things, that
if $R$ and $S$ are left-coherent rings, then
\[
\D^b(\proj R)\cong \D^b(\proj S) \quad\Longleftrightarrow\quad
\D^b(\mod R)\cong \D^b(\mod S)\ .
\]
To elaborate, in words rather than symbols:
the $R$-perfect complexes are triangle equivalent with
the $S$-perfect complexes
if and only if
the bounded complexes of finitely generated
$R$-modules are triangle equivalent to the
bounded complexes of finitely generated
$S$-modules.
And the short summary is that the article
\cite{Neeman18A}
does three things, the first of which is:

\smr{S900.1}
In \cite[Definition~\ref{D20.1}]{Neeman18A}
we introduce \emph{good metrics} on
triangulated categories.
Given a triangulated category
$\cs$, together with a good metric
$\{\cm_i,\,i\in\nn\}$,
then
\cite[Definition~\ref{D20.11}(iii)]{Neeman18A}
produces a new category
$\fs(\cs)$---this category
depends on the choice of metric
$\{\cm_i,\,i\in\nn\}$, but the notation
suppresses this dependence. Writing
something like
$
\fs_{\{\cm_i,\,i\in\nn\}}^{}(\cs)
$
would be too cumbersome.

Next
\cite[Definition~\ref{D28.110}]{Neeman18A}
goes on to define what are the
distinguished triangles in
the category $\fs(\cs)$,
and finally
\cite[Theorem~\ref{T28.128}]{Neeman18A}
proves that
the category $\fs(\cs)$ is
triangulated.
\esmr

The next step is to establish the relevance
to Rickard's old theorem, and
for this we proved:

\smr{S900.3}
Let $R$ be a left-coherent ring.
There is an equivalence class
of metrics on the triangulated
category $\D^b(\proj R)$,
and an equivalence class of
metrics on $\D^b(\mod R)\op$, as
well as
triangle equivalences
\[
\fs\big(\D^b(\proj R)\big)\cong
\D^b(\mod R)\qquad\text{ and }\qquad
\fs\big(\D^b(\mod R)\op\big)\cong
\D^b(\proj R)\op\ .
\]
The relevant metrics are obtained
by specializing 
\cite[Example~\ref{E20.3}]{Neeman18A}
to the case where $\ct=\D(\Mod R)$
and therefore $\ct^c=\D^b(\proj R)$
and $\ct^b_c=\D^b(\mod R)$,
and the triangle equivalences
above
are obtained, again by
specializing to the case
$\ct=\D(\Mod R)$, the computations
of
\cite[Example~\ref{E22.3} and
  Proposition~\ref{P29.9}]{Neeman18A}
(respectively).
\esmr

Since Rickard's theorem makes no
mention of any metrics, the final
step in 
\cite{Neeman18A}
is:

\smr{S900.5}
There are recipes that start with a triangulated
category and produce metrics, well-defined
up to equivalence. And in particular
the two metrics used
in Summary~\ref{S900.3}, the one on
$\D^b(\proj R)$ and the one
on $\D^b(\mod R)\op$,
can each be given by such a recipe.
Thus the metrics of 
Summary~\ref{S900.3} do not
constitute added structure.

For the reference: the ``intrinsicness''
of the metric on $\D^b(\proj R)$
follows from
\cite[Remark~\ref{R22.105.5} and
  Proposition~\ref{P22.109}]{Neeman18A},
while the
``intrinsicness''
of the metric on $\D^b(\mod R)\op$
is proved in
\cite[Proposition~\ref{P1.105}]{Neeman18A}.
And once again: what is proved in
\cite{Neeman18A} is far more general,
we are specializing the results of
\cite{Neeman18A} to the case
$\ct=\D(\Mod R)$.
\esmr

\rmk{R900.7}
Let $\cs$ be either
$\D^b(\proj R)$ or $\D^b(\mod R)$.
From the intrinsic construction of the
metrics, of Summary~\ref{S900.5},
it follows that any triangulated
autoequivalence of $\cs$ carries the
metric to an equivalent one, and hence
doesn't change $\fs(\cs)$.
Thus an immediate consequence, of
the theorems of \cite{Neeman18A},
is that the group of autoequivalences
of $\cs$ is isomorphic to the group
of autoequivalences of $\fs(\cs)$.
\ermk

As we said in Summary~\ref{S900.5}:
the results
of \cite{Neeman18A} are far more
general, they apply to many categories
$\cs$ which are neither
$\D^b(\proj R)$ nor $\D^b(\mod R)$.
And the statement
in Remark~\ref{R900.7}
is
equally generalizable.
This is evident already
from the generality of
the results
in \cite{Neeman18A},
and the current article will
generalize the theory far beyond
the previous scope.

In view of the phenomena discovered
in \cite{Neeman18A}, it is natural to ask
the following question:

\plm{P900.9}
Are there reasonable conditions guaranteeing that,
if the triangulated
category $\cs$ has a unique 
enhancement, then so does
the triangulated category $\fs(\cs)$?
\eplm

This is the question that motivated this article.
As we will show in a sequel
(joint with Canonaco and Stellari)
the answer turns out to be a resounding Yes.

\rmk{R900.11}
Of course: to thoroughly
prepare the ground we need
to understand
those triangulated categories
$\cs$,
which admit good metrics
$\{\cm_i,\,i\in\nn\}$, for which
the passage from $\cs$ to $\fs(\cs)\op$
is an involution. Being a good metric
doesn't suffice, it is easy to
come up with counterexamples.
Explicitly:
\cite[Example~\ref{E20.2}]{Neeman18A}
shows that every triangulated category
$\cs$ can be endowed with a dumb metric,
and
\cite[Example~\ref{E20.1098}]{Neeman18A}
teaches us that, for the dumb metric,
we obtain $\fs(\cs)=\{0\}$. Therefore
the dumb metric cannot
be involutive unless $\cs=\{0\}$,
and so the
involutive property
can only be true 
for a special class of metrics.
This article sets out to
understand this class.

From a philosophical point of view,
this article takes an
approach diametrically opposite
to the
one in
\cite{Neeman18A}. In
\cite{Neeman18A} the focus in on recovering
and generalizing the old theorem of
Rickard's, which predates good
metrics and makes no mention of
them. And hence \cite{Neeman18A} goes to
some length to show that
the main metrics of interest
are intrinsic, as explained in
Summary~\ref{S900.5}:
there is a recipe that
obtains the relevant metrics from only the
triangulated category. In this article
the focus is on the metrics, and on
trying to identify what singles out the
involutive ones.
\ermk

The new results in the current paper
formalize this, starting with the following.

\pro{P900.13}
Let $\cs$ be a triangulated category,
and suppose we are given on $\cs$
a good metric $\{\cm_i\subset\cs,\, i\in\nn\}$.
Then in Definition~\ref{D3.-3} we introduce
a string of subcategories
$\{\cn_i\subset\fs(\cs),\,i\in\nn\}$,
and Lemma~\ref{L3.1} proves that
this is a good metric on
$\fs(\cs)$.
\epro

We note in passing: it is obvious
from the old
\cite[Definition~\ref{D20.1}]{Neeman18A}
that a good metric on $\cs$
is the same thing as a good
metric on
$\cs\op$. Therefore
$\{\cn\op_i,\,i\in\nn\}$
is a good metric
on $\fs(\cs)\op$.

This allows us to state a revised
version of Problem~\ref{P900.9}:

\plm{P900.15}
What are the conditions one needs to
impose on the pair
\[
\cs,\quad\{\cm_i,\,i\in\nn\}
\]
of a triangulated category $\cs$ and a good
metric $\{\cm_i,\,i\in\nn\}$,
to guarantee that the passage to the pair
\[
\fs(\cs)\op,\quad\{\cn\op_i,\,i\in\nn\}
\]
is an involution?
\eplm

And the main results of the article
can be summarized by saying that we
find a very general set of conditions
which are sufficient for this to happen.
That is we prove

\thm{T900.17}
Inside the class of good metrics on
triangulated categories there is
a subclass, called the \emph{excellent}
metrics,
which satisfy
\be
\item
If $\cs$ is a triangulated
category and 
$\{\cm_i,\,i\in\nn\}$ is an excellent metric
on $\cs$, then
$\{\cn\op_i,\,i\in\nn\}$
is an excellent metric
on $\fs(\cs)\op$.
\item
For triangulated catgeories
$\cs$ with excellent metrics
$\{\cm_i,\,i\in\nn\}$,
the
passage from the pair
\[
\cs,\quad\{\cm_i,\,i\in\nn\}
\]
to the pair
\[
\fs(\cs)\op,\quad\{\cn\op_i,\,i\in\nn\}
\]
is \emph{almost} an involution.
\setcounter{enumiv}{\value{enumi}}
\ee
\emph{We need to elaborate, and explain what it
means for a map to be
``almost'' an involution.}
The meaning is:
\be
\setcounter{enumi}{\value{enumiv}}
\item
There is a fully faithful,
triangulated functor 
$\wh Y:\cs\la\fs\big(\fs(\cs)\op\big)\op$,
such that every object of
$\fs\big(\fs(\cs)\op\big)\op$
is a direct summand of an object
in the image of $\wh Y$.
\item
Let $\cs$ be a triangulated
category with an excellent
metric
$\{\cm_i,\,i\in\nn\}$.
The passage takes the pair
\[
\cs,\quad\{\cm_i,\,i\in\nn\}
\]
to the pair
\[
\fs(\cs)\op,\quad\{\cn\op_i,\,i\in\nn\}
\]
which in turn gets taken to a
pair we will denote
\[
\fs\big(\fs(\cs)\op\big)\op,\quad\{\wh\cm_i,\,i\in\nn\}\ .
\]
The functor
$\wh Y:\cs\la\fs\big(\fs(\cs)\op\big)\op$,
of (iii) above, satisfies $\wh Y(\cm_i)\subset\wh\cm_i$.
Moreover: every object in
$\wh\cm_i$ is the direct summand of
an object in $\wh Y(\cm_i)$.
\item
Finally: this ``idempotent completion''
process terminates after the first step.
That is: if the category $\cs$, together
with its metric, is equal
to $\fs(\car)$ for some triangulated
category $\car$ with an excellent metric, then
the functor
$\wh Y:\cs\la\fs\big(\fs(\cs)\op\big)\op$
is a triangle equivalence of
triangulated categories, and is an isometry
with respect to the metrics.
\ee
\ethm

We have so far explained the bare essentials,
the key highlights of the current article.
But the article contains much more, and
we will now attempt to summarize this.

\rmd{R900.19}
Let $\cs$ be a triangulated category,
and let $\{\cm_i,\,i\in\nn\}$
be a good metric on $\cs$.
Going back to \cite{Neeman18A}: an
important tool we used was to consider
the Yoneda embedding $Y:\cs\la\MMod\cs$,
the functor sending an object
$X\in\cs$ to the representable functor
$\Hom(-,X):\cs\op\la\ab$, which we view
as an object in $\MMod\cs$.
And then the idea was to study
various subcategories of
$\MMod\cs$ and the relations among them.

To be specific:
back
in \cite[Definition~\ref{D20.11}]{Neeman18A}
we introduced 
the full subcategory $\fl(\cs)\subset\MMod\cs$,
which 
has for objects the colimits
of the image under $Y$ of
Cauchy sequences
in $\cs$. And the full subcategory
$\fs(\cs)\subset\fl(\cs)$ is given by the formula
\[
\fs(\cs)\eq\fl(\cs)\cap\bigcup_{n\in\nn}Y(\cm_i)^\perp\ .
\]
In other word: the category
$\fs(\cs)$ of the last few pages  isn't
mysterious, we have just spelled out
what it is
as a subcategory of $\MMod\cs$.
\ermd

\dis{D900.21}
Now the metric on $\cs$ is given by
full subcategories $\cm_n\subset\cs$,
which in several surveys on the subject
have been interpreted as
the balls, centered at the origin, of
radius $2^{-n}$. The reader is
referred to
\cite{Neeman19,Neeman22}
for elaboration on this intuition.
Now Definition~\ref{D3.-3}
of the current article
forms full subcategories
$\cl_n\subset\fl(\cs)$ and $\cn_n\subset\fs(\cs)$,
respectively,
by the rules
\[
\cl_n\eq
\left\{X\in\MMod\cs
\left|
\begin{array}{c}
X=\colim\,Y(x_*)\\
\text{with }x_*\subset\cm_n\text{ a Cauchy sequence}
\end{array}
\right.\right\}
\]
and
\[
\cn_n=\fs(\cs)\cap\cl_n\ .
\]
And then Lemma~\ref{L3.1} teaches us that,
on the triangulated category $\fs(\cs)$,
the subcategories $\{\cn_n,\,n\in\nn\}$
define a good metric.
Note that $\fl(\cs)$ is rarely
a triangulated category, and hence there
is no
notion of a good metric on
$\fl(\cs)$.

These definitions combine to make
$\fs(\cs)$, together with its
good metric $\{\cn_n,\,n\in\nn\}$,
a metric triangulated category. This render
meaningful the question posed in
Problem~\ref{P900.15}: we can wonder when
the passage is involutive.
\edis

In \cite[Definition~\ref{D28.907}]{Neeman18A}
we introduced the notion
of \emph{strong triangles}
in the category $\fl(\cs)$;
those are the colimits of the image
under the functor $Y:\cs\la\MMod\cs$
of Cauchy sequences of
distinguished triangles in $\cs$.
We showed that any morphism
$F:A\la B$, in the category
$\fl(\cs)$, can be completed to
a strong triangle
$A\stackrel F\la B\la C\la\T A$ in at
least one way. But
the category $\fl(\cs)$ isn't
triangulated, and
this completion
to a strong triangle
is not in general
unique up
to isomorphism. This turns
out to be a technical
headache. 

\dis{D900.23}
Fortunately for
us there is a way to work
around this difficulty. If we are given
not only the morphism
$F:A\la B$ in $\fl(\cs)$,
but in addition
a pair of Cauchy sequences
$a'_*$ and $b'_*$ in $\cs$ with
$A=\colim\,Y(a'_*)$
and
$B=\colim\,Y(b'_*)$,
then the following is true
\be
\item
Let $n>0$ be an
integer.  
If there exists 
one Cauchy sequence of triangles
$a_*\la b_*\la c_*\la\T a_*$,
such that $a_*$ is a subsequence
of $a'_*$ and $b_*$ is
a subsequence of $b'_*$,
and the Cauchy sequence
$c_*$ belongs
to $\cm_n$, then every
choice of
Cauchy sequence of
triangles as above has
a subsequence satisfying
the hypothesis.
\setcounter{enumiv}{\value{enumi}}
\ee
And this technical tool turns out to
be enough. It gives
sufficient control on
the strong triangles
$A\stackrel F\la B\la C\la\T A$
to allow the theory to proceed.
It leads us to a definition:
\be
\setcounter{enumi}{\value{enumiv}}
\item
Let $F:A\la B$ be a morphism
in the category $\fl(\cs)$, and
choose Cauchy sequences
$a'_*$ and $b'_*$ in $\cs$ with
$A=\colim\,Y(a'_*)$
and
$B=\colim\,Y(b'_*)$.
The morphism $F:A\la B$
is declared to be
of \emph{type-$n$ with respect
to $(a'_*,b'_*)$}
if a Cauchy sequence
of triangles as in (i) exists.
\ee
With this technical baggage established,
we are in a position to define
excellent metrics on triangulated categories.
Let $\cs$ be a triangulated
category and let $\{\cm_n,\,n\in\nn\}$
be a good metric on $\cs$. The
metric
is declared \emph{excellent} if
it satisfies the
hypotheses of Definition~\ref{D3.3}.
Let us not
spell out the definition here,
the notable part is that it asserts the
existence of some type-$n$ morphisms
in the category $\fl(\cs)$,
as in (ii) above.
\edis

\rmk{R900.25}
Now just because the category
$\fs(\cs)$ is defined
to be some subcategory of $\MMod\cs$,
the Yoneda functor taking
$B\in\MMod\cs$ to the functor
$\Hom_{\MMod\cs}^{}(B,-)\big|_{\fs(\cs)}^{}:\fs(\cs)\la\ab$
is a map we will denote
\[\xymatrix{
\wh Y\ar@{}[r]|-{\ds:} &
\big(\MMod\cs\big)\op\ar[rrr]
&&&
\MMod{\fs(\cs)\op}
}\]
And because the category $\fs(\cs)$ is
a full subcategory of $\MMod\cs$, the composite
\[\xymatrix{
\fs(\cs)\op\,\,\ar@{^{(}->}[rrr] && &
\big(\MMod\cs\big)\op\ar[rrr]^-{\wh Y}
&&&
\MMod{\fs(\cs)\op}
}\]
agrees with the usual Yoneda map
$\fs(\cs)\op\la\MMod{\fs(\cs)\op}$.
\ermk

An important technical result, which is crucial
if we are going to iterate this
process, is the following.

\lem{L900.27}
Suppose we are given a triangulated
category $\cs$ with an excellent metric
$\{\cm_i,\,i\in\nn\}$. 
Then the functor 
$\wh Y:\big(\MMod\cs\big)\op\la\MMod{\fs(\cs)\op}$
satisfies the following
two properties:
\be
\item
The functor $\wh Y$ induces
an equivalence of
$\fl(\cs)\op\subset\big(\MMod\cs\big)\op$
with $\fl\big(\fs(\cs)\op\big)\subset\MMod{\fs(\cs)\op}$.
\item
Furthermore:
$\wh Y$
takes strong triangles
in $\fl(\cs)$ to strong triangles in
$\fl\big(\fs(\cs)\op\big)$.
\ee
\elem

\rmk{R900.29}
Let us go back:
the functor $\wh Y$
permits us to compare
the constructions done
in $\MMod\cs$, which build
on the subcategory
$\cs$ and its
good metric $\{\cm_n,\,n\in\nn\}$,
with the constructions performed
in $\MMod{\fs(\cs)\op}$, building
on the subcategory $\fs(\cs)\op$ and
its 
good metric $\{\cn\op_n,\,n\in\nn\}$.
We have already seen a glimmer of this
in Lemma~\ref{L900.27}.
In pictures we obtain the diagram
where the vertical and slanted arrows are inclusions
{\color{black}
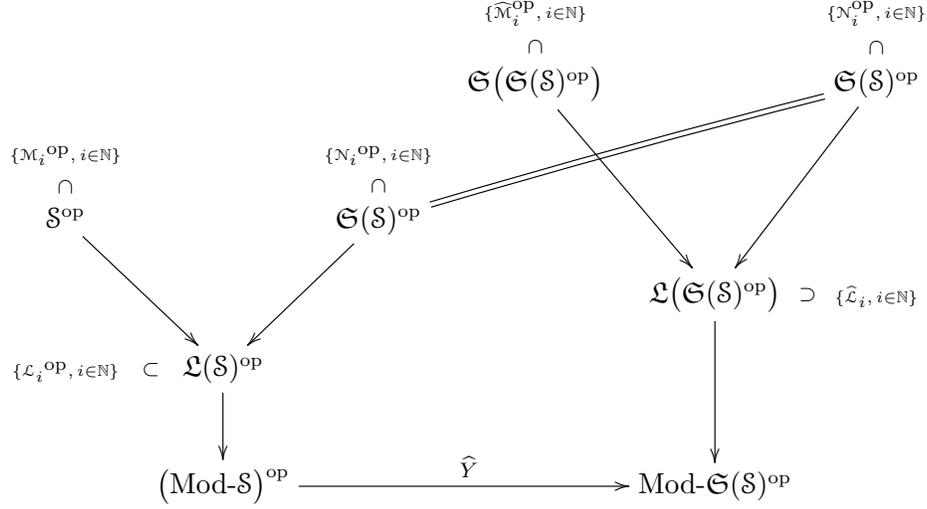
{  
\[\xymatrix@C-18pt@R-15pt{
& &  &{\scriptscriptstyle{\color{black}{\{\wh\cm\op_i,\,i\in\nn\}}}}
\ar@{}[d]|-{\color{black}{\text{\rotatebox[origin=c]{-90}{$\ds\subset$}}}}
& &{\scriptscriptstyle{\color{black}{\{\cn\op_i,\,i\in\nn\}}}}
\ar@{}[d]|-{\color{black}{\text{\rotatebox[origin=c]{-90}{$\ds\subset$}}}}
\\
& & &{\color{black}{\fs\big(\fs(\cs)\op\big)}} \ar@[black][dddr]
\ar@[white][ddddddr]
& &{\color{black}{\fs(\cs)}}{\color{black}{\op}}
\ar@[black][dddl]
\ar@[white][ddddddl]
\\
{\scriptscriptstyle\{\cm_i{\color{black}{\op}},\,i\in\nn\}}
\ar@{}[d]|-{\color{black}{\text{\rotatebox[origin=c]{-90}{$\ds\subset$}}}}
& &{\scriptscriptstyle{\color{black}{\{\cn_i{\color{black}{\op}},\,i\in\nn\}}}}
\ar@{}[d]|-{\color{black}{\text{\rotatebox[origin=c]{-90}{$\ds\subset$}}}}
& &
\\
\cs{\color{black}{\op}} \ar@[black][ddr]\ar@[white]@/^1pc/[rrruu]
\ar@[white][ddddr]_{\color{white}{Y}}
& &{\color{black}{\fs(\cs)}}{\color{black}{\op}}
\ar@[black][ddl]
\ar@[white][ddddl]
\ar@{=}@[black]@/^0.1pc/[rrruu]
& &
\\
& &&
& {\color{black}{\fl\big(\fs(\cs)\op\big)}}
\ar@[black][ddd]
&{\color{black}{\scriptscriptstyle\{\wh\cl_i,\,i\in\nn\}}}
\ar@{}[l]|-{\color{black}{\text{\rotatebox[origin=c]{180}{$\ds\subset$}}}}
\\
{\color{black}{\scriptscriptstyle\{\cl_i{\color{black}{\op}},\,i\in\nn\}}}
& {\color{black}{\fl(\cs)}}{\color{black}{\op}}
\ar@[black][dd]\ar@[white][rrru]^-{\color{white}{\ds\sim}}
\ar@{}[l]|-{\color{black}{\text{\rotatebox[origin=c]{0}{$\ds\subset$}}}}
&  & & &
  \\
& & \\
&{\color{black}{\big(}}{\color{black}{\MMod\cs}}{\color{black}{\big)\op}}
\ar@[black][rrr]^{\color{black}{\wh Y}} &&&
{\color{black}{\MMod{\fs(\cs)\op}}}
}\]}}
And what is unconditional, and free of any
hypotheses on the metric, is that the
functor
$\wh Y:\big(\MMod\cs\big)\op\la\MMod{\fs(\cs)\op}$
induces an equivalence
from the full subcategory
$\fs(\cs)\op\subset\big(\MMod\cs\big)\op$
to the full subcategory 
$\fs(\cs)\op\subset\MMod{\fs(\cs)\op}$,
and this is an isometry with respect to
the metrics $\{\cn\op_i,\,i\in\nn\}$. In the diagram this is
labeled an equality.
\ermk

And now we come to what can be proved about
the
diagram if
the metric is excellent.

\thm{T900.31}
If the metric $\{\cm_i\,i\in\nn\}$
on the category $\cs$ is excellent,
then so is the metric
$\{\cn\op_i,\,i\in\nn\}$ on the
category $\fs(\cs)\op$. Moreover:
the diagram above extends to
{\color{black}{  
\[\xymatrix@C-18pt@R-15pt{
& &  &{\scriptscriptstyle{\color{black}{\{\wh\cm\op_i,\,i\in\nn\}}}}
\ar@{}[d]|-{\color{black}{\text{\rotatebox[origin=c]{-90}{$\ds\subset$}}}}
& &{\scriptscriptstyle{\color{black}{\{\cn\op_i,\,i\in\nn\}}}}
\ar@{}[d]|-{\color{black}{\text{\rotatebox[origin=c]{-90}{$\ds\subset$}}}}
\\
& & &{\color{black}{\fs\big(\fs(\cs)\op\big)}} \ar@[black][dddr]
\ar@[white][ddddddr]
& &{\color{black}{\fs(\cs)}}{\color{black}{\op}}
\ar@[black][dddl]
\ar@[white][ddddddl]
\\
{\scriptscriptstyle\{\cm_i{\color{black}{\op}},\,i\in\nn\}}
\ar@{}[d]|-{\color{black}{\text{\rotatebox[origin=c]{-90}{$\ds\subset$}}}}
& &{\scriptscriptstyle{\color{black}{\{\cn_i{\color{black}{\op}},\,i\in\nn\}}}}
\ar@{}[d]|-{\color{black}{\text{\rotatebox[origin=c]{-90}{$\ds\subset$}}}}
& &
\\
\cs{\color{black}{\op}} \ar@[black][ddr]\ar@[black]@/^1pc/[rrruu]
\ar@[white][ddddr]_{\color{white}{Y}}
& &{\color{black}{\fs(\cs)}}{\color{black}{\op}}
\ar@[black][ddl]
\ar@[white][ddddl]
& &
\\
& &&
& {\color{black}{\fl\big(\fs(\cs)\op\big)}}
\ar@[black][ddd]&
{\color{black}{\scriptscriptstyle\{\wh\cl_i,\,i\in\nn\}}}
\ar@{}[l]|-{\color{black}{\text{\rotatebox[origin=c]{180}{$\ds\subset$}}}}
\\
{\color{black}{\scriptscriptstyle\{\cl_i{\color{black}{\op}},\,i\in\nn\}}}
& {\color{black}{\fl(\cs)}}{\color{black}{\op}}
\ar@[black][dd]\ar@[black][rrru]^-{\color{black}{\ds\sim}}
\ar@{}[l]|-{\color{black}{\text{\rotatebox[origin=c]{0}{$\ds\subset$}}}}
&  & & &
  \\
& & \\
&{\color{black}{\big(}}{\color{black}{\MMod\cs}}{\color{black}{\big)\op}}
\ar@[black][rrr]^{\color{black}{\wh Y}} &&&
{\color{black}{\MMod{\fs(\cs)\op}}}
}\]}}
\nin
We have already seen, in Lemma~\ref{L900.27},
that the functor $\wh Y$ induces
an equivalence of
$\fl(\cs)\op\subset\big(\MMod\cs\big)\op$
with $\fl\big(\fs(\cs)\op\big)\subset\MMod{\fs(\cs)\op}$.
Furthermore: this restricts to equivalences of
the subcategories $\cl\op_i\subset\fl(\cs)\op$
with the subcategories
$\wh\cl_i\subset\fl\big(\fs(\cs)\op\big)$.
One level up the diagram, the map $\wh Y$
only induces a fully faithful
triangulated functor from
$\cs\op$ into $\fs\big(\fs(\cs)\op\big)$,
and the balls $\cm\op_i$ are sent
into the balls
$\wh\cm\op_i$. And both of
these inclusions
are idempotent completions.
\ethm

And the final part of the article works out what
happens in the presence of good extensions $F:\cs\la\ct$
with respect to the metric. The notion
was introduced in \cite[Definition~\ref{D21.1}]{Neeman18A},
the point being that in the presence
of good extensions it becomes much easier to
compute what is the triangulated
category $\fs(\cs)$.

\exm{E900.?!?}
Suppose the category $\cs$ has an
enhancement---meaning there exists either some
pretriangulated dg category $\mathbf{S}$, or some
stable $\infty$-category $\mathbf{S}$, with
$\cs=H^0(\mathbf{S})$. Then the category
$\ct=H^0(\MMod{\mathbf{S}})$ is always a good
extension of $\cs$, with respect to any metric.
After all: the image of the fully faithful
inclusion $F:\cs\la\ct$ is contained
in the subcategory $\ct^c\subset\ct$ of
compact objects, and
\cite[Example~\ref{E21.3}]{Neeman18A}
permits us to deduce that $F$ is a good extension.
\eexm

\rmd{R900.33}
We recall: suppose
$\ct$ is a triangulated category
and $F:\cs\la\ct$ is
a good extension. In
\cite[Definition~\ref{D21.7}]{Neeman18A}
we met the full subcategory
$\fl'(\cs)\subset\ct$,
whose objects are the homotopy
colimits, in $\ct$, of the images under
the functor $F$ of Cauchy sequences in
$\cs$. And the functor
$\cy:\ct\la\MMod\cs$, introduced
in 
\cite[Definition~\ref{N21.-100}]{Neeman18A},
has the property that it carries
$\fl'(\cs)\subset\ct$ into
$\fl(\cs)\subset\MMod\cs$, and is almost
an equivalence of categories. It is essentially
surjective, full, and the morphisms
in $\fl'(\cs)$ annihilated by the
functor $\cy$ form a square-zero ideal.

Next we form the diagram
\[\xymatrix@C+20pt{
\cy^{-1}\fs(\cs)\,\,
\ar@{^{(}->}[r]
\ar[d] &
\fl'(\cs)\,\,\ar[d]^\cy
\ar@{^{(}->}[r]&
\ct
\\
\fs(\cs)\,\,
\ar@{^{(}->}[r]&
\fl(\cs)&
}\]
where the square is a pullback. 
And the computational highlight,
back in
\cite[Theorem~\ref{T20.17}]{Neeman18A},
can be stated as follows.
\be
\item
In the pullback
square above
the vertical morphism
on the left, meaning the functor
$\cy:\cy^{-1}\fs(\cs)\la\fs(\cs)$, is
an equivalence of categories.
\item
The composite along the top  
row, that is the inclusion
$\cy^{-1}\fs(\cs)\la\ct$,
is a triangulated functor.
\ee
We will commit the abuse of notation
of writing
$\fs(\cs)$ both for the subcategory
$\fs(\cs)\subset\fl(\cs)$, and for the
(equivalent) subcategory $\cy^{-1}\fs(\cs)\subset\fl'(\cs)$.
Accordingly
the fully faithful, triangulated
inclusion $\cy^{-1}\fs(\cs)\la\ct$ will,
with
the same abuse of notation, be
labeled $\wh F:\fs(\cs)\la\ct$.
And the point that was stressed in
the article \cite{Neeman18A}
is that this allows us to compute
$\fs(\cs)$, and its triangulated
structure, from the embedding into $\ct$.
\ermd

Everything
in Reminder~\ref{R900.33},
is valid for any triangulated category
with a good metric; excellence is not necessary
in Reminder~\ref{R900.33}.

\plm{P900.35}
This leads us to
the following, immediate
questions:
\be
\item
In the presence of a good extension, is it
any easier to recognize when a good  
metric is in fact excellent?
\item
Suppose that the functor $F:\cs\la\ct$ 
is a good extension with respect
to the metric $\{\cm_i,\,i\in\nn\}$.
Then, with the embedding
$\wh F:\fs(\cs)\la\ct$ as
in
Reminder~\ref{R900.33},
are there reasonable sufficient
conditions to ensure that
$\wh F\op:\fs(\cs)\op\la\ct\op$
is a good extension with
respect to the metric $\{\cn\op_i,\,i\in\nn\}$?
\item
Suppose we are in the
situation where both
$F:\cs\la\ct$ and $\wh F\op:\fs(\cs)\op\la\ct\op$
are good extensions, and furthermore
the metric on $\cs$
is assumed excellent.
To what extent do the pretty
diagrams of Remark~\ref{R900.29}
and Theorem~\ref{T900.31}
lift via the functor $\cy$, from the category
$\MMod\cs$ to the category $\ct$?
\ee
\eplm

\dis{D900.37}
Problem~\ref{P900.35}(i)
is easy to answer and in
Proposition~\ref{P94.3} we will,
under mild hypotheses,
give necessary and sufficient conditions 
for excellence, in terms of the good
extension. The reader should note that,
in the absence of good extensions,
the definition of excellence involves the
existence of type-$n$ morphisms in the
category $\fl(\cs)$, something that does
not lend itself to verification. Hence
finding checkable characterizations, even
ones that presuppose some given
good extension, is
unmistakeable progress.

Now for Problem~\ref{P900.35}~(ii)
and (iii). To begin with:
the subcategories $\cl_n\subset\fl(\cs)$
of Discussion~\ref{D900.21} have
obvious analogs in $\fl'(\cs)$,
we declare that $\cl'_n\subset\fl'(\cs)$
is the full subcategory whose objects
are the homotopy colimits
of Cauchy sequences in $\cm_n$.
And it can easily be shown that
\be
\item
$\cl'_n=\cy^{-1}(\cl_n)$, where the inverse
image is under the functor $\cy:\fl'(\cs)\la\fl(\cs)$.  
\item
The intersection of $\cl'_n\subset\fl'(\cs)$,
with the subcategory
$\fs(\cs)\subset\fl'(\cs)$ as embedded in
Reminder~\ref{R900.33},
agrees with the subcategory $\cn_n\subset\fs(\cs)$
of Discussion~\ref{D900.21}.
\ee
Hence
this much of the
extra data, in the diagrams of
Remark~\ref{R900.29}
and Theorem~\ref{T900.31}, is easy to lift via
the functor $\cy$.

But for a general good extension
$F:\cs\la\ct$, there is no reason for
$\wh F\op:\fs(\cs)\op\la\ct\op$ to
be a good extension. And I see no
reason for the lovely diagrams
of Remark~\ref{R900.29}
and Theorem~\ref{T900.31} to lift.
What makes these diagrams special, and
allows us to work with them,
is that we never stray far away from
the category $\cs$ and its
metric $\{\cm_i,\,i\in\nn\}$;
all the categories constructed
are subcategories of either
$\MMod\cs$ or $\MMod{\fs(\cs)\op}$.
\edis

\plm{P900.37}
The attempt to remain close,
to the category $\cs$ and its metric,
leads us to formulate the following questions:
\be
\item
Are there reasonable conditions,
on the good embedding $F:\cs\la\ct$,
which guarantee that the subcategory
$\fl'(\cs)\subset\ct$ is a triangulated
subcategory, and that $\{\cl'_i,\,i\in\nn\}$
is a good metric on $\fl'(\cs)$?
\item
Assuming the answer to question (i) is
Yes: are there reasonable conditions
that guarantee that  
the embedding $\cs\la\fl'(\cs)$ is also
a good extension?
\item
Assume we are
in the situation of (ii), and assume
further that the metric on $\cs$ is excellent.
Do the diagrams
of Remark~\ref{R900.29}
and Theorem~\ref{T900.31}
have parallels in $\fl'(\cs)$?
\ee
\eplm

\dis{D900.39}
When phrased this way, the answer turns out to
be Yes. Of course: the necessary conditions
on $\cs$ and on its metric, and on the good
extension $F:\cs\la\ct$, are all relatively
restrictive. But the highlight is
that when $F:\cs\la\cl'(\cs)$ and
$\wh F\op:\fs(\cs)\op\la\cl'(\cs)\op$ are both
good extensions, then in fact
$\fl'(\cs)\op=\fl'\big(\fs(\cs)\op\big)$,
and the metrics agree.
Moreover: the identification
$\fl'(\cs)\op=\fl'\big(\fs(\cs)\op\big)$,
as subcategories of $\ct$,
is such that the diagram below commutes
\[\xymatrix@C+20pt@R-10pt{
\fl'(\cs)\op
\ar[r]\ar@{=}[dd]
& 
\fl(\cs)\op\,\,
\ar@{^{(}->}[r]
\ar[dd]^\wr &
(\MMod\cs)\op
\ar[dd]^{\wh Y} \\
 & & \\
\fl'\big(\fs(\cs)\op\big)
\ar[r]
 &
\fl\big(\fs(\cs)\op\big)
\ar@{^{(}->}[r] &
\MMod{\fs(\cs)\op}
}\]
In this diagram, if we delete the left column, then
the commutative square that remains is
(a rotation of)
the bottom part of
the much more elaborate
diagram of Theorem~\ref{T900.31}.
\edis

\exm{E900.41}
Assume $\ct$ is a weakly approximable,
coherent triangulated category.
The metrics of \cite[Example~\ref{E20.3}]{Neeman18A},
on the subcategories $\ct^c\subset\ct$
and $\big(\ct^b_c\big)\op\subset\ct\op$,
are the ones that motivated the study---in
a very clear sense they are the examples
being generalized when we study
excellent metrics. That said:
there are undoubtedly many other examples of
excellent metrics out there, but we focus
only on the one that will be important
in the upcoming applications.

Let $\ct$ be a weakly approximable triangulated
category, but now drop the hypothesis that
$\ct$ is coherent.
Choose any \tstr\ $\tst\ct$ in the preferred
equivalence class, and
let $\cs=\big(\ct^b\big)\op$ be the
bounded part of the induced \tstr\
on $\ct\op$. Then with
$\cm_n=\big(\ct^{\leq-n}\cap\ct^b\big)\op$,
the pair
\[
\cs,\qquad\{\cm_n,\,n\in\nn\}
\]
is a triangulated category with an excellent
metric. Let $G\in\ct$ be
a compact generator.
Then the category $\fs(\cs)\op$, which
we will christen $\tsb$, is given by
the formula
\[
\tsb\eq\bigcup_{n\in\nn}\ogenu G{}{-n,n}\ .
\]
And the induced excellent
metric $\{\cn\op_n,\,n\in\nn\}$,
on the category $\tsb=\fs(\cs)\op$, is given
by the formula
\[
\cn\op_n\eq\ct^{\leq-n}\cap\tsb\ .
\]
All of this will be proved in Section~\ref{S96}.
\eexm

\exm{E900.43}
Let us return to earth, with a more concrete
version of the above. Let $R$ be
a ring, and we specialize Example~\ref{E900.41}
to the case where $\ct=\D(\Mod R)$ is the unbounded
derived category of all complexes of $R$--modules.

In this case the category $\tsb$ turns out to
be $\tsb=\K^b(\Proj R)$, the homotopy
category of all bounded complexes of
projective modules. And the excellent
metric on $\tsb$ is given by $\cm_i=K^b(\Proj R)^{\leq-i}$,
the bounded complexes of projective $R$-modules
that vanish in degree $>-i$.

Note that this is always an excellent metric,
no hypothesis on the ring $R$ is necessary. We can look
at the restriction of the metric to
$\K^b(\proj R)\subset\K^b(\Proj R)$, the subcategory of bounded
complexes of \emph{finitely generated}
projective modules. This is always a good metric,
but its excellence is equivalent to the category
$\D(\Mod R)$ being a coherent triangulated category.
If $R$ is a left-coherent ring then $\D(\Mod R)$ is a
coherent triangulated category,
but not vice versa.
See \cite[Example~\ref{E29098793.567}]{Neeman18A}.
\eexm

\rmk{R900.45}
In the case where $\ct=\D(\Mod R)$, the subcategory
$\tsb\subset\ct$ turned out to be a familiar,
old friend. If $X$ is a quasicompact
quasiseparated scheme, and $Z\subset X$
is a closed subset with quasicompact
complement, then the category $\ct=\Dqcs Z(X)$
is weakly approximable by
\cite[Theorem~3.2(iv)]{Neeman22A}, but
as far as I know the subcategory
$\tsb$ has not featured in the
literature. Similarly: when $\ct$ is the
homotopy category of spectra, then
\cite[Example~4.2]{Neeman24}
shows that $\ct$ is weakly approximable. And
once again the subcategory $\tsb$ seems
to be largely unstudied.
\ermk

\nin
{\bf Acknowledgements.}\ \ The author is grateful to 
Paul Balmer, Alberto Canonaco, Patrick Lank and Paolo Stellari  for improvements on an earlier version.

\section{The metrics on $\fl(\cs)$ and $\fs(\cs)$}
\label{S3}

\rmd{R3.-111222333}
Let $\cs$ be a triangulated category,
and assume $\{\cm_i,\,i\in\nn\}$
is a good metric on $\cs$.
In \cite[Definition~\ref{D20.11}(i)]{Neeman18A}
we met the category
$\fl(\cs)$, it is the completion
of $\cs$ with respect to the metric,
and in
\cite[Definition~\ref{D28.907}]{Neeman18A}
we introduced strong triangles
in $\fl(\cs)$, those are the colimits
of Cauchy sequences of distinguished triangles.
And then \cite[Remark~\ref{R28.905}]{Neeman18A}
showed that strong triangles in $\fl(\cs)$ are
pre-triangles in the sense of  \cite[Definition~\ref{D28.7}]{Neeman18A}.
Recall: this means that,
for any object $s\in\cs$,
the functor $\Hom_{\MMod\cs}^{}\big(Y(s),-\big)$
takes strong triangles to long exact sequences.
\ermd

Since it is always helpful to have ways to produce
long exact sequences, we begin the
section by proving the couple of lemmas below, each of
which produces more long exact sequences.
The lemmas will be used starting
in the next section, but are included here
because they work without any
additional hypothesis on the good metric.
And, starting from the next section, our
good metrics will no longer be general.

\lem{L3.305}
Let $\cs$ be a triangulated category, and let
$\{\cm_i\mid i\in\nn\}$ be a good metric on
$\cs$. If $A\la B\la C\la\T A$ is a
distinguished triangle in $\fs(\cs)$ and
$F\in\fl(\cs)$ is arbitrary, then $\Hom(F,-)$
takes $A\la B\la C\la\T A$ to a long exact sequence.
\elem

\prf
It suffices to prove the exactness
at a single point, that is we
prove the exactness
of
\[\xymatrix{
\Hom(F,A)\ar[r] &
\Hom(F,B)\ar[r] &
\Hom(F,C)\ .
}\]
Because $A,B,C$ all belong to $\fs(\cs)$,
we may choose an integer $n>0$ with
$A,B,C$ all in $\cm_n^\perp$.
Let $f_*$ be a
Cauchy sequence in $\cs$ with
$F=\colim\, Y(f_*)$.
Replacing $f_*$ by a subsequence, we
may assume that in every triangle $f_i\la f_{i+1}\la x_i$
we have $\Tm x_i,x_i\in\cm_{n}$. 
Hence in the
diagram below
\[\xymatrix{
\Hom(f_{i+1}^{},A)\ar[r]\ar[d] &
\Hom(f_{i+1}^{},B)\ar[r]\ar[d] &
\Hom(f_{i+1}^{},C)\ar[d]\\
\Hom(f_{i}^{},A)\ar[r] &
\Hom(f_{i}^{},B)\ar[r] &
\Hom(f_{i}^{},C)
}\]
the rows are exact and the vertical maps
are isomorphisms, and the Lemma follows just
by taking inverse limits over $i\in\nn$.
\eprf

\lem{L3.303}
Let $\cs$ be a triangulated category, and let
$\{\cm_i\mid i\in\nn\}$ be a good metric on
$\cs$. Let $F\in\fs(\cs)$ be
any object. Then the functor
$\Hom(-,F)$ takes any strong triangle
$A\la B\la C\la\T A$ in the category
$\fl(\cs)$ to a long exact sequence.
\elem

\prf
It suffices to prove the exactness
at a single point, that is we
prove the exactness
of
\[\xymatrix{
\Hom(C,F)\ar[r] &
\Hom(B,F)\ar[r] &
\Hom(A,F)\ .
}\]
Let $a_*\la b_*\la c_*\la\T a_*$ be a
Cauchy sequence of triangles in $\cs$ with
colimit $A\la B\la C\la\T A$.
As $F$ belongs to $\fs(\cs)\subset\fl(\cs)$,
\cite[Remark~\ref{R20.13}]{Neeman18A}
tells us that $\Hom(-,F)$ is a homological
functor on $\cs$. 
As $F$ belongs to $\fs(\cs)\subset\fc(\cs)$ there
must exist an integer $n>0$ with $F\in\cm_n^\perp$.
And because the sequence $a_*\la b_*\la c_*$ is
Cauchy, we may (after replacing by a subsequence)
assume that, for all $i\in\nn$, the triangles
\[
a_i^{}\la a_{i+1}^{}\la d_{i}^{}\,,\qquad
b_i^{}\la b_{i+1}^{}\la \ov d_{i}^{}\,,\qquad
c_i^{}\la c_{i+1}^{}\la \wh d_{i}^{}
\]
all have $\Tm d_i^{}$, $d_i^{}$, $\Tm \ov d_i^{}$, $\ov d_i^{}$, $\Tm \wh d_i^{}$ and $\wh d_i^{}$
in $\cm_{n}$. Hence in the
diagram below
\[\xymatrix{
\Hom(c_{i+1}^{},F)\ar[r]\ar[d] &
\Hom(b_{i+1}^{},F)\ar[r]\ar[d] &
\Hom(a_{i+1}^{},F)\ar[d]\\
\Hom(c_{i}^{},F)\ar[r] &
\Hom(b_{i}^{},F)\ar[r] &
\Hom(a_{i}^{},F)
}\]
the rows are exact and the vertical maps
are isomorphisms, and the Lemma follows just
by taking inverse limits over $i\in\nn$.
\eprf

So far we haven't introduced any new concepts or
constructions. This is about to change.

\dfn{D3.-3}
Let $\cs$ be a triangulated category with a good
metric $\{\cm_i\mid i\in\nn\}$.
\be
\item
In the
category $\fl(\cs)$ of \cite[Definition~\ref{D20.11}(i)]{Neeman18A},
define the sequence of subcategories
$\{\cl_i\subset\fl(\cs)\mid i\in\nn\}$
by the formula
\[
\cl_i\eq\left\{F\in\fl(\cs)\left|
\begin{array}{c}
  \text{\rm there exists in $\cs$ a Cauchy sequence $E_1\sr E_2\sr E_3\sr\cdots$}\\
\text{\rm with $E_n\in\cm_i$ for all $i$, and with 
$F=\colim\,Y(E_n)$} 
\end{array}
\right.\right\}\ .
\]
\item
In the category $\fs(\cs)$, define
the sequence of subcategories
$\{\cn_i\subset\fs(\cs)\mid i\in\nn\}$
by the formula
$\cn_i=\cl_i\cap\fs(\cs)$.
\ee
\edfn

We begin with the trivial little lemma.

\lem{L3.-1}
Let $\cs$ be a triangulated category with a good
metric $\{\cm_i\mid i\in\nn\}$.
Then the categories
$\{\cl_i\subset\cl(\cs)\mid i\in\nn\}$
of Definition~\ref{D3.-3}(i)
satisfy
\be
\item
$\Tm\cl_{i+1}\cup\cl_{i+1}\cup\T\cl_{i+1}\subset\cl_i$.
\item
$0\in\cl_i$ for all $i\in\nn$.
\ee
\elem

\prf
Because $0\in\cm_i$ for all $i\in\nn$, the
Cauchy sequence $0\sr 0\sr 0\sr\cdots$ is
contained in $\cm_i$. Hence its colimit
$0$ is contained in $\cl_i$.

Next we prove that
$\Tm\cl_{i+1}\cup\cl_{i+1}\cup\T\cl_{i+1}\subset\cl_i$.
Suppose therefore that $F\in\cl_{i+1}$. Then there
exists a Cauchy sequence
$E_1\sr E_2\sr E_3\sr\cdots$, with $E_n\in\cm_{i+1}$ for
all $n\in\nn$, and such that $F=\colim\, Y(E_n)$.
Then the inclusions
$\Tm\cm_{i+1}\cup\cm_{i+1}\cup\T\cm_{i+1}\subset\cm_i$
of \cite[Definition~\ref{D20.1}(i)]{Neeman18A}
tell us that
\[
\Tm F=\colim\, Y(\Tm E_n),\qquad
F=\colim\, Y(E_n)\qquad\text{and}\qquad
\T F=\colim\, Y(\T E_n)
\]
all belong to $\cl_i$.
\eprf

Slightly less trivial is the next lemma.

\lem{L3.1}
Let $\cs$ be a triangulated category with a good
metric $\{\cm_i\mid i\in\nn\}$.
Then
the $\cn_i\subset\fs(\cs)$ of
Definition~\ref{D3.-3}(ii) define
a  good metric
$\{\cn_i\mid i\in\nn\}$
on the triangulated category
$\fs(\cs)$ of \cite[Theorem~\ref{T28.128}]{Neeman18A}.
\elem

\prf
Lemma~\ref{L3.-1} teaches us that $0\in\cl_i$ for all
$i\in\nn$, and that
$\Tm\cl_{i+1}\cup\cl_{i+1}\cup\T\cl_{i+1}\subset\cl_i$.
Intersecting with $\fc(\cs)$ we have that
$0\in\cn_i$ for all $i\in\nn$, and that
$\Tm\cn_{i+1}\cup\cn_{i+1}\cup\T\cn_{i+1}\subset\cn_i$.

It remains to prove that
$\{\cn_i\mid i\in\nn\}$
satisfies \cite[Definition~\ref{D20.1}(ii)]{Neeman18A}, that is
in the triangulated category $\fs(\cs)$ we
have
$\cn_i*\cn_i\subset\cn_i$.
Suppose therefore that $A\la B\la C\stackrel f\la \T A$ is
a distinguished
triangle in $\fs(\cs)$ with $A,C\in\cn_i$. Then
there must exist Cauchy sequences $a_*$ and $c_*$
in $\cm_i$, such that $A=\colim\,Y(a_*)$ and
$C=\colim\,Y(c_*)$. Now 
\cite[Lemma~\ref{L28.19}]{Neeman18A} tells us that we
may choose a subsequence $a'_*$ of $a_*$ and a map
of sequences $f_*:c_*\la\T a'_*$ such
that $f=\colim\,Y(f_*)$, and then continue in $\cs$
to a Cauchy sequence of triangles
$a'_*\la b_*\la c_*\stackrel{f_*}\la\T a'_*$.
Because $b_n\in\cm_i*\cm_i\subset\cm_i$ for all $n\in\nn$,
we have that $\ov B=\colim\,Y(b_*)$ belongs to
$\cn_i$. But now the colimit of
the Cauchy sequence of distinguished triangles
in $\cs$ is a strong triangle
$A\la \ov B\la C\stackrel f\la\T A$ with
$A,C\in\fs(\cs)$, 
and \cite[Remark~\ref{R200976}]{Neeman18A} tells
us that it is a distinguished triangle
in $\fs(\cs)$. Therefore
it must be isomorphic (non-canonically) to the
given distinguished triangle
$A\la B\la C\stackrel f\la \T A$,
and we deduce that $B\cong\ov B$ must belong
to $\cn_i$.
\eprf

Also useful will be the
orthogonals of the $\cm_i$, $\cl_i$
and $\cn_i$. The next lemma studies
some elementary relations among them.

\lem{L3.2.5}
Let $\cs$ be a triangulated category, and let
$\{\cm_i\mid i\in\nn\}$ be a good metric on
$\cs$. Then
\be
\item
${^\perp Y(\cm_i)}\cap Y(\cs)={^\perp\cl_i}\cap Y(\cs)$.
\item
In the category $\MMod\cs$ we have   
$Y(\cm_i)^\perp=\cl_i^\perp$, and hence
\[
Y(\cs)\cap Y(\cm_i)^\perp= Y(\cs)\cap \cl_i^\perp\qquad
\fs(\cs)\cap Y(\cm_i)^\perp= \fs(\cs)\cap \cl_i^\perp
\]
\ee
\elem

\prf
Note that, as $Y(\cm_i)\subset\cl_i$, we have that
$\cl_i^\perp\subset Y(\cm_i)^\perp$ and $^\perp\cl_i\subset {^\perp Y(\cm_i)}$.
It's the reverse inclusions that need proof.

To prove (i) note that, for any object $s\in\cs$, the object
$Y(s)$ is finitely presented in the category $\MMod\cs$,
meaning that $\Hom_{\MMod\cs}^{}\big(Y(s),-\big)$
respects filtered colimits. If $b_*$ is a Cauchy sequence in
$\cm_i$ and $Y(s)\in{^\perp Y(\cm_i)}\cap Y(\cs)$, then
\[
\Hom\big(Y(s)\,,\,\colim\, Y(b_*)\big)\eq
\colim\,\big(Y(s)\,,\,Y(b_*)\big)\eq 0\ .
\]
Hence $Y(s)\in{^\perp\cl_i}\cap Y(\cs)$.

Now for (ii).
Clearly $Y(\cm_i)\subset{^\perp\big(Y(\cm_i)^\perp)}$,
but as ${^\perp\big(Y(\cm_i)^\perp)}$ is closed under colimits
we deduce that $\cl_i\subset{^\perp\big(Y(\cm_i)^\perp)}$. Hence
$\Hom\big(\cl_i\,,\,Y(\cm_i)^\perp\big)=0$ and
$Y(\cm_i)^\perp\subset\cl_i^\perp$.
\eprf

We will frequently use the following, immediate consequence

\rmk{R3.2??}
With $\fc(\cs)$ as in \cite[Definition~\ref{D20.11}(ii)]{Neeman18A},
we have that
\cite[Observation~\ref{O28.-1}]{Neeman18A} gives the first equality in
\[
\fc(\cs)\eq\bigcup_{i\in\nn} Y(\cm_i)^\perp\eq
\bigcup_{i\in\nn} \cl_i^\perp\ ,
\]
while the second comes from Lemma~\ref{L3.2.5}(ii).

Pairing the formula
above, with
the formula
for the subcategory $\fs(\cs)\subset\fc(\cs)$
in
\cite[Definition~\ref{D20.11}(iii)]{Neeman18A},
the combination yields the formula
\[
\fs(\cs)\eq\fl\cs)\cap\fc(\cs)\eq\fl(\cs)\cap
\bigcup_{i\in\nn} \cl_i^\perp\ .
\]
\ermk

\dis{D3.whoknows}
We will be interested in constructing,
in the category $\fl(\cs)$,
strong triangles $A\la B\la C\la\T A$,
where for some integer $i>0$ the object $C$
belongs to
$\cl_i\subset\fl(\cs)$.
More precisely: we will ask ourselves
what hypothesis one
needs, to place on a morphism
$F:A\la B$ in the category $\fl(\cs)$,
to guarantee that
it can be completed to such a strong
triangle.

Of course: we can choose Cauchy
sequences $a_*,b_*$ in $\cs$
such that $A=\colim\,Y(a_*)$
and $B=\colim\,Y(b_*)$,
and then (up to replacing
$b_*$ by a subsequence) we can
choose a morphism of
sequences $f_*:a_*\la b_*$
in such a way that $F=\colim\, Y(f_*)$.
And then we can complete to a Cauchy
sequence of triangles, obtaining
as a colimit a strong triangle
$A\la B\la \wt C\la\T A$.
Of course: this strong
triangle depends on all the
choices we have made.
Now suppose that $\wt C\in\cl_i$ for
some choice of Cauchy sequences.
Is it in $\cl_i$ for every other choice?

I 
have no idea. Without
solving this problem, we
will introduce
a method that
gets around this difficulty.
The idea is to consider
not just the morphism
$F:A\la B$ in $\fl(\cs)$,
but also a pair of
Cauchy sequences
$a_*,b_*$ in $\cs$
such that $A=\colim\,Y(a_*)$
and $B=\colim\,Y(b_*)$.
With this much data chosen in
advance, the remaining choices
turn out not to matter.
We will now set out to make
this precise.
\edis

In Discussion~\ref{D3.whoknows}
we raced ahead of ourselves,
for now we make the following definition.

\dfn{D3.333}
Let $\cs$ be a triangulated category
with a good metric
$\{\cm_i\mid i\in\nn\}$. If
$f_*:a_*\la b_*$ is a map of
Cauchy sequences in $\cs$,
then $f_*$ is declared to be a
\emph{type-$n$ morphism} if, in any
completion to
a Cauchy sequence of distinguished
triangles
$a_*\stackrel{f_*}\la b_*\stackrel{g_*}\la c_*\stackrel{h_*}\la \T a_*$
in $\cs$,
we have that $c_i\in\cm_n$ for all
$i\gg0$.
\edfn

\rmk{D3.335}
The reader should note that, for each $i>0$, the
object $c_i$ is determined up to
isomorphism by the morphism
$f_i:a_i\la b_i$.
Up to (non-canonical) isomorphism there
is only one distinguished triangle
$a_i\stackrel{f_i}\la b_i\stackrel{g_i}\la c_i\stackrel{h_i}\la \T a_i$
completing the map $f_i:a_i\la b_i$.
The subtlety in the proof of
\cite[Lemma~\ref{L28.19}(iii)]{Neeman18A} was to show that the
morphisms of triangles can be chosen
to render $c_*$ into a Cauchy sequence.
\ermk

\lem{L3.337}
Let $\cs$ be a triangulated category
with a good metric
$\{\cm_i\mid i\in\nn\}$.
If a map of Cauchy sequences $f'_*:a'_*\la b'_*$
is a type-$n$ morphism then so is any
subsequence $f_*:a_*\la b_*$. Less
trivially, the
strong reverse direction is also true:
if there exists a subsequence $f_*:a_*\la b_*$
which is a type-$n$ morphism,
then $f'_*:a'_*\la b'_*$
is also a type-$n$ morphism.
\elem

\prf
The first implication is trivial: 
if $f'_*:a'_*\la b'_*$
is a type-$n$ morphism then so is
any subsequence $f_*:a_*\la b_*$.

Now assume there exists a subsequence
$f_*:a_*\la b_*$ which is a type-$n$ morphism.
Then
\cite[Lemma~\ref{L28.19}(iii)]{Neeman18A} permits us to
complete $f'_*:a'_*\la b'_*$ to a Cauchy
sequence of triangles
$a'_*\stackrel{f'_*}\la b'_*\stackrel{g'_*}\la c'_*\stackrel{h'_*}\la \T a'_*$.
Because $c'_*$ is Cauchy, there exists an
integer $N>0$ such that, if $m<m'$ is a pair
of integers both $>N$, then the triangle
$c'_m\la c'_{m'}\la d_{m,m'}$ has $d_{m,m'}\in\cm_n$.
But because $f_*:a_*\la b_*$ is a type-$n$ morphism
we may choose an integer $m>N$ such that, in the
triangle
$a_m\stackrel{f_m}\la b_m\stackrel{g_m}\la c_m\stackrel{h_m}\la \T a_m$,
we have $c_m\in\cm_n$. 
On the other hand $f_*:a_*\la b_*$
is a subsequence of $f'_*:a'_*\la b'_*$, 
and hence $f_m:a_m\la b_m$ is
equal to $f'_{\ell+m}:a'_{\ell+m}\la b'_{\ell+m}$
for some $\ell\geq0$, and we deduce that $c'_{\ell+m}\cong c_m$
must also belong to $\cm_n$. But now,
for any $i>\ell+m>N$, the triangle
$c'_{\ell+m}\la c'_i\la d_{\ell+m,i}$ exhibits
$c'_i$ as belonging to $\cm_n*\cm_n\subset\cm_n$.
\eprf

\dfn{D3.343}
Let $\cs$ be a triangulated category
with a good metric
$\{\cm_i\mid i\in\nn\}$. 
Let $A$ and $B$ be objects of $\fl(\cs)$,
let $F:A\la B$ be a morphism in $\fl(\cs)$, and
let $a'_*$ and $b'_*$ be Cauchy
sequences in $\cs$ with $A=\colim\,Y(a'_*)$ and
$B=\colim\,Y(b'_*)$.

Then the set $L(F,a'_*,b'_*)$ is defined by
\[
L(F,a'_*,b'_*)\eq
\left\{f_*:a_*\la b_*\left|
\begin{array}{c}
\text{\rm $a_*$ is a subsequence of $a'_*$,}\\
\text{\rm $b_*$ is a subsequence of $b'_*$,}\\
\text{\rm and }F=\colim\,Y(f_*)
\end{array}
\right.\right\}
\]
The set $L(F,a'_*,b'_*)$ is non-empty by
\cite[Lemma~\ref{L28.19}(i)]{Neeman18A}.
\edfn

\lem{L3.339}
Let $\cs$ be a triangulated category
with a good metric
$\{\cm_i\mid i\in\nn\}$. 
Let $a'_*$ and $b'_*$ be Cauchy sequences in $\cs$,
let $A=\colim\,Y(a'_*)$ and $B=\colim\,Y(b'_*)$,
and let $F:A\la B$ be a morphism
in $\fl(\cs)$.

If one element $f_*:a_*\la b_*$, of the set
$L(F,a'_*,b'_*)$ of
Definition~\ref{D3.343}, is a
type-$n$ morphism, then all
the elements of $L(F,a'_*,b'_*)$
are type-$n$ morphisms.
\elem

\prf
Suppose $f_*:a_*\la b_*$
and $\wt f_*:\wt a_*\la \wt b_*$
are two elements of the set
$L(F,a'_*,b'_*)$,
with $f_*:a_*\la b_*$ a type-$n$ morphism.
We need to show that
$\wt f_*:\wt a_*\la \wt b_*$ is also
a type-$n$ morphism.

We are given a Cauchy sequence $a'_*$ with
two subsequences $a_*$ and $\wt a_*$,
and a Cauchy sequence $b'_*$ with
two subsequences $b_*$ and $\wt b_*$.
And now we propose to construct
a third pair of Cauchy sequences $\wh a_*$ and
$\wh b_*$, and a map $\wh f_*:\wh a_*\la\wh b_*$
with $F=\colim\,Y(\wh f_*)$.

To begin with: the map
$\wh f_1:\wh a_1^{}\la \wh b_1$
is declared to be
$f_1:a_1^{}\la b_1$.
And now we proceed by induction, in such a way that
\be
\item
$\wh a_*$ is a subsequence of $a'_*$
and $\wh b_*$ is a subsequence of $b'_*$.
\item
The odd subsequence
$\wh f_{2{\ds *}-1}: \wh a_{2{\ds *}-1}^{}\la\wh b_{2{\ds *}-1}$
is some subsequence of
$f_*:a_*\la b_*$.
\item
The even subsequence
$\wh f_{2{\ds *}}: \wh a_{2{\ds *}}^{}\la\wh b_{2{\ds *}}$
is some subsequence of
$\wt f_*:\wt a_*\la \wt b_*$
\ee
We have constructed the map
$\wh f_1:\wh a_1^{}\la \wh b_1$
and need to prove the induction step.
As the odd and even cases are very similar
we treat only the case where we have
constructed $\wh f_{2m-1}: \wh a_{2m-1}^{}\la\wh b_{2m-1}$,
which agrees with 
$f_n:a_n\la b_n$ for some
$n$ large enough, and want to construct
$\wh f_{2m}: \wh a_{2m}^{}\la\wh b_{2m}$.

Now because $\wh a_{2m-1}^{}=a_n=a'_{\ell+n}$ and
$\wh b_{2m-1}^{}=b_n=b'_{\ell'+n}$ for some
$\ell,\ell'\geq0$, then for any integer
$N>\max(\ell+n,\ell'+n)$
we have that $\wt a_N^{}=a'_{N+r}$
and $\wt b_N^{}=b'_{N+r'}$
and we may form a diagram
\[\xymatrix@C+30pt{
\wh a_{2m-1}^{}\ar[d]_\alpha\ar[r]^-{\wh f_{2m-1}}
  &\wh b_{2m-1}\ar[d]^\beta\\
\wt a_N^{} \ar[r]^-{\wt f_N^{}}\ar[d] &
    \wt b_N^{}\ar[d] \\
A\ar[r]^-F & B
}\]
where the map $\alpha:\wh a_{2m-1}^{}\la \wt a_N^{}$
is just the morphism $a'_{\ell+n}\la a'_{N+r}$
of the Cauchy sequence $a'_*$, and the map
$\beta:\wh b_{2m-1}^{}\la \wt b_N^{}$
is just the morphism $b'_{\ell+n}\la b'_{N+r'}$
of the Cauchy sequence $b'_*$.
Of course: there is no reason for the top square
to commute, we are not given any compatibility
between the maps
$f_*:a_*\la b_*$
and
$\wt f_*:\wt a_*\la \wt b_*$
beyond the hypothesis that
$\colim\,Y(f_*)=F=\colim\,Y(\wt f_*)$.
But if we delete the top row
then the bottom square
commutes, and if we delete the middle row then
the large square commutes. We deduce the
equality of the composites from top left to
bottom right in
\[\xymatrix@C+30pt{
\wh a_{2m-1}^{}\ar[d]_\alpha\ar[r]^-{\wh f_{2m-1}}
  &\wh b_{2m-1}\ar[d]^\beta\\
\wt a_N^{} \ar[r]^-{\wt f_N^{}} &
\wt b_N^{}\ar[d] \\
& B
}\]
Therefore the map $\wt b_N^{}\la B=\colim\,\wt b_*$
coequalizes the two composites in the top square, and
there must exists an integer $M>N$ such that
the two composites from top right to bottom
left in
\[\xymatrix@C+30pt{
\wh a_{2m-1}^{}\ar[d]_{\alpha}\ar[r]^-{\wh f_{2m-1}}
  &\wh b_{2m-1}\ar[d]^{\beta}\\
\wt a_N^{} \ar[r]^-{\wt f_N^{}} &
\wt b_N^{}\ar[d] \\
& \wt b_M^{}
}\]
agree. But then the square
\[\xymatrix@C+30pt{
\wh a_{2m-1}^{}\ar[d]_{\alpha'}\ar[r]^-{\wh f_{2m-1}}
  &\wh b_{2m-1}\ar[d]^{\beta'}\\
\wt a_M^{} \ar[r]^-{\wt f_M^{}} &
\wt b_M^{} 
}\]
commutes, permitting us to define 
$\wh f_{2m}: \wh a_{2m}^{}\la\wh b_{2m}$ to be
$\wt f_{M}: \wt a_{M}^{}\la\wt b_{M}$.

And now to finish the proof: we are given that
$f_*:a_*\la b_*$ is a type-$n$ morphism of
Cauchy sequences. Hence so is the
subsequence
$\wh f_{2{\ds *}-1}: \wh a_{2{\ds *}-1}^{}\la\wh b_{2{\ds *}-1}$.
But Lemma~\ref{L3.337} tells us that if a
subsequence is a type-$n$ morphism
then the entire sequence is,
and we deduce that
$\wh f_*:\wh a_*\la\wh  b_*$
must be a type-$n$ morphism, as is its
even subsequence
$\wh f_{2{\ds *}}: \wh a_{2{\ds *}}^{}\la\wh b_{2{\ds *}}$.
But this is a subsequence of 
$\wt f_*:\wt a_*\la\wt  b_*$,
and by using Lemma~\ref{L3.337} again
we conclude that 
$\wt f_*:\wt a_*\la\wt  b_*$
must be a type-$n$ morphism.
\eprf

\dfn{D3.341}
Let $\cs$ be a triangulated category
with a good metric
$\{\cm_i\mid i\in\nn\}$. 
Let $A$ and $B$ be objects of $\fl(\cs)$, and
let $a_*$ and $b_*$ be Cauchy
sequences in $\cs$ with $A=\colim\,Y(a_*)$ and
$B=\colim\,Y(b_*)$.

A morphism $F:A\la B$, in the category $\fl(\cs)$,
is defined to be \emph{of type-$n$ with respect to
$(a_*,b_*)$} if some
element $f'_*:a'_*\la b'_*$
of the non-empty set $L(F,a_*,b_*)$ of
Definition~\ref{D3.343} is a type-$n$ morphism.
\edfn

By Lemma~\ref{L3.339}
this
implies that every element of
$L(F,a_*,b_*)$ is
a type-$n$ morphism.

\lem{L3.341}
Let $\cs$ be a triangulated category
with a good metric
$\{\cm_i\mid i\in\nn\}$. 
Let $A$, $B$ and $C$ be objects of $\fl(\cs)$
and let $a_*$, $b_*$ and $c_*$ be Cauchy
sequences in $\cs$ with $A=\colim\,Y(a_*)$,
$B=\colim\,Y(b_*)$ and $C=\colim\,Y(c_*)$.

Now suppose that $A\stackrel F\la B\stackrel G\la C$
are two composable morphisms in $\fl(\cs)$.
Then
\be
\item
If $F$ is a type-$n$ morphism with respect
to $(a_*,b_*)$ and $G$ is a type-$n$ morphism
with respect to $(b_*,c_*)$, then
$GF$ is a type-$n$ morphism with respect
to $(a_*,c_*)$.
\item
If 
$GF$ is a type-$n$ morphism with respect
to $(a_*,c_*)$, and $F$ is a type-$(n+1)$
morphism with respect
to $(a_*,b_*)$ [repsectively: $G$ is a type-$(n+1)$ morphism
  with respect to $(b_*,c_*)$], then
$G$ is a type-$n$ morphism
with respect to $(b_*,c_*)$
[respectively: $F$ is a type-$n$ morphism with respect
to $(a_*,b_*)$].
\ee
\elem

\prf
By \cite[Lemma~\ref{L28.19}(i)]{Neeman18A} we may, up to replacing
$b_*$ by a subsequence, choose a map
$f_*:a_*\la b_*$ with $F=\colim\,Y(f_*)$. And
applying \cite[Lemma~\ref{L28.19}(i)]{Neeman18A} again we may,
this time replacing $c_*$ by a subsequence,
choose a map $g_*:b_*\la c_*$ with
$G=\colim\,Y(g_*)$. And by
Lemma~\ref{L3.339} the type of
the morphisms $f_*$, $g_*$ and $g_*f_*$ is independent
of the choices.

For each integer $\ell>0$ we may complete the composable
morphisms $a_\ell^{}\la b_\ell^{}\la c_\ell^{}$
to an octahedron
\[\xymatrix@C+20pt{
a_\ell^{}\ar@{=}[d]\ar[r]^-{f_\ell^{}} &
   b_\ell^{}\ar[d]^-{g_\ell^{}}\ar[r] &
       d_\ell^{}\ar[d]\\
a_\ell^{}\ar[r]^-{g_\ell^{}f_\ell^{}} &
  c_\ell^{}\ar[r]\ar[d] &
    \wh d_\ell^{}\ar[d]\\
 &\ov d_\ell^{}\ar@{=}[r] &\ov d_\ell^{}
}\]
If $\ell\gg0$ then our assumption in (i) is
that $d_\ell^{}$ and $\ov d_\ell^{}$ lie in
$\cm_n$, and the triangle
$d_\ell^{}\la \wh d_\ell^{}\la\ov d_\ell^{}$
allows us to conclude
that $\wh d_\ell^{}\in\cm_n*\cm_n\subset\cm_n$.
This makes $GF$ a type-$n$ morphism with respect to
$(a_*,c_*)$. 

In (ii) our assumptions are that $\wh d_\ell^{}$
lies in $\cm_n$, and that one of $d_\ell^{}$ or
$\ov d_\ell^{}$ belongs to $\cm_{n+1}$,
depending on the respective cases.
Assume $F$ is a type-$(n+1)$ morphism
with respect to $(a_*,b_*)$; then using
the triangle
$\wh d_\ell^{}\la\ov d_\ell^{}\la\T d_\ell^{}$
coupled with the fact that
$\T d_\ell^{}\in\T\cm_{n+1}\subset\cm_{n}$,
we conclude that $\ov d_\ell^{}\in\cm_n*\cm_n\subset\cm_n$
making $G$ a type-$n$ morphism
with respect to $(b_*,c_*)$.

Assume $G$ is a type-$(n+1)$ morphism
with respect to $(b_*,c_*)$; then using
the triangle
$\Tm\ov d_\ell^{}\la d_\ell^{}\la \wh d_\ell^{}$,
coupled with the fact that
$\Tm\ov d_\ell^{}\in\Tm\cm_{n+1}\subset\cm_{n}$,
we conclude that $d_\ell^{}\in\cm_n*\cm_n\subset\cm_n$,
making $F$ a type-$n$ morphism
with respect to $(a_*,b_*)$.
\eprf

It is natural enough to wonder if being a type-$n$
morphism depends on the choice of Cauchy sequences. To spell
it out: Definition~\ref{D3.341} tells us what it
means, for a morphism $F:A\la B$ in
the category $\fl(\cs)$, to
be of type-$n$ \emph{with respect to a pair of
Cauchy sequences}
$(a_*,b_*)$. Suppose we change the Cauchy sequences.
Can the type change?

\exm{E3.977}
The answer is Yes, at least for sufficiently ridiculous
good metrics. The author would like to thank Kabeer
Manali-Rahul for the example below.

Let $\cs$ be an
essentially small triangulated category, and
assume its Grothendieck group $K_0(\cs)$ is
non-zero. Let the good metric $\{\cm_i,\,i\in\nn\}$
be defined by
\[
\cm_i=\{X\in\cs\mid \text{In the Grothendieck group }K_0(\cs)\text{, we have }[X]=0\}\ .
\]
Note that, by definition, $\cm_i$ is independent of $i\in\nn$
and is a triangulated subcategory of $\cs$.

Choose any object $X\in\cs$ with $[X]\in K_0(\cs)$
nonzero.
Then the following two sequences $a_*$
and $b_*$ below are both Cauchy
\[\xymatrix@R-20pt{
a_*\ar@{}[r]|-{\ds:}& & 0\ar[r] &0\ar[r] &0\ar[r] &0\ar[r] &0\ar[r] &\cdots \\
b_*\ar@{}[r]|-{\ds:}& & X\ar[r]^-{0} &X\ar[r]^-{0} &X\ar[r]^-{0} &X\ar[r]^-{0} &X\ar[r] &\cdots 
}\]
and $\colim\,Y(a_*)=0=\colim\,Y(b_*)$.
But now: the identity map $\id:0\la 0$,
in the category $\cl(\cs)$, is of type-$n$ with respect to
$(a_*,a_*)$ for any $n\geq1$, but is not of
type-$n$ with respect to $(a_*,b_*)$ for any
$n\geq1$.
\eexm

There is a sense in which the pathology of
Example~\ref{E3.977} is the only problem. We prove:

\lem{L3.979}
Let $\cs$ be a triangulated category with a
good metric $\{\cm_i\mid i\in\nn\}$.
Assume further that if 
\[\xymatrix@R-20pt{
c_1^{}\ar[r]^-{0} &c_2^{}\ar[r]^-{0} &c_3^{}\ar[r]^-{0} &c_4^{}\ar[r]^-{0} &c_5^{}\ar[r] &\cdots 
}\]
is a Cauchy sequence in which the maps all
vanish, then for any $n>0$ there exists an
$N\gg0$ with $c_i\in\cm_n$ for all $i>N$.

Given any object $F\in\fl(\cs)$, and
two Cauchy sequences $a_*$ and $b_*$ in $\cs$
with $\colim\,Y(a_*)=F=\colim\,Y(b_*)$,
then the identity map $\id:F\la F$ is
of type-$n$ with respect to $(a_*,b_*)$
for any $n\geq1$.
\elem

\prf
We begin by recalling that
\cite[Lemma~\ref{L28.19}(i)]{Neeman18A}
permits us to replace $a_*$ and $b_*$ by subsequences
and choose a map $\ph_*:a_*\la b_*$ such
that $\id=\colim\,Y(\ph_*)$. And
\cite[Lemma~\ref{L28.19}(ii)]{Neeman18A} permits us to extend
to a Cauchy sequence of triangles
$a_*\stackrel{\ph_*}\la b_*\la c_*\la\T a_*$.
Taking colimits gives a strong triangle
$F\stackrel\id\la F\la G\la \T F$
in the category $\fl(\cs)$, but by
\cite[Remark~\ref{R28.905}]{Neeman18A} strong triangles are
pre-triangles, making them exact sequences
in the abelian category $\MMod\cs$. Hence $G=0$, and
the Cauchy sequence $c_*$ is such that $\colim\,Y(c_*)=0$.

But then, for any integer $i>0$, the map $c_i\la\colim\,Y(c_*)$
vanishes. This means that there must exist an integer
$j>i$ such that the natural map $c_i\la c_j$ vanishes.
Thus, up to replacing $c_*$ by a subsequence, we may assume
it to be the sequence
\[\xymatrix@R-20pt{
c_1^{}\ar[r]^-{0} &c_2^{}\ar[r]^-{0} &c_3^{}\ar[r]^-{0} &c_4^{}\ar[r]^-{0} &c_5^{}\ar[r] &\cdots 
}\]
For any $n>0$,
by assumption all but finitely many $c_i$ belongs
to $\cm_n$. This makes the identity
map $\id:F\la F$ a type-$n$ morphism with
respect to $(a_*,b_*)$.
\eprf

\cor{C3.981}
Keeping the assumptions as in Lemma~\ref{L3.979},
we have that the
type of a morphism in $\fl(\cs)$ is
independent of any choice of Cauchy sequences.
Precisely: Let $F:A\la B$ be a morphism in $\fl(\cs)$.
Assume there exist Cauchy sequences $a_*$ and $b_*$
in $\cs$ with $A=\colim\,Y(a_*)$ and
$B=\colim\,Y(b_*)$. If $F$ is a type-$n$ morphism
with respect to $(a_*,b_*)$ then it is a
type-$n$ morphism with respect to any pair of
Cauchy sequences $a'_*,b'_*$ with $A=\colim\,Y(a'_*)$
and with $B=\colim\,Y(b'_*)$.
\ecor

\prf
Just consider the composite
\[\xymatrix@R-20pt{
A\ar[r]^-{\id}\ar@{=}[d]&
A\ar[r]^-{F}\ar@{=}[d]&
B\ar[r]^-{\id}\ar@{=}[d]&
B\ar@{=}[d]\\
\colim\,Y(a'_*)&
\colim\,Y(a_*)&
\colim\,Y(b_*)&
\colim\,Y(b'_*)
}\]
and apply Lemmas~\ref{L3.341} and \ref{L3.979}.
\eprf

The next little lemma allows us to recognize morphisms
$F:A\la B$ which are type-$n$ with respect to
a pair $(a_*,b_*)$.

\lem{L3.307}
Let $\cs$ be a triangulated category with a
good metric $\{\cm_i\mid i\in\nn\}$.
Let $F:A\la B$ be a morphism in the category $\fl(\cs)$,
and choose in $\cs$ Cauchy sequences
$a_*$ and $b_*$ with
$A=\colim\,Y(a_*)$ and $B=\colim\,Y(b_*)$. Assume, for
some integer
$n>0$, that the following hold:
\be
\item
We are given some positive integers $k,\ell$
and a commutative square
\[\xymatrix{
Y(a_k^{})\ar[r]^-{Y(\ph)}\ar[d] & Y(b_\ell^{})\ar[d] \\
A\ar[r]^-\ph & B
}\]
\item
The integers $k,\ell$ are chosen large enough so
that, for all  $j\geq i\geq k$, the triangles
$a_i\la a_{j}\la d_{i,j}$ have $\T d_{i,j}\in\cm_{n}$,
and, for all
$j\geq i\geq \ell$, the triangles
$b_i\la b_{j}\la \wh d_{i,j}$ have $\wh d_{i,j}\in\cm_{n}$.
\item
If we complete $a_k^{}\la b_\ell^{}$ to a triangle
$a_k^{}\la b_\ell^{}\la c'\la\T a_k^{}$,
then we have $c'\in\cm_n$.
\ee
Then $F:A\la B$ is a type-$n$
morphism with respect to
the pair $(a_*,b_*)$.
\elem

\prf
First of all: define the
sequence $\wt a_*$ by putting $\wt a_i^{}=a_{i+k-1}^{}$, and
replace $b_*$ by the subsequence $b'_*$ defined
by $b'_i=b_{i+\ell-1}^{}$. And then we apply
\cite[Lemma~\ref{L28.19}(i)]{Neeman18A}:
starting with the given commutative square
of (i) above, which rewrites as
\[\xymatrix{
\wt a_1^{}\ar[r]^-{\ph}\ar[d] & b'_1\ar[d] \\
A\ar[r]^-\ph & B
}\]
we may choose a subsequence $\wt b_*$ of $b'_*$,
with $\wt b_1^{}=b'_1$, and a map of
sequences $f_*:\wt a_*\la\wt b_*$
such that $f_1=\ph$ and $F=\colim\,Y(f_*)$.
And then, as in
\cite[Lemma~\ref{L28.19}(ii)]{Neeman18A},
we may extend to a Cauchy sequence of triangles
$\wt a_*\stackrel{f_*}\la\wt b_*\la\wt c_*\la\T\wt a_*$
in the category $\cs$.

What needs proof is that, for all $\ell\in\nn$, we have
$\wt c_\ell^{}\in\cm_n$. For this observe that 
the commutative square
\[\xymatrix{
\wt a_1^{}\ar[r]\ar[d] & \wt b_1\ar[d]\\
\wt a_\ell^{}\ar[r] &\wt b_\ell^{}
}\]
may be completed to a diagram where the rows and
columns are triangles
\[\xymatrix{
\wt a_1^{}\ar[r]\ar[d] &\wt b_1\ar[r]\ar[d] 
  &\wt c_1^{}\ar[r]\ar[d]& \wt a_1^{}\ar[d]\\
\wt a_\ell^{}\ar[r]\ar[d] &\wt b_\ell^{}\ar[r]\ar[d] &
\wt c_\ell^{}\ar[r]\ar[d] & \T\wt a_\ell^{}\ar[d] \\
d \ar[r] & \wh d \ar[r] &
y \ar[r] &   \T d
}\]
In hypothesis (ii) of the
current Lemma we are given that
$\wh d$ and $\T d$ both belong to $\cm_{n}$,
and hence the triangle $\wh d\la y\la\T d$
in the third row exhibits
that $y\in\cm_n*\cm_n\subset\cm_n$.
But hypothesis (iii) says that $\wt c_1^{}=c'$
also lies
in $\cm_n$, and the triangle
$\wt c_1^{}\la \wt c_\ell^{}\la y$ of
the third column shows that
$\wt c_\ell^{}\in\cm_n*\cm_n\subset\cm_n$.
\eprf

\section{Excellent metrics and the equivalence
$\wh Y:\fl(\cs)\op\la\fl\big[\fs(\cs)\op\big]$}
\label{S70}

So far this article has remained in the general
setting introduced in its
predecessor, more specifically back in
\cite[Section~\ref{S28}]{Neeman18A}: we have
until now stayed in the context where
the good metrics on our triangulated categories
are unrestricted. This is about to change.

\dfn{D3.3}
Let $\cs$ be a triangulated category with a
good metric $\{\cm_i\mid i\in\nn\}$.
The metric is called \emph{excellent} if
the following holds:
\be
\item
$\cs=\cup_{i\in\nn}\,{^\perp\cm_i}$.
\item
For every integer $m\in\nn$ there
exists an integer $n>m$ such that any object
$F\in\cs$ admits, in $\cs$,
a triangle $E\la F\la D\la\T E$ with $E\in {^\perp\cm_n}$
and with $D\in\cm_m$.
\item
For every integer $m\in\nn$ there
exists an integer $n>m$ such that any object
$F\in\cs$ admits, in $\fl(\cs)$,
a type-$m$ morphism $Y(F)\la D$  with respect to
$(F,d_*)$, where in the notation $(F,d_*)$ the
Cauchy sequence $F$ is taken to be 
the constant sequence
$F\stackrel\id\la F\stackrel\id\la F\stackrel\id\la\cdots$,
and the only
restriction, on the Cauchy sequence $d_*$ and on 
the object $D=\colim\,Y(d_*)$,
is that
$D\in\fs(\cs)\cap\cl_n^\perp$.
\ee
\edfn

\rmk{R3.!?!}
There is some redundancy in the
statement of Definition~\ref{D3.3}(iii),
more specifically in the assertions
about $D$. After all: we are
given that $D=\colim\,Y(d_*)$
belongs to $\fl(\cs)$. If we
also know that $D\in\cl_n^\perp$ for some
integer $n>0$, then it automatically follows
that
\[
D\quad\in\quad\fl(\cs)\cap\cl_n^\perp\sub
\fl(\cs)\cap\bigcup_{i\in\nn}\cl_i^\perp\eq\fs(\cs)\ ,
\]
where the equality is from Remark~\ref{R3.2??}.
\ermk

The motivating example, which led to this definition,
is the following.

\exm{E3.39393901}
Let $\ct$ be a coherent, weakly approximable
triangulated category as in \cite{Neeman18A}.
Let $\cs=\ct^c$, and give it the metric
$\{\cm_i,\,i\in\nn\}$ of
\cite[Example~\ref{E20.3}(i)]{Neeman18A}.
We remind the reader: we choose
in $\ct$ a
\tstr\ $\tst\ct$ in the preferred equivalence
class, and the metric is defined by
the formula $\cm_i=\ct^c\cap\ct^{\leq-i}$.
We assert that this is an excellent metric
on $\cs$.
\eexm

\prf
We check that the conditions
in Definition~\ref{D3.3} are
satisfied. For every object
$F\in\ct^c$, the vanishing of
$\Hom(F,\ct^{\leq i})$ for some $i\gg0$ comes
from
\cite[Remark~0.24 coupled with Lemma~2.8]{Neeman24}.
That is
\[
\ct^c\sub\bigcup_{i\in\nn}^{}{^\perp\ct^{\leq i}}\sub\bigcup_{i\in\nn}^{}{^\perp\cm_i}\ ,
\]
which proves that Definition~\ref{D3.3}(i)
is satisfied.

To prove Definition~\ref{D3.3}(ii),
pick an integer $m>0$.
Let the compact generator
$G\in\ct$ and
the integers $A,B>0$ be as
in
\cite[Proposition~2.14]{Neeman24}.
Any object $F\in\ct^c\subset\ct^-$ belongs
to $\ct^{\leq C}$ for some integer $C>0$,
and \cite[Proposition~2.14]{Neeman24}
allows us to form a triangle
$E\la F\la D$ with
$D\in\ct^c\cap\ct^{\leq-m}$, and with
\[
E\quad\in\quad\genu G{}{1-m-B,B+C}
\sub{^\perp\ct^{\leq m-A-B}}\sub{^\perp\cm_{m+A+B}}\ .
\]
This proves that $n=m+A+B$ works, with $n>m$ 
satifying Definition~\ref{D3.3}(ii)

And finally we come to the proof of
Definition~\ref{D3.3}(iii), which is where
the coherence of the
triangulated category $\ct$ enters.
And for the purpose of the proof
we adopt the notation
of \cite{Neeman18A}; the inclusion
$\cs=\ct^c\la\ct$ will be denoted $F:\cs\la\ct$,
and the Yoneda functors $\cy:\ct\la\Mod\cs$ and
$Y:\cs\la\Mod\cs$ will be as in
\cite{Neeman18A}.

Let the integer $N>0$ be chosen as in
\cite[Definition~\ref{D29.1}]{Neeman18A},
which is
the definition of the coherence of $\ct$.
Given any integer $m>0$, let $B\in\ct^c\subset\ct^-_c$
be any object.
By 
\cite[Definition~\ref{D29.1}]{Neeman18A}
there must exist in $\ct$ a triangle
$A\la F(B)\la C$, with $A\in\ct^{\leq-m+1}$ and
with $C\in\ct^-_c\cap\ct^{\geq-m-N+1}$.

Because $C$ belongs to $\ct^-_c$, we may
choose in $\ct^c$ a Cauchy sequence
$c_*$ with $C=\hoco F(c_*)$, with
the notation as
in \cite{Neeman18A}. Since
$F(B)$ is compact in $\ct$, the morphism
$F(B)\la C=\hoco F(c_*)$ must factor
through some $F(c_i)$ for $i>0$; replacing
$c_*$ by a subsequence, we may
assume that the map factors
through $F(c_1)$. We have a map
of Cauchy sequences $B\la c_*$
in the category
$\cs=\ct^c$. Now
$C=\hoco F(c_*)$
belongs to
\[
\ct^{\geq-m-N+1}\eq(\ct^{\leq-m-N})^\perp
\sub F(\cm_{m+N})^\perp\ ,
\]
and the combination of the inclusion
$C\in F(\cm_{m+N})^\perp$
and 
\cite[Observation~\ref{O21.-1}]{Neeman18A}
allows us to
deduce that $\cy(C)\in Y(\cm_{m+N})^\perp$.
And now Lemma~\ref{L3.2.5}(ii) permits
us to conclude that
\[
\colim\,Y(c_*)\eq\cy(C)\quad\in\quad\fl(\cs)\cap\cl_{m+N}^\perp\eq
\fs(\cs)\cap\cl_{m+N}^\perp\ .
\]
It remains to show that the map of
Cauchy sequences $B\la c_*$ is of
type-$m$.

Now in the notation of
\cite[Lemma~8.5]{Neeman24}
we may replace the
Cauchy sequence
$F(c_*)$ by a subsequence
which is a strong $\ct^c$--approximating
sequence. And \cite[Lemma~8.5(iii)]{Neeman24}
informs us that it is a strong
$\ct^c$ approximating sequence
for $C=\hoco F(c_*)$.
In particular: for $i\gg0$ we have
that, in the triangle $F(c_i)\la C\la d_i$, the
inclusion $d_i\in\ct^{\leq-m-1}$ holds.
We are also given
the triangle $A\la F(B)\la C\la\T A$ with
$\T A\in\ct^{\leq-m}$. Now the octahedron
on the comopsable morphisms
$F(B)\la F(c_i)\la C$ gives the diagram
\[\xymatrix{
&
\Tm d_i\ar@{=}[r]\ar[d] &  
\Tm d_i \ar[d] \\
F(b)\ar[r]\ar@{=}[d]&
F(c_i)\ar[r]\ar[d] &
F(\ov d_i)\ar[d] \\
F(b)\ar[r] & C\ar[r] & \T A
}\]
and the triangle $\Tm d_i\la F(\ov d_i)\la\T A$,
coupled with the fact that
$\Tm d_i$ and $\T A$ both belong
to $\ct^{\leq-m}$, allows us to
deduce that
$F(\ov d_i)\in\ct^{\leq-m}*\ct^{\leq-m}\subset\ct^{\leq-m}$.
This makes $F(\ov d_i)$ an object in
$F(\cs)\cap\ct^{\leq-m}$, meaning
$\ov d_i\in\cm_m$ for all $i\gg0$.
\eprf

Example~\ref{E3.39393901} delivers one class of
excellent metrics. In this article we will study
the formal properties of these metrics, and
at the end we will find another major
class of examples.

Let us not get too far ahead of ourselves.
We note first that
Definition~\ref{D3.3}~(ii) and (iii) are about objects
$F\in\cs$, but both
extend easily to objects $F\in\fl(\cs)$.

\lem{L3.4}
Let $\cs$ be a triangulated category, and assume that 
$\{\cm_i\mid i\in\nn\}$
is an excellent metric on $\cs$,
as in Definition~\ref{D3.3}.
Choose any object $F\in\fl(\cs)$,
and express it as $F=\colim\,Y(f_*)$ for some
Cauchy sequence $f_*$ in $\cs$.
Given any integer $m>0$, then
the integer $n>m$ of Definition~\ref{D3.3}(ii)
has the property that the object $F\in\fl(\cs)$
admits a type-$m$ morphism $Y(E)\la F$ with respect to
$(E,f_*)$ analoguous
to that of
Definition~\ref{D3.3}(ii). That is: first of
all $E\in\cs$ belongs to ${^\perp\cm_n}$, and
secondly the Cauchy
sequence $E$ is the trivial
$E\stackrel\id\la E\stackrel\id\la E\stackrel\id\la\cdots$.
\elem

\prf
We are given an object $F\in\fl(\cs)$ and
Cauchy sequence $f_*$ in $\cs$ with
$F=\colim\,Y(f_*)$. And, with the given
integer $m$ as in the statement of the Lemma and
the integer 
$n>m$ as in Definition~\ref{D3.3}(ii), we may replace
$f_*$ by a subsequence such that, for any integers
$j\geq i\geq1$, when we complete $f_i\la f_j$ to a triangle
$f_i\la f_{j}\la y$ we have $y\in\cm_{m}$.

Now apply Definition~\ref{D3.3}(ii) to the
object $f_1\in\cs$ to obtain, in the category $\cs$,
a distinguished
triangle $E\stackrel{\ph_1^{}}\la f_1\la d_1$ with $E\in{^\perp\cl_n}$ and
with $d_1\in\cm_{m}$. Let $\ph:Y(E)\la F$ be the composite
\[\xymatrix@C+30pt{
Y(E)\ar[r]^-{Y(\ph_1^{})} & Y(f_1)\ar[r] &\colim\,Y(f_*)\ar@{=}[r]  & F\ .
}\]
Lemma~\ref{L3.307} now applies, with $k=\ell=1$
and to the commutative square
\[\xymatrix@C+20pt{
Y(E)\ar[r]\ar@{=}[d] & Y(f_1)\ar[d] \\
Y(E)\ar[r]^-{\ph} & F
}\]
to tells us that $\ph:Y(E)\la F$ is a type-$m$ morphism
with respect to $(E,f_*)$.
\eprf

\rmk{R3.?!?}
In the statement of Lemma~\ref{L3.4} we formulated
the condition on $E\in\cs$ as being $E\in{^\perp\cm_n}$.
But $Y$ is fully faithful, and hence
\[
\Hom(E,\cm_n)\cong\Hom\big(Y(E),Y(\cm_n)\big)\ .
\]
This gives the first of the two
equivalences below
\[
\{E\in{^\perp\cm_n}\}\quad\Longleftrightarrow\quad
\big\{Y(E)\in{^\perp Y(\cm_n)}\big\}\quad\Longleftrightarrow\quad
\big\{Y(E)\in{^\perp\cl_n}\big\}\ ,
\]
while the second equivalence follows from
the equality
$Y(\cs)\cap {^\perp Y(\cm_n})=Y(\cs)\cap{^\perp\cl_n}$
of Lemma~\ref{L3.2.5}(i).

To sum it up: the condition on $E\in\cs$
in the statement of Lemma~\ref{L3.4}
could
be rephrased as $Y(E)\in{^\perp\cl_n}$.
\ermk

\rmk{R3.?!?2}
Let $\cs$ be a triangulated category, and assume that 
$\{\cm_i\mid i\in\nn\}$
is an excellent metric on $\cs$.
By Definition~\ref{D3.3}(i) we know that
$\cs=\bigcup_{n\in\nn}{^\perp\cm_n}$,
and from
the
equivalence 
\[
\{E\in{^\perp\cm_n}\}\quad\Longleftrightarrow\quad
\big\{Y(E)\in{^\perp\cl_n}\big\}\ ,
\]
of Remark~\ref{R3.?!?}, we deduce that
\[
Y(\cs)\sub \bigcup_{n\in\nn}{^\perp\cl_n}\ .
\]
Also: given any object
$m\in\cm_n$, then we can let $m_*$ be the
constant Cauchy sequence
$m\stackrel\id\la m\stackrel\id\la m\stackrel\id\la\cdots$\ \ 
Hence $Y(m)=\colim\,Y(m_*)$ belongs
to $\cl_n$. 
We conclude that
\[
Y(\cm_n)\sub Y(\cs)\cap\cl_n\sub\left[\bigcup_{i\in\nn}{^\perp\cl_i}\right]\cap\cl_n\eq\bigcup_{i\in\nn}\big[{^\perp\cl_i}\cap\cl_n\big] \ .
\]
\ermk

In case the reader is wondering why, in the statement
of Lemma~\ref{L3.4}, we insisted on taking the constant
Cauchy sequence $e_*=E$ for $Y(E)$, we include the
following consequence.

\cor{C3.?!?!}
Let $\cs$ be a triangulated category, and assume that 
$\{\cm_i\mid i\in\nn\}$
is an excellent metric on $\cs$.
Let $\ell>0$ be an integer, 
choose an object $F\in\cl_\ell\subset\fl(\cs)$,
and express it as $F=\colim\,Y(f_*)$ for some
Cauchy sequence $f_*$ in $\cm_\ell$.
Given any integer $m\geq\ell+1$, then
the integer $n>m$ of Definition~\ref{D3.3}(ii)
has the property that, in the
morphism $Y(E)\la F$ produced by
Lemma~\ref{L3.4}, we automatically have
that $E\in{^\perp\cm_n}\cap\cm_\ell$.
\ecor

\prf
Lemma~\ref{L3.4} produces for us a type-$m$ morphism
$\ph:Y(E)\la F$ with respect to the Cauchy sequences
$(E,f_*)$. By
\cite[Lemma~\ref{L28.19}(i)]{Neeman18A} we may, after
passing to subsequences, construct a map
of Cauchy sequences $\ph_*:E\la f_*$ with
$\ph=\colim\,Y(\ph_*)$, and
then
\cite[Lemma~\ref{L28.19}(ii)]{Neeman18A} permits us to continue
this to a Cauchy sequence of triangles
$E\stackrel{\ph_*}\la f_*\la d_*\la\T E$.
And the assertion that $\ph$ is of type-$m$
says that, for $i\gg0$, we have $d_i\in\cm_m\subset\cm_{\ell+1}$,
where the inclusion is because $m\geq \ell+1$.

Now choose $i\gg0$ large. As $d_i\in\cm_{\ell+1}$ it
follows that $\Tm d_i\in\cm_\ell$. In the triangle
$\Tm d_i\la E\la f_i\la d_i$ we have that $\Tm d_i$
and $f_i$ both belong to $\cm_\ell$, and hence 
$E\in\cm_\ell*\cm_\ell\subset\cm_\ell$.

Of course: the fact that $E$ belongs to $^\perp\cm_n$
is given as part of the assertion of Lemma~\ref{L3.4}.
\eprf

\lem{L3.5}
Let $\cs$ be a triangulated category, and assume that 
$\{\cm_i\mid i\in\nn\}$
is an excellent metric on $\cs$,
as in Definition~\ref{D3.3}.
Choose any object $F\in\fl(\cs)$,
and express it as $F=\colim\,Y(f_*)$ for some
Cauchy sequence $f_*$ in $\cs$.
Given any integer $m>0$, then
the integer $n>0$ of Definition~\ref{D3.3}(iii)
has the property that there exists in $\fl(\cs)$
a type-$m$ morphism $F\la D$ with respect
to $(f_*,d_*)$, where the only
condition on the Cauchy sequence
$d_*$ and its colimit $D=\colim\,Y(d_*)$
is that $D\in\fs(\cs)\cap\cl_n^\perp$.
\elem

\prf
We are given an object $F\in\fl(\cs)$ and
Cauchy sequence $f_*$ in $\cs$ with
$F=\colim\,f_*$.
Let the given integer $m>0$ be as in the assertion
of the Lemma, and the integer $n>m$ 
as in Definition~\ref{D3.3}(iii).
Replacing $f_*$ by a subsequence we may assume
that, for
all pairs of integers $j\geq i\geq 1$, the
triangles
$f_i\la f_{j}\la y$ are such that $\Tm y$ and $y$
both lie in $\cm_{n}$.

Now, using Definition~\ref{D3.3}(iii) for the
object $f_1^{}\in\cs$, we construct a
type-$m$ morphism $\ph:Y(f_1)\la D$ with respect
to $(f_1,d_*)$, where $D\in\cl_n^\perp$ and
$f_1$ is the constant sequence.
Replacing $d_*$ by a subsequence, we may assume
that for all $1\leq k<\ell$ and every
triangle $d_k\la d_\ell\la z$, we have
$z\in\cm_m$.
Now
as in
\cite[Lemma~\ref{L28.19}(i)]{Neeman18A} we pass
to a subsequence yet again, and express
$Y(f_1)\la D$ as $\colim\,(\ph_*)$ for a map
of Cauchy sequences $\ph_*:f_1\la d_*$.
The fact that $\ph:Y(f_1)\la D$ is
a type-$m$ morphism with respect
to $(f_1,d_*)$ tells us
that, in any completion to a Cauchy sequence
of triangles $e_*\la f_1\la d_*\la\T e_*$,
we have that $\T e_i^{}\in\cm_m$ for $i\gg0$.
Replacing by a subsequence we may assume
that $e_1^{}\la f_1\la d_1\la\T e_1^{}$ has
$\T e_1^{}\in\cm_m$.

Now consider the composite
\[\xymatrix@C+30pt{
Y(f_1)\ar[r]^-{Y(\ph_1)} & Y(d_1)\ar[r] &\colim\,Y(d_*)\ar@{=}[r]  & D\ .
}\]
Because in the triangles $\Tm y_i\la f_i\la f_{i+1}\la y_i$
we have that $\T y_i,y_i\in\cm_n$,
and furthermore $D\in\cl_n^\perp=\cm_n^\perp$ is
a homological functor on $\cs$, 
\cite[Lemma~\ref{L28.102}]{Neeman18A}
applies and allows us to factor
the map $Y(f_1)\la D$, uniquely, as
$Y(f_1)\la F\stackrel\psi\la D$.
We are
in a position to apply Lemma~\ref{L3.307} again,
once again with $k=\ell=1$, and
this time to the commutative square
\[\xymatrix@C+20pt{
Y(f_1)\ar[r]\ar[d] & Y(d_1)\ar[d] \\
F\ar[r]^-{\psi} & D
}\]
We deduce that $\psi:F\la D$ is a type-$m$ morphism
with respect to $(f_*,d_*)$.
\eprf

Lemmas~\ref{L3.4} and \ref{L3.5} produce for us
a plethora of type-$m$ morphisms. We observe first:

\lem{L3.309}
Let $\cs$ be a triangulated category, and assume that 
$\{\cm_i\mid i\in\nn\}$
is a good metric on $\cs$.
Assume we
are given in $\fl(\cs)$ a type-$(m+1)$ morphism
$\psi:F\la G$, with respect to some
pair of Cauchy sequences $(f_*,g_*)$.
Then there exists in $\fl(\cs)$ a strong triangle
$E\la F\stackrel\psi\la G\la\T E$ with
both $E$ and $\T E$ in $\cl_m$.
\elem

\prf
By hypothesis we may,
up to replacing $f_*$ and $g_*$ by
subsequences, choose
a type-$(m+1)$ morphism of
Cauchy sequences $\psi_*:f_*\la g_*$
with $\psi=\colim\,(\psi_*)$. Now
\cite[Lemma~\ref{L28.19}(ii)]{Neeman18A} permits us
to extend to a Cauchy sequence
of triangles $e_*\la f_*\stackrel{\psi_*}\la g_*\la\T e_*$,
and the hypothesis that $\psi_*$ is a type-$(m+1)$
morphism says
that, after replacing by a subsequence, we may assume
$\T e_*$ lies in $\cm_{m+1}$, and hence $e_*$ lies
in $\cm_m$.

Taking colimits delivers a strong triangle
$E\la F\stackrel\psi\la G\la\T E$ with
both $E$ and $\T E$ in $\cl_m$.
\eprf

From this we deduce a useful consequence:

\cor{C3.378}
Let $\cs$ be a triangulated category, and assume that 
$\{\cm_i\mid i\in\nn\}$
is a good metric on $\cs$.
Assume we
are given
in $\fl(\cs)$ a type-$(m+1)$ morphism $\psi:F\la G$,
with respect to some Cauchy sequences $(f_*,g_*)$.
Then $\Hom(A,-)$ takes
$\psi$ to an isomorphism whenever $A\in{^\perp\cl_{m}}\cap\cs$
and $\Hom(-,B)$ takes $\psi$ to an
isomorphism whenever $B\in\fs(\cs)\cap\cl_{m}^\perp$.
\ecor

\prf
By Lemma~\ref{L3.309} the morphism
$\psi:F\la G$ may be extended
in $\fl(\cs)$ to a strong triangle
$E\la F\stackrel\psi\la G\la\T E$ with
both $E$ and $\T E$ in $\cl_m$, 
and both $E$ and $\T E$ are annihilated by
$\Hom(A,-)$ when $A\in{^\perp\cl_{m}}$ and
by $\Hom(-,B)$ when $B\in\cl_{m}^\perp$.

Thus the Corollary comes from the fact that the
functor $\Hom(A,-)$ takes strong
triangles to long exact sequences as
long as $A\in\cs$, and the functor
$\Hom(-,B)$ takes strong triangles
to long exact sequences whenever $B\in\fs(\cs)$.
See Reminder~\ref{R3.-111222333} for the assertion
about $\Hom(A,-)$, and 
Lemma~\ref{L3.303} for the statement about
$\Hom(-,B)$.
\eprf

\lem{L3.7}
Let $\cs$ be a triangulated category, and assume that 
$\{\cm_i\mid i\in\nn\}$
is an excellent metric on $\cs$.
Consider also the category $\fs(\cs)\op$,
with its good metric
$\{\cn\op_i,\,i\in\nn\}$
of Lemma~\ref{L3.1}.

Then there
exists an increasing sequence of integers
\[
2=n_1^{}<n_2^{}<n_3<\cdots
\]
such that the following hold.
\be
\item
Given an object $F\in\fl(\cs)$,
a Cauchy sequence $f_*$ in $\cs$ with
$F=\colim\,Y(f_*)$,
and an integer $i>0$, there exist
an object $g_i^{}\in\fl(\cs)\cap\cl^\perp_{n_{i+1}^{}}$
and a Cauchy sequence $g_{i,*}^{}$ in $\cs$
with $g_i=\colim\,Y(g_{i,*}^{})$,
as well as, in the category $\fl(\cs)$,
a type-$(1+n_i)$
morphism
$\psi_i:F\la g_i$
with respect to $(f_*,g_{i,*}^{})$.
\item
Furthermore:
if
the sequence $f_*$ lies in $\cm_k$
for some $k>0$, then $g_i$
can be chosen to lie in $\cn_k$ for all $i\geq k$.
\item
Finally: the $g_i$ can be assembled, uniquely,
into an inverse sequence $\cdots\la g_3^{}\la g_3^{}\la g_1^{}$
which is Cauchy with respect to the metric
$\{\cn\op_j,\,j\in\nn\}$ on $\fs(\cs)\op$,
in such a way that the
maps $\psi_i:F\la g_i$ assemble to
a morphism $\psi:F\la\clim g_*$.
\ee
\elem

\prf
The sequence of integers is easy to construct:
we are given that $n_1^{}=2$ and,
assuming we have constructed $n_i$, the recipe
for $n_{i+1}^{}$ is to set $m=1+n_i$ and let
$n_{i+1}^{}=n$ be the integer $n>m$ whose existence is
guaranteed by 
Definition~\ref{D3.3}(iii).

In (i) we fix an
object $F\in\fl(\cs)$ and a Cauchy sequence
$f_*$ in $\cs$ with $F=\colim\, Y(f_*)$.
Now apply
Lemma~\ref{L3.5} to the object $F=\colim\, f_*$ with
$m=1+n_i$. By the choice of 
the integer $n_{i+1}^{}>1+n_i$,
there exists,
in the category
$\fl(\cs)$, a type-$(1+n_i)$ morphism
$F\la g_i^{}$ with respect to a pair of
Cauchy sequences $(f_*,g_{i,*}^{})$, and such that
$g_i^{}=\colim\,Y(g_{i,*}^{})$ satisfies
$g_i^{}\in\cl_{n_{i+1}^{}}^\perp\cap\fs(\cs)$.

Now for (ii): observe that
\cite[Lemma~\ref{L28.19}(i)]{Neeman18A}
tells us that,
after replacing
$(f_*,g_{i,*}^{})$ by subsequences, we may
express $F\la g_i^{}$ as a colimit of the 
image under $Y$ of a Cauchy sequence of morphisms
$\ph_*:f_*\la g_{i,*}^{}$. The fact
that $F\la g_i^{}$ is a type-$(1+n_i)$ morphism
with respect to $(f_*,g_{i,*}^{})$ means
that, for any choice of $\ph_*:f_*\la g_{i,*}^{}$,
if we use
\cite[Lemma~\ref{L28.19}(ii)]{Neeman18A}
to complete $\ph_*$  to a Cauchy
sequence of triangles $e_{i,*}\la f_*\la g_{i,*}\la\T e_{i,*}$,
the objects
$\T e_{i,\ell}^{}$ will be contained in
$\cm_{1+n_i}\subset\cm_{1+i}$ for $\ell\gg0$. 
Passing to a subsequence we may assume that this is
true for all $\ell$.

Note: we are assuming that $f_\ell\in\cm_k$ for all $\ell$,
and by the above $\T e_{i,\ell}^{}\in\cm_{i+1}$ for every $\ell$.
If $i\geq k-1$, the
triangle $f_\ell\la g_{i,\ell}\la\T e_\ell^{}$
tells us that  $g_{i,\ell}^{}\in\cm_k*\cm_k\subset\cm_k$.
Therefore $g_i=\colim\,Y(g_{i,*})$ lies in $\cl_k$.
But as it lies in $\fs(\cs)$ by (i), it belongs
to $\cn_k=\cl_k\cap\fs(\cs)$.

It remains to prove (iii).
We started out with an object $F$ in the
category of $\fl(\cs)$,
and a Cauchy sequence $f_*$ in
$\cs$ with $F=\colim\,Y(f_*)$.
And back in (i) we constructed,
for every integer $i\geq1$, an object
$g_i\in\cl_{n_{i+1}^{}}^\perp$, a Cauchy sequence
$g_{i,*}$ in $\cs$ such that $g_i=\colim\,g_{i,*}$,
and a type-$(1+n_{i})$ 
morphism
$\psi_{i}:F\la g_{i}^{}$
with respect to $(f_*,g_{i,*})$.
Applying $\Hom(-,g_i)$ to the morphism
$\psi_{i+1}:F\la g_{i+1}^{}$ gives an isomorphism by
Corollary~\ref{C3.378}, and hence
$\psi_i:F\la g_i$ factors uniquely
as 
\[\xymatrix@C+50pt@R-10pt{
F\ar[rr]^-{\psi_i} \ar[rd]_-{\psi_{i+1}^{}} & &  g_i \\
 & g_{i+1}^{}\ar[ur]_-{\exists!\gamma_i} &
}\]
We know that $\psi_{i}=\gamma_i\circ\psi_{i+1}^{}$
is a type-$(1+n_i)$ morphism
with respect to $(f_*,g_{i,*})$,
while $\psi_{i+1}^{}$ is 
a type-$(1+n_{i+1}^{})$ 
morphism 
with respect to $(f_*,g_{i+1,*}^{})$.
Lemma~\ref{L3.341}(ii) tells us first
that $\gamma_i$ must be a type-$(1+n_i)$
morphism with respect to
$(g_{i+1,*}^{},g_{i,*})$, and when we combine with
Lemma~\ref{L3.309} we deduce that there
exists in $\fl(\cs)$ a strong
triangle $E_i\la g_{i+1}\la g_{i}^{}\la\T E_i$ with
$E$ and $\T E$ both in $\cl_{n_i}$.
But $g_i$ and $g_{i+1}^{}$ both
belong to $\fs(\cs)$, and by
\cite[Remark~\ref{R200976}]{Neeman18A}
the strong
triangle $E_i\la g_{i+1}\la g_{i}^{}\la\T E_i$ must
be a distinguished triangle in $\fs(\cs)$.
Hence we conclude that $E_i$ must
belong to $\fs(\cs)\cap\cl_{n_i}=\cn_{n_i}$,
making the sequence $g_*$ a Cauchy sequence
in $\fs(\cs)\op$.

In the category
$\MMod\cs$ we have produced a commutative diagram
\[\xymatrix@C+10pt@R-0pt{
& & &F\ar[lld]|-{\psi_5^{}} 
\ar[ld]|-{\psi_4^{}} \ar[d]|-{\psi_3^{}}
\ar[rd]|-{\psi_2^{}} \ar[rrd]|-{\psi_1^{}}\\
\cdots\ar[r]& g_5\ar[r]_-{\gamma_5^{}} &
g_4\ar[r]_-{\gamma_3^{}} &
g_3\ar[r]_-{\gamma_2^{}} &
g_2\ar[r]_-{\gamma_1^{}} &
g_1
}\]
giving an inverse sequence $g_*$ 
and a map $\psi:F\la\clim g_*$. This completes the
proof of (iii).
\eprf

\dis{D3.6.5}
We remind the reader of standard notation: 
if $\ca$ is an additive category
then $\MMod\ca$ stands for
the category of right $\ca$--modules,
that is additive functors $\ca\op\la\ab$, while
$\Mod\ca$ stands for the category
of left $\ca$--modules,
meaning additive functors $\ca\la\ab$.
Thus $\MMod{\ca\op}=\Mod\ca$.

Now let $\cs$ be a triangulated category, and assume that 
$\{\cm_i\mid i\in\nn\}$
is a good metric on $\cs$. Then we may form
the category $\fs(\cs)\subset\MMod\cs$, and
deduce a functor
\[\xymatrix@C-10pt{
\wh Y\ar@{}[r]|-{\ds:} &
\big[\MMod\cs\big]\op\ar[rrrr] &&&& \MMod{\fs(\cs)\op}
\ar@{=}[r]&\Mod{\fs(\cs)}
}\]
by taking an object $M\in\MMod\cs$ to the
restriction to $\fs(\cs)\subset\MMod\cs$
of the functor $\Hom_{\MMod\cs}^{}(M,-)$.
In symbols:
\[\xymatrix@C-10pt{
\wh Y(M)\ar@{=}[r] &\Hom_{\MMod\cs}^{}(M,-)\Big|^{}_{\fs(\cs)}\ar@{}[r]|-{\ds:} &
\fs(\cs)\ar[rrr] &&& \ab\ .
}\]
\edis

\lem{L3.90909}
Let $\cs$ be a triangulated category, and assume that 
$\{\cm_i\mid i\in\nn\}$
is an excellent metric on $\cs$.
Consider also the category $\fs(\cs)\op$,
with its good metric
$\{\cn\op_i,\,i\in\nn\}$
of Lemma~\ref{L3.1}.

Suppose we are given:
\be
\item
An object $F\in\fl(\cs)$.
\item
A Cauchy sequence $g_*$ in the category
$\fs(\cs)\op$.
\item
For each integer $i>0$ we are
given a morphism
$\psi_i:F\la g_i$,
which assemble to a commutative
diagram
\[\xymatrix@C+10pt@R-0pt{
& & &F\ar[lld]|-{\psi_5^{}} 
\ar[ld]|-{\psi_4^{}} \ar[d]|-{\psi_3^{}}
\ar[rd]|-{\psi_2^{}} \ar[rrd]|-{\psi_1^{}}\\
\cdots\ar[r]& g_5\ar[r]_-{\gamma_5^{}} &
g_4\ar[r]_-{\gamma_3^{}} &
g_3\ar[r]_-{\gamma_2^{}} &
g_2\ar[r]_-{\gamma_1^{}} &
g_1
}\]
Assume further that each morphism $\psi_i:F\la g_i$ can
be completed in $\fl(\cs)$ to a
strong triangle $E_i\la F\stackrel{\psi_i}\la g_i\la \T E_i$
with $E_i$ and $\T E_i$ both in $\cl_i$.
\setcounter{enumiv}{\value{enumi}}
\ee
Then the map $\psi:F\la\clim g_*$ is an isomorphism
in $\MMod\cs$, while the map
$\ph:\colim\,\wh Y(g_*)\la \wh Y(F)$ is
an isomorphism in $\Mod{\fs(\cs)}=\MMod{\fs(\cs)\op}$.
\elem

\prf
In the category
$\MMod\cs$ we are given a commutative diagram,
providing an inverse sequence $g_*$ 
and compatible maps $\psi_i:F\la g_i$. We need
to prove two things:
\be
\setcounter{enumi}{\value{enumiv}}
\item
For every object $s\in\cs$, the 
functor $\Hom\big(Y(s),-)$ takes
$\psi$ to an isomorphism.
\item
For object $Z\in\fs(\cs)$,
the functor $\Hom\big(-,Z)$ takes $\ph$ to an
isomorphism.
\ee
Let us begin with
(iv): choose any object $s\in\cs$.
Because the metric
$\{\cm_i,\,i\in\nn\}$ is excellent,
Definition~\ref{D3.3}(i) tells us that there
must exist an integer $\ell>0$ with $s\in{^\perp\cm_\ell}$.
Lemma~\ref{L3.2.5}(i) informs
us that $Y(s)\in{^\perp\cl_\ell}$, but if $i\geq\ell$ then
the strong triangle
$E_i\la F\stackrel{\psi_i}\la g_i\la \T E_i$
has both $E_i$ and $\T E_i$ lying in $\cl_i\subset\cl_\ell$.
Hence $\Hom\big(Y(s),-)$ annihilates
both $E_i$ and $\T E_i$. As $\Hom\big(Y(s),-)$
takes any strong triangle to a long exact
sequence, we deduce that it must
take $\psi_i$ to an isomorphism for all
$i\geq\ell$. The assertion of (iv) comes
by passing to the limit.

Now for (v). Choose any object $Z\in\fs(\cs)$.
By the definition of $\fs(\cs)$ the object
$Z$ must lie
in $\cm_\ell^\perp$ for some $\ell>0$, but
Lemma~\ref{L3.2.5}(ii) teaches us that
$\cm_\ell^\perp=\cl_\ell^\perp$.
If $i\geq\ell$ then
the strong triangle
$E_i\la F\stackrel{\psi_i}\la g_i\la \T E_i$
is such that
$E_i$ and $\T E_i$ both belong to $\cl_i\subset\cl_\ell$, and
hence $\Hom\big(-,Z)$ annihilates
both $E_i$ and $\T E_i$. Lemma~\ref{L3.303}
tells us that the functor $\Hom\big(-,Z)$
takes any strong triangle to a long exact
sequence, and we deduce that it must
take $\psi_i$ to an isomorphism for all
$i\geq\ell$. The assertion of (v) 
comes
by passing to the colimit.
\eprf

\cor{C3.91919}
Combining Lemmas~\ref{L3.7} and \ref{L3.90909},
we obtain:
\be
\item
The functor $\wh Y:\big[\MMod\cs\big]\op\la\MMod{\fs(\cs)\op}$
takes the subcategory $\fl(\cs)\op\subset\big[\MMod\cs\big]\op$ into the subcategory
$\fl\big[\fs(\cs)\op\big]\subset\MMod{\fs(\cs)\op}$.
\item
If $\wh\cl_n\subset\fl\big[\fs(\cs)\op\big]$
is just the construction of Definition~\ref{D3.-3}(ii),
done on the triangulated category
$\fs(\cs)\op$ with respect to its good metric
$\{\cn\op_i,\,i\in\nn\}$,
then $\wh Y(\cl_n\op)$ is contained in $\wh\cl_n$.
\ee
\ecor

\prf
Given an object $F\in\fl(\cs)$, Lemma~\ref{L3.7}
constructs for us the Cauchy sequence
$g_*$ in the category $\fs(\cs)$,
with maps $\psi_i:F\la g_i$, with the data
satisfying the hypotheses of 
Lemma~\ref{L3.90909}. From the conclusion of
the Lemma we deduce that
$\wh Y(F)=\colim\,\wh Y(g_i)$
lies in $\fl\big(\fs(\cs)\op\big)$.

If $F$ belongs to $\cl_\ell$, then the Cauchy sequence
$g_*$ constructed in Lemma~\ref{L3.7}(ii) may be
chosen to lie in $\cn_\ell$. Hence the object
$\wh Y(F)=\colim\,\wh Y(g_i)$ must lie in $\wh\cl_\ell$.
\eprf

\lem{L3.92929}
Let $\cs$ be a triangulated category, and assume that 
$\{\cm_i\mid i\in\nn\}$
is an excellent metric on $\cs$.
Suppose we are given in the category
$\fs(\cs)\op$ a Cauchy sequence $g_*$
with respect to the metric $\{\cn\op_i,\,i\in\nn\}$.
Translating from the category $\fs(\cs)\op$ 
to the category $\fs(\cs)$: we are given in
$\fs(\cs)$ an
inverse sequence $\cdots\la g'_3\la g'_2\la g'_1$
such that, for any integer $n>0$, there
exists an integer $N>0$ such that, for all $j>i\geq N$, 
in every triangle $d'_{i,j}\la g'_{j}\la g'_i$
we have $d'_{i,j}\in\cn_n$.

Then, up to replacing $g'_*$ by a subsequence $g_*$,
we may find an object $F\in\fl(\cs)$ and
a commutative diagram
\[\xymatrix@C+10pt@R-0pt{
& & &F\ar[lld]|-{\psi_5^{}} 
\ar[ld]|-{\psi_4^{}} \ar[d]|-{\psi_3^{}}
\ar[rd]|-{\psi_2^{}} \ar[rrd]|-{\psi_1^{}}\\
\cdots\ar[r]& g_5\ar[r]_-{\gamma_5^{}} &
g_4\ar[r]_-{\gamma_3^{}} &
g_3\ar[r]_-{\gamma_2^{}} &
g_2\ar[r]_-{\gamma_1^{}} &
g_1
}\]
satisfying the hypotheses of Lemma~\ref{L3.90909}.
Moreover: if the $g'_i$ all belong to $\cn_\ell$ for some
$\ell>0$, then $F$ can be constructed to belong to
$\cl_\ell$.
\elem

\prf
As in Lemma~\ref{L3.7} we begin by choosing an
increasing sequence
of integers 
\[
2=n_1^{}<n_2^{}<n_3<\cdots
\]
only this time we do it as in Definition~\ref{D3.3}(ii).
The integer $n_1^{}=2$ is given, and the construction is
that, given the integer $n_i$, we let $m=1+n_i$
and choose $n_{i+1}^{}=n$ to be the number $n>m$
whose existence is guaranteed
in Definition~\ref{D3.3}(ii).

Now: because the
sequence $g'_*$ is Cauchy, for each
integer $i>0$ there exists an integer $N_i\gg0$
such that $j>i\geq N_i$ implies that,
in the triangles
$g'_{j}\la g'_i\la\T d'_{i,j}$, we have 
that $\T d'_{i,j}$ is contained
in $\cn_{1+n_{i+1}^{}}$. And we may of course
choose
the $N_i$ to form a strictly  increasing sequence
of integers. Let $g_i=g'_{N_i}$. This
has defined for us the subsequence $g_*$ of
$g'_*$.

Next,
choose any Cauchy sequence $g_{1,*}^{}$ in $\cs$
with $g_1^{}=\colim\,Y(g_{1,*}^{})$.
If the sequence $g_*$ lies in $\cn_\ell$,
then we begin the construction with $g_\ell^{}$
instead of $g_1^{}$, and assume the sequence
$g_{\ell,*}^{}$ with $g_\ell^{}=\colim\,Y(g_{\ell,*}^{})$
is chosen to lie in $\cm_\ell$.
We will construct by induction Cauchy
sequences
$g_{i,*}^{}$ in $\cs$, so that
$g_i^{}=\colim\,Y(g_{i,*}^{})$.
Suppose that, in the category $\cs$,
the sequence $g_{i,*}^{}$ has been constructed.
Remember the triangle $g_{i+1}^{}\la g_i\stackrel\ph\la\T d_i$
in $\fs(\cs)$; by construction $\T d_i$ belongs
to $\cn_{1+n_{i+1}^{}}$, and hence we can express it
as $\colim\,Y(\T d_{i,*})$ with $\T d_{i,*}$ a
Cauchy sequence in $\cm_{1+n_{i+1}^{}}$.
If $i\geq\ell$ then
we have $\cm_{1+n_{i+1}^{}}\subset\cm_{1+\ell}$.
Passing to subsequences we may by
\cite[Lemma~\ref{L28.19}(i)]{Neeman18A} choose
a morphism of Cauchy sequences $\ph_*:g_{i,*}^{}\la\T d_{i,*}$
with $\ph=\colim\,Y(\ph_*)$,
and by
\cite[Lemma~\ref{L28.19}(ii)]{Neeman18A}
we may complete in $\cs$ to a Cauchy sequence
of triangles
$d_{i,*}\la g_{i+1,*}^{}\la g_{i,*}\stackrel{\ph_*}\la
\T d_{i,*}$. Note that
if $i\geq\ell$ then $d_{i,*}$ and $\T d_{i,*}$ both
lie in $\cm_\ell$, and if the sequence
$g_{i,*}$ was constructed to lie in $\cm_\ell$ then
$g_{i+1,*}^{}$ lies in $\cm_\ell*\cm_\ell\subset\cm_\ell$.
Applying the functor $Y$ to this Cauchy sequence
of triangles and taking colimits gives a
strong triangle
$d_i\la \ov g_{i+1}^{}\la g_i\stackrel\ph\la\T d_i$
in the category $\fl(\cs)$, but
as the objects $d_i$ and $g_i$ both lie
in $\fs(\cs)$,
\cite[Remark~\ref{R200976}]{Neeman18A}
guarantees that this is a distinguished triangle
in $\fs(\cs)$. Hence the two triangles
\[
d_i\la \ov g_{i+1}^{}\la g_i\stackrel\ph\la\T d_i\qquad\text{and}\qquad
d_i\la g_{i+1}^{}\la g_i\stackrel\ph\la\T d_i
\]
must be an isomorphic, and we have expressed
$g_{i+1}^{}$ as the colimit of the image
under $Y$ of a Cauchy sequence
$g_{i+1,*}^{}$, in such a way that the morphism
$\gamma_i:g_{i+1}^{}\la g_i$ is of type-$(1+n_{i+1}^{})$ with
respect to $(g_{i+1,*}^{},g_{i,*}^{})$.
And if the sequence $g_*$ lies in $\cn_\ell$ then
our construction started with
$g_\ell^{}=\colim\,Y(g_{\ell,*}^{})$,
and for $i\geq\ell$ all the Cauchy sequences
$g_{i,*}^{}$ lie in $\cm_\ell$.

And now we apply Lemma~\ref{L3.4} to each
$g_i=\colim\,Y(g_{i ,*})$, with $m=1+n_i$.
For each $i>0$ we obtain an object
$f_i\in\cs\cap{^\perp\cl_{n_{i+1}^{}}}$,
as well as a type-$(1+n_i)$ morphism
$\Psi_i:Y(f_i)\la g_i$ with respect to
the Cauchy sequences $(f_i,g_{i ,*})$, with
$f_i$ the constant sequence. In the
special case where $i\geq\ell$ and the
sequence $g_{i,*}$ was constructed to lie
in $\cm_\ell$, the fact that the
morphism $\Psi_i:Y(f_i)\la g_i$ is of type-$(1+n_i)$
with respect to
the Cauchy sequences $(f_i,g_{i ,*})$
tells us that there must exist an integer
$n\gg0$ and, in the category $\fl(\cs)$,
a commutative square
\[\xymatrix{
Y(f_i)
\ar[r]\ar@{=}[d] & Y(g_{i ,n})\ar[d]\\
Y(f_i)\ar[r] &g_i
}\]
such that, if we complete $f_i\la g_{i ,n}$ to
a triangle $\wh d\la f_i\la g_{i ,n}\la\T\wh d$ in $\cs$,
the object $\T\wh d$ belongs to $\cm_{1+n_i}\subset\cm_{1+\ell}$.
Therefore the object $\wh d$ belongs to $\cm_\ell$, and
hence the object $f_i$ has to belong to
$\cm_\ell*\cm_\ell\subset\cm_\ell$.

Now
for each $i$ we have composable morphisms
\[\xymatrix@C+50pt@R-10pt{
Y(f_{i+1})\ar[r]^-{\Psi_{i+1}}
& g_{i+1}^{}\ar[r]^-{\gamma_i} & g_i
}\]
with $\Psi_{i+1}$ of type-$(1+n_{i+1}^{})$ with
respect to $(f_{i+1},g_{i+1,*}^{})$ and
$\gamma_{i}$ of type-$(1+n_{i+1}^{})$ with
respect to $(g_{i+1,*}^{},g_{i,*}^{})$.
Lemma~\ref{L3.341}(i)
tells us that the composite
is a morphism $\gamma_i\Psi_{i+1}:Y(f_{i+1})\la g_i$
of type-$(1+n_{i+1}^{})$ with
respect to $(f_{i+1},g_{i,*}^{})$.
But as  $Y(f_i)\in Y(\cs)\cap{^\perp\cl_{n_{i+1}^{}}}$,
Corollary~\ref{C3.378} tells
us that $\Hom\big(Y(f_i),-\big)$ must take
$\gamma_i\Psi_{i+1}:Y(f_{i+1})\la g_i$
to an isomorphism: there is a unique
morphism $\alpha_i:f_i\la f_{i+1}$
rendering
commutative the diagram
\[\xymatrix@C+50pt@R-10pt{
Y(f_i)\ar[rrr]^-{\Psi_i} \ar[rd]_-{\exists!Y(\alpha_{i}^{})} && &  g_i \\
  & Y(f_{i+1})\ar[r]_-{\Psi_{i+1}}& g_{i+1}^{}\ar[ur]_-{\gamma_i} &
}\]
Since $\Psi_i$ is a type-$(1+n_i)$ morphism
with respect to
the Cauchy sequences $(f_i,g_{i ,*})$ and
$\gamma_i\Psi_{i+1}:f_{i+1}\la g_i$ is of type-$(1+n_{i+1}^{})$ with
respect to $(f_{i+1},g_{i,*}^{})$,
Lemma~\ref{L3.341}(ii) permits us to deduce that
$Y(\alpha_i):Y(f_i)\la Y(f_{i+1})$ is a type-$(1+n_i)$ morphism
with respect to
$(f_i,f_{i+1})$. For each integer $i>0$
we have produced a morphism
$\alpha_i:f_i\la f_{i+1}$, in the category
$\cs$. And the fact that
$Y(\alpha_i):Y(f_i)\la Y(f_{i+1})$ is a type-$(1+n_i)$ morphism
with respect to the \emph{constant sequences}
$(f_i,f_{i+1})$ tells us that, when we complete
to a triangle $f_i\stackrel{\alpha_i}\la f_{i+1}\la\ov d_i$
in the category $\cs$, the object $\ov d_i$ belongs of
$\cm_{1+n_i}$. Define the sequence $f_*$ to be
\[\xymatrix@C+10pt{
f_1\ar[r]^-{\alpha_1^{}} &
f_2\ar[r]^-{\alpha_2^{}} &
f_3\ar[r]^-{\alpha_3^{}} &
f_4\ar[r]^-{\alpha_4^{}} &
f_5\ar[r] &\cdots
}\]
and the construction is such that we have
have produced in $\fl(\cs)$ a
commutative diagram
\[\xymatrix@C+10pt{
Y(f_1)\ar[d]^{\Psi_1}\ar[r]^-{Y(\alpha_1^{})} &
Y(f_2)\ar[d]^{\Psi_2}\ar[r]^-{Y(\alpha_2^{})} &
Y(f_3)\ar[d]^{\Psi_3}\ar[r]^-{Y(\alpha_3^{})} &
Y(f_4)\ar[d]^{\Psi_4}\ar[r]^-{Y(\alpha_4^{})} &
Y(f_5)\ar[d]^{\Psi_5}\ar[r] &\cdots\\
g_1&\ar[l]^-{\gamma_1^{}} 
g_2&\ar[l]^-{\gamma_2^{}} 
g_3&\ar[l]^-{\gamma_3^{}} 
g_4&\ar[l]^-{\gamma_4^{}} 
g_5&\ar[l] \cdots
}\]
Now let $F=\colim\,Y(f_*)$; the
commutative diagram above allows us to factor
the map
$\Psi_i:Y(f_i)\la g_i$ as the composite
$Y(f_i)\la F\stackrel{\psi_i}\la g_i$.

The sequence $f_*$ is Cauchy in $\cs$ with
respect to the metric $\{\cm_i,\,i\in\nn\}$,
and its colimit $F$ belongs to $\fl(\cs)$.
If the Cauchy sequence $g'_*$ in $\fs(\cs)$
that we started with lies in $\cn_\ell$ for
some $\ell>0$, then the Cauchy sequence $f_*$ which
we constructed lies in $\cm_\ell$, and hence
$F=\colim\,Y(f_*)$ belongs to $\cl_\ell$.

Next we want to apply
Lemma~\ref{L3.307} to the
morphism $Y(f_i)\la F$.
We have Cauchy sequences $f_i$ and $f_*$ with
$Y(f_i)=\colim\,Y(f_i)$ and with
$F=\colim\,Y(f_*)$, we choose the pair integers
$k,\ell$ of Lemma~\ref{L3.307}(ii) to be 
$k=1$ and $\ell=i$, and we have
that for all $h>j\geq k$ the triangles
$f_i\stackrel\id\la f_i\la 0$ most definitely
satisfy $\T 0= 0\in\cm_{1+n_i}$,
and for all $h>j\geq\ell$ the triangles
$f_j\la f_h\la \wh d_{j,h}$ satisfy
$\wh d_{j,h}\in\cm_{1+n_i}$.
As in Lemma~\ref{L3.307}(i)
we also have the 
commutative square
\[\xymatrix{
Y(f_i)
\ar@[r][r]\ar@{=}[d] & Y(f_i)\ar[d]\\
Y(f_i)\ar[r] &F
}\]
and, for Lemma~\ref{L3.307}(iii) we
again observe that
the triangle $f_i\stackrel\id\la f_i\la 0$
satisfies the hypothesis that 
$0\in\cm_{1+n_i}$. Thus 
the conditions of Lemma~\ref{L3.307} hold, and
the conclusion of Lemma~\ref{L3.307} is that
the morphism $Y(f_i)\la F$ is of
type-$(1+n_i)$ with respect to
$(f_i,f_*)$.

Now we have a pair of composable morphisms
$Y(f_i)\la F\stackrel{\psi_i}\la g_i$. We know that
$Y(f_i)\la F$ is of type-$(1+n_i)$ with respect to
$(f_i,f_*)$, while the composite $\Psi_i:Y(f_i)\la g_i$
is of type-$(1+n_i)$ with respect to
$(f_i,g_{i,*})$. By Lemma~\ref{L3.341}(ii) the morphism
$\psi_i:F\la g_i$ is of type-$n_i$ with respect to
$(f_*,g_{i,*})$. Lemma~\ref{L3.309} permits
us to extend the morphism $\psi_i:F\la g_i$ to
a strong triangle
$E_i\la F\stackrel{\psi_i}\la g_i\la \T E_i$
with $E_i$ and $\T E_i$ both in $\cl_{n_i-1}\subset\cl_i$,
where the inclusion is because $n_i\geq i+1$.

This completes the proof of Lemma~\ref{L3.92929}.
\eprf

Combining Corollary~\ref{C3.91919}, with 
Lemmas~\ref{L3.92929} and \ref{L3.90909},
we obtain the following.

\cor{C3.93939}
The statements below hold:
\be
\item
The essential image of the subcategory
$\fl(\cs)\op\subset\big[\MMod\cs\big]\op$,
under the functor
$\wh Y:\big[\MMod\cs\big]\op\la\MMod{\fs(\cs)\op}$,
is precisely 
the subcategory
$\fl\big[\fs(\cs)\op\big]\subset\MMod{\fs(\cs)\op}$.
\item
The induced functor 
$\wh Y:\fl(\cs)\op\la \fl\big[\fs(\cs)\op\big]$
is full; that is every morphism of the form
$\wh Y(A)\la \wh Y(B)$
is equal to  $\wh Y(h)$ for some $h:B\la A$.
\item
The essential image
of $\cl\op_\ell\subset\fl(\cs)\op$,
under the functor
$\wh Y:\fl(\cs)\op\la \fl\big[\fs(\cs)\op\big]$,
is precisely
$\wh\cl_\ell\subset\fl\big[\fs(\cs)\op\big]$.
\ee
\ecor

\prf
Corollary~\ref{C3.91919}(i) shows that the essential
image of $\fl(\cs)\op$ under the functor $\wh Y$
is contained in $\fl\big[\fs(\cs)\op\big]$.
But every object in $\fl\big[\fs(\cs)\op\big]$
is isomorphic to $\colim\,\wh Y(g'_*)$ for some
Cauchy sequence $g'_*\in \fs(\cs)\op$,
and Lemmas~\ref{L3.92929} and \ref{L3.90909}
combine to exihibit an object
$F\in\fl(\cs)$ and a subsequence $g_*$ of $g'_*$
with $\wh Y(F)=\colim\,\wh Y(g_*)$. Thus the essential
image of $\fl(\cs)\op$ under the functor $\wh Y$
contains all of $\fl\big[\fs(\cs)\op\big]$.
This completes the proof of (i).

Now for (ii). By (i), which we have already
proved, the \emph{objects} of
the categories $\fl(\cs)\op$ and 
$\fl\big[\fs(\cs)\op\big]$ are identified,
up to isomorphism, by the functor $\wh Y$.
Start with objects $A',B'\in\fl\big[\fs(\cs)\op\big]$,
and choose objects $A,B\in\fl(\cs)\op$ and isomorphisms
$\wh Y(A)\cong A'$ and $\wh Y(B)\cong B'$.
It suffices to exhibit any morphism
$k:A'\la B'$,
in the category $\fl\big[\fs(\cs)\op\big]$,
as being isomorphic to
$\wh Y(\ph)$ for some $\ph:A\la B$,
with the isomorphisms as above.

Let $A',B'$ be objects in 
$\fl\big[\fs(\cs)\op\big]$; by
definition we can choose
Cauchy sequences $a'_*$ and $b'_*$ in
the category
$\fs(\cs)\op$, with $A'=\colim\,\wh Y(a'_*)$ and
with $B'=\colim\,\wh Y(b'_*)$.
Applying the
combination of Lemmas~\ref{L3.92929} and
\ref{L3.90909} allows us to find,
in the category $\fl(\cs)$, subsequences
$a_*$ of $a'_*$ and $b_*$ of $b'_*$, and objects
$A=\clim a_*$ and $B=\clim b_*$ with
$\wh Y(A)=\colim\,\wh Y(a_*)=A'$ and
$\wh Y(B)=\colim\,\wh Y(b_*)=B'$.
But
\cite[Lemma~\ref{L28.19}(i)]{Neeman18A}
tells us that any morphism
$k:\wh Y(A)\la\wh Y(B)$, in the category
$\MMod{\cs\op}$, is (up to replacing $a_*$ and
$b_*$ by subsequences) of the form
$\colim\,\wh Y(\ph_*)$, where $\ph_*:a_*\la b_*$ is
a map of Cauchy sequences in $\fs(\cs)\op$.

Translating from $\fs(\cs)\op$ to $\fs(\cs)$:
we are given \emph{inverse sequences} $a_*$ and
$b_*$ in $\fs(\cs)$, and a map of
inverse sequences $\ph\op_*:b_*\la a_*$.
Passing to the limit in the category
$\MMod\cs$ gives a morphism
$h=\clim\ph\op_*:\clim b_*\la\clim a_*$,
and as the functor
$\wh Y:\big[\MMod\cs\big]\op\la\MMod{\fs(\cs)\op}$
respects the colimits in question,
we have that
$\wh Y(h)=\wh Y(\clim\ph\op_*)=\colim\,\wh Y(\ph_*)$.
We have produced a map $h:B\la A$ such that
$\wh Y(h):\wh Y(A)\la\wh Y(B)$ agrees
with the given map
$k:A'\la B'$.
This completes the proof of (ii).

And finally for (iii). In
Corollary~\ref{C3.91919}(ii) we showed that the
functor
$\wh Y:\big[\fl(\cs)\big]\op\la\fl\big[\fs(\cs)\op\big]$
takes $\cl\op_\ell$ into $\wh\cl_\ell$.
To prove that the essential image
is all of $\wh\cl_\ell$ we appeal again to
the combination of 
Lemmas~\ref{L3.92929} and \ref{L3.90909}:
any object $B\in\wh\cl_\ell$ is the
colimit of a Cauchy sequence $b'_*$ contained in
$\cn\op_\ell\subset\fs(\cs)\op$. That is:
$b'_*$ is an inverse sequence in $\cn_\ell$.
Lemmas~\ref{L3.92929} and \ref{L3.90909}
allow us to construct a subsequence
$b_*$ of $b'_*$, and the object $F=\clim b_*$
can also be expressed in 
the form $F=\colim\,Y(f_*)$ for a
Cauchy sequence $f_*$ contained in $\cm_\ell$.
Hence $F$ belongs to $\cl_\ell$.
But Lemmas~\ref{L3.92929} and \ref{L3.90909} also tell us
that $\wh Y(F)=\colim\,\wh Y(b_*)=B$.
Thus $\wh Y(F)=B$ for $F\in\cl_\ell$, proving that
every object $B\in\wh\cl_\ell$ is in the essential
image of $\cl\op_\ell$ under the functor $\wh Y$.
\eprf

\pro{P3.94949}
Let $\cs$ be a triangulated category, and assume that 
$\{\cm_i\mid i\in\nn\}$
is an excellent metric on $\cs$.
Then the functor
$\wh Y:\big[\MMod\cs\big]\op\la\MMod{\big[\fs(\cs)\op\big]}$,
of Definition~\ref{D3.6.5},
restricts to an \emph{equivalence of the subcategories}
\[\xymatrix@C-10pt{
\wh Y\ar@{}[r]|-{:}&
\fl(\cs)\op\ar[rrr]^-\sim &&& \fl\big[\fs(\cs)\op\big]\ .
}\]
Moreover: this equivalence restricts, for every $\ell>0$,
to an equivalence
\[\xymatrix@C-10pt{
\wh Y\ar@{}[r]|-{:}&
\cl\op_\ell\ar[rrr]^-\sim &&& \wh\cl_\ell\ .
}\]
\epro

\prf
Corollary~\ref{C3.93939} tells us that
the restriction of
$\wh Y:\big[\MMod\cs\big]\op\la\MMod{\big[\fs(\cs)\op\big]}$
delivers a diagram
\[\xymatrix@C-10pt@R-20pt{
\wh Y\ar@{}[r]|-{:}&
\cl\op_\ell\ar[rrr]
\ar@{}[d]|-{\text{\rotatebox[origin=c]{-90}{$\ds\subset$}}}&&&
\wh\cl_\ell
\ar@{}[d]|-{\text{\rotatebox[origin=c]{-90}{$\ds\subset$}}}\\
\wh Y\ar@{}[r]|-{:}&
\fl(\cs)\op\ar[rrr] &&& \fl\big[\fs(\cs)\op\big]\ ,
}\]
and what we know so far is that the horizontal 
functors are full and essentially surjective. It remains to
prove the faithfulness.

Suppose therefore that $f:A\la B$ is a morphism
in $\fl(\cs)$ with $\wh Y(f)=0$; we need to show
that $f=0$. And because this is a statement
about the morphism $f:A\la B$ in the
category $\MMod\cs$, it suffices to
prove that, for every object $s\in\cs$ and for every
morphism $\s:Y(s)\la A$, the composite
$Y(s)\stackrel\s\la A\stackrel f\la B$ vanishes.

Because the metric $\{\cm_i,\,i\in\nn\}$ is
excellent, Definition~\ref{D3.3}(i) tells us that there must
exist an integer $n>0$ with $s\in{^\perp\cm_n}$. And
Lemma~\ref{L3.2.5}(i) shows
therefore that $Y(s)\in{^\perp\cl_n}$.
Now choose any Cauchy sequence $b_*\in\cs$ with
$B=\colim\,Y(b_*)$. Lemma~\ref{L3.5}
permits us to find a type-$(n+1)$ morphism
$B\la D$ with respect to Cauchy sequences $(b_*,d_*)$,
where we also know that $D\in\fs(\cs)$; the
more precise information about $D$ is irrelevant here.
Lemma~\ref{L3.309} establishes that
the morphism $B\la D$ extends to a strong
triangle $E\la B\la D\la\T E$ with $E$ and $\T E$ both
in $\cl_n$. And now consider the composite
\[\xymatrix{
Y(s) \ar[r]^-\s & A\ar[r]^-f & B\ar[r] & D\ .
}\]
We are given that $\wh Y(f)=0$, meaning that the
composite $A\stackrel f\la B\la D$ must vanish,
and
hence so does the longer composite displayed above.
We abbreviate this vanishing composite to
$Y(s)\stackrel{f\s}\la B\la D$.
By
Lemma~\ref{L3.303} the functor $\Hom\big(Y(s),-\big)$
takes the strong triangle $E\la B\la D\la\T E$ to a
long exact sequence. The vanishing
of $Y(s)\stackrel{f\s}\la B\la D$ means that $f\s:Y(s)\la B$
must factor as $Y(s)\la E\la B$, but as $Y(s)$ belongs
to $^\perp\cl_n$ and $E\in\cl_n$ the map $Y(s)\la E$ must
vanish. This completes the proof.
\eprf

\section{Strong triangles in $\fl(\cs)\op\cong\fl\big[\fs(\cs)\op\big]$}
\label{S77}

With the functor $\wh Y$ as defined in
Discussion~\ref{D3.6.5},
the main result in this section will be:

\pro{P77.1}
Let $\cs$ be a triangulated category, and assume that 
$\{\cm_i\mid i\in\nn\}$
is an excellent metric on $\cs$.
If $A\stackrel f\la B\stackrel g\la C\stackrel h\la\T A$ is a strong triangle
in $\fl(\cs)$, then
$\wh Y(A)\stackrel{\wh Y(f)}\longleftarrow\wh Y(B)\stackrel{\wh Y(g)}\longleftarrow\wh Y(C)\stackrel{\wh Y(h)}\longleftarrow\wh Y(\T A)$ is
a strong triangle in $\fl\big[\fs(\cs)\op\big]$.
\epro

\rmk{R77.3}
In view of Lemmas~\ref{L3.92929} and \ref{L3.90909},
the proof reduces to exhibiting the strong triangle
$A\stackrel f\la B\stackrel g\la C\stackrel h\la\T A$
as the inverse limit in $\MMod\cs$ of a Cauchy sequence
of distinguished triangles in $\fs(\cs)$. After
all: the functor $\wh Y$ takes inverse
limits, of Cauchy inverse sequences
in the category $\fs(\cs)\subset\MMod\cs$ with respect
to the metric $\{\cn_i,\,i\in\nn\}$,
to colimits in $\MMod{\fs(\cs)\op}$.

Since
the details are involved, we break up the proof
into a sequence of Lemmas. But before we launch into
the proof we note the following refinement,
which we will also prove:
\ermk

\rfn{R77.5}
\begin{em}
Assume throughout that we are given not only a strong
triangle $A\stackrel f\la B\stackrel g\la C\stackrel h\la\T A$
in the category $\fl(\cs)$, but also we choose in
advance
\be
\item
A Cauchy sequence of
distinguished triangles
$a_*\stackrel{f_*}\la b_*\stackrel{g_*}\la c_*\stackrel{h_*}\la\T a_*$ in the category $\cs$,
such that in the
commutative diagram below we have vertical isomorphisms
\[\xymatrix@C+30pt@R-10pt{
\colim\,Y(a_*)\ar@{=}[d]\ar[r]^-{\colim\,Y(f_*)}& \colim\,Y(b_*)\ar@{=}[d]
\ar[r]^-{\colim\,Y(g_*)}& \colim\,Y(c_*)\ar@{=}[d]
\ar[r]^-{\colim\,Y(h_*)}&\colim\,Y(\T a_*)\ar@{=}[d]\\
A\ar[r]_-f & B\ar[r]_-g & C\ar[r]_-h &\T A
}\]
\setcounter{enumiv}{\value{enumi}}
\ee
So far the assumptions are global,
they will hold in all the
statements in this section. Occasionally,
we might wish
to add one or both of the next two hypotheses:
\be
\setcounter{enumi}{\value{enumiv}}
\item
The object $B$ belongs to $\fs(\cs)$.
\item
There exists an integer $\ell>0$ such that
the Cauchy sequence $c_*$, in the
category $\cs$, satisfies
the hypothesis that $c_n\in\cm_\ell$
for all $n\gg0$.
\setcounter{enumiv}{\value{enumi}}
\ee
As we noted in Remark~\ref{R77.3},
the proof of
Proposition~\ref{P77.1} will amount
to exhibiting in $\fs(\cs)\op$ a
Cauchy sequence of triangles which, when we
view it as an inverse sequence
in $\fs(\cs)$ and write it as
$\wt a_*\stackrel{\wt f_*}\la \wt b_*\stackrel{\wt g_*}\la\wt c_*\stackrel{\wt h_*}\la\T\wt a_*$,
delivers in $\MMod\cs$ isomorphisms
\[\xymatrix@C+40pt@R-10pt{
A\ar@{=}[d]\ar[r]^-f & B\ar@{=}[d]\ar[r]^-g & C\ar@{=}[d]\ar[r]^-h &\T A\ar@{=}[d]\\
\clim \wt a_*\ar[r]_-{\clim \wt f_*}& \clim \wt b_*\ar[r]_-{\clim \wt g_*}& \clim \wt c_*\ar[r]_-{\clim \wt h_*}&\clim \T\wt a_*
}\]
The refinement we assert is that, if we 
assume (ii) and/or (iii) above, then the
Cauchy sequence $\wt a_*\stackrel{\wt f_*}\la \wt b_*\stackrel{\wt g_*}\la\wt c_*\stackrel{\wt h_*}\la\T\wt a_*$,
of
triangles in $\fs(\cs)\op$, can be chosen so that:
\be
\setcounter{enumi}{\value{enumiv}}
\item
Assume (ii) holds. Then,
after a finite initial segment, the sequence 
$\wt b_*$ becomes the constant
sequence $\cdots\la B\stackrel\id\la B\stackrel\id\la B\stackrel\id\la B$.
\item
Assume (iii) holds. Then,
for the integer $\ell>0$ of (iii),
the Cauchy sequence $\wt c_*$
satisfies
the hypothesis that $\wt c_n\in\cn_\ell$
for all $n\gg0$.
\ee
\end{em}
\erfn

\lem{L77.7}
With the notation as in Refinement~\ref{R77.5}(i),
we can, for
every integer
$i>0$,
\be
\item
Construct in $\fl(\cs)$ a commutative square
\[\xymatrix@C+10pt{
A\ar[r]^-f\ar[d]_{\alpha_i} &B\ar[d]^{\beta_i}\\
\wt a_i\ar[r]^-{\wt f_i} &\wt b_i
}\]
\item
We also construct Cauchy sequences
$\wt a_{i,*}$ and $\wt b_{i,*}$,
in the category $\cs$, with
$\wt a_i=\colim\,Y(\wt a_{i,*})$ and
$\wt b_i=\colim\,Y(\wt b_{i,*})$.
\setcounter{enumiv}{\value{enumi}}
\ee
For the next assertions, let the
integers $n_i$ be
as chosen in Lemma~\ref{L3.7}.
\be
\setcounter{enumi}{\value{enumiv}}
\item
The object
$\wt a_i$ in the commutative
square of \emph{(i)}
belongs to $\fs(\cs)\cap\cl_{n_{i+2}^{}}^\perp$,
and the morphism $\alpha_i:A\la \wt a_i$  in the commutative
square of \emph{(i)} is
of type-$(1+n_{i+1}^{})$ with respect to
to the Cauchy sequences $(a_*,\wt a_{i,*})$.
\item
The object
$\wt b_i$  in the commutative
square of \emph{(i)} belongs
to $\fs(\cs)\cap\cl_{n_{i+1}^{}}^\perp$,
and the morphism $\beta_i:B\la \wt b_i$
in the commutative
square of \emph{(i)}
is
of type-$(1+n_i)$ with respect to $(b_*,\wt b_{i,*})$.
\item
If the object $B$ belongs to $\fs(\cs)$, then
for $i\gg0$ 
the commutative
square of \emph{(i)} may be chosen to be
\[\xymatrix@C+10pt{
A\ar[r]^-f\ar[d]_{\alpha_i} &B\ar@{=}[d]\\
\wt a_i\ar[r]^-{\wt f_i} & B
}\]
and the Cauchy sequence $\wt b_{i,*}$ of \emph{(ii)}
may be chosen to be
equal to $b_*$.
\setcounter{enumiv}{\value{enumi}}
\ee
\elem

\prf
Let us begin with (v). If $B$ belongs to $\fs(\cs)$,
then there exists an integer $N>0$ with $B\in\cl_N^\perp$.
Therefore, if $i$ is large enough, then setting
$\wt b_i=B$, letting $\beta_i:B\la B$ be the identity,
and choosing $\wt\beta_{i,*}$ to be the
sequence $b_*$ certainly meets all the requirements
placed on $\beta_i:B\la\wt b_i$ in (ii) and (iii).

More generally,
since $A$ and $B$ are objects in $\fl(\cs)$ we may apply to
each of them 
Lemma~\ref{L3.7}(i).
For every
integer $i\geq1$ this gives:
\be
\setcounter{enumi}{\value{enumiv}}
\item
An object
$\wt a_i\in\fs(\cs)\cap\cl_{n_{i+2}^{}}^\perp$,
a Cauchy sequence $\wt a_{i,*}$ with
$\wt a_i=\colim\,Y(\wt a_{i,*})$, 
and a morphism $\alpha_i:A\la \wt a_i$ is
of type-$(1+n_{i+1}^{})$ with respect to $(a_*,\wt a_{i,*})$.
\item
An object
$\wt b_i\in\fs(\cs)\cap\cl_{n_{i+1}^{}}^\perp$,
a Cauchy sequence $\wt b_{i,*}$ with
$\wt b_i=\colim\,Y(\wt b_{i,*})$,
and a morphism $\beta_i:B\la \wt b_i$ is
of type-$(1+n_i)$ with respect to $(b_*,\wt b_{i,*})$.
\setcounter{enumiv}{\value{enumi}}
\ee
Thus (ii), (iii), (iv) and (v) pose
no serious problem. The issue is that,
until now, all we have done is  construct the diagram
\[\xymatrix@C+10pt{
A\ar[r]^-f\ar[d]_{\alpha_i} &B\ar[d]^{\beta_i}\\
\wt a_i &\wt b_i
}\]
We need to show that there exists an
$\wt f_i:\wt a_i\la\wt b_i$ completing it to
the commutative square in (i).

Since $\alpha_i:A\la\wt a_1^{}$ is a
type-$(1+n_{i+1}^{})$ morphism
with respect to $(a_*,\wt a_{i,*})$,
while the object $\wt b_i$ belongs to
$\fl(\cs)\cap\cl_{n_{i+1}^{}}^\perp$,
Corollary~\ref{C3.378} tells us that the
$\Hom(-,\wt b_i)$ takes the
morphism $\alpha_i:A\la\wt a_i$ to an isomorphism. Therefore
there is an (even unique) map $\wt f_i:\wt a_i\la\wt b_i$
rendering commutative the square
\[\xymatrix@C+10pt{
A\ar[r]^-f\ar[d]_{\alpha_i} &B\ar[d]^{\beta_i}\\
\wt a_i\ar[r]^-{\wt f_i} &\wt b_i
}\]
This completes the proof of the Lemma.
\eprf

\dis{D77.9}
For each $i>0$, Lemma~\ref{L77.7} produced
for us (among other things) a morphism
$\wt f_i:\wt a_i\la\wt b_i$ in the category
$\fs(\cs)$.
This morphism $\wt f_i:\wt a_i\la\wt b_i$ may be completed,
in the triangulated category $\fs(\cs)$, to
a distinguished triangle. We obtain in $\fl(\cs)$
the diagram
\[\xymatrix@C+10pt{
A\ar[r]^-{f}\ar[d]^{\alpha_i} &
B \ar[d]^{\beta_i}\ar[r]^-{g} &
C\ar[r]^-{h} &
\T A
\\
\wt a_i\ar[r]^-{\wt f_i} &\wt b_i\ar[r]^-{\wt g_i^{}} &
\wt c_i\ar[r]^-{\wt h_i} &\wt \T a_i &
}\]
where the top row is the strong triangle
we started out with
in Proposition~\ref{P77.1}, the bottom
row is a distinguished triangle
in $\fs(\cs)$, and the commutative square on the
left is as given in Lemma~\ref{L77.7}(i).
And the objective of the rest of the
proof of Proposition~\ref{P77.1},
and of its Refinement~\ref{R77.5}, is to show
that
\be
\item
We may choose a morphism $\gamma_i:C\la\wt c_i$
rendering commutative the diagram
\[\xymatrix@C+10pt{
A\ar[r]^-{f}\ar[d]^{\alpha_i} &
B \ar[d]^{\beta_i}\ar[r]^-{g} &
C\ar[d]^{\gamma_i}\ar[r]^-{h} &
\T A\ar[d]^{\T\alpha_i}
\\
\wt a_i\ar[r]^-{\wt f_i} &\wt b_i\ar[r]^-{\wt g_i^{}} &
\wt c_i\ar[r]^-{\wt h_i} &\wt \T a_i &
}\]
\item
We can do it in such a way that the
distinguished triangles
$\wt a_i\stackrel{\wt f_i}\la \wt b_i\stackrel{\wt g_i}\la\wt c_i\stackrel{\wt h_i}\la\T\wt a_i$,
in the category
$\fs(\cs)$, assemble to an inverse sequences, this
sequence is Cauchy, and the inverse limit is
the given strong triangle
$A\stackrel f\la B\stackrel g\la C\stackrel h\la\T A$.
\ee
\edis

In Discussion~\ref{D77.9} we sketched our long-term
strategy, mapping out the goals of this section.
Let us step back to where we are right now.

\con{C77.11}
So far we have the following diagram in
the category $\fl(\cs)$
\[\xymatrix@C+10pt{
A\ar[r]^-{f}\ar[d]^{\alpha_i} &
B \ar[d]^{\beta_i}\ar[r]^-{g} &
C\ar[r]^-{h} &
\T A
\\
\wt a_i\ar[r]^-{\wt f_i} &\wt b_i & &
}\]
where the top row is the strong triangle 
fixed for the proof of Proposition~\ref{P77.1},
and
where the square is the commutative square of
Lemma~\ref{L77.7}. But we are also given
Cauchy sequences
\[
a_*,\quad \wt a_{i,*},\quad
b_*,\quad \wt b_{i,*},\quad
c_*
\]
and isomorphisms
\[\xymatrix@C+0pt@R-10pt{
\colim\,Y(a_*) \ar@{=}[d] &
\colim\,Y(\wt a_{i,*}) \ar@{=}[d] &
\colim\,Y(b_*) \ar@{=}[d] &
\colim\,Y(\wt b_{i,*}) \ar@{=}[d] &
\colim\,Y(c_*) \ar@{=}[d] \\
A &\wt a_i &B &\wt b_i & C
}\]
And, in this construction, we will pass to
subsequences that facilitate the rest of the
argument. 

First of all, by
\cite[Lemma~\ref{L28.19}(i)]{Neeman18A}
we may, after replacing
$\wt a_{i,*}$ and $\wt b_{i,*}$ by subsequences, choose
maps of sequences $\alpha_{i.*}:a_*\la\wt a_{i,*}$
and $\beta_{i.*}:b_*\la\wt \beta_{i,*}$
with $\alpha_i=\colim\,\alpha_{i,*}$ and
with $\beta_i=\colim\,\beta_{i,*}$.
Now, the fact that $\alpha_i:A\la \wt a_i$
of type-$(1+n_{i+1}^{})$ with respect to $(a_*,\wt a_{i,*})$
and that
$\beta_i:B\la \wt b_i$
of type-$(1+n_i)$ with respect to $(b_*,\wt b_{i,*})$
tells us that there exists an integer $N\gg0$ such that
\be
\item
If $\ell>N$ is any integer,
then the triangles
\[
a_\ell\la \wt a_{i,\ell}\la \wh d_\ell 
\qquad{and}\qquad
b_\ell\la \wt b_{i,\ell}\la d_\ell 
\]
have $\wh d_\ell\in\cm_{1+n_{i+1}^{}}$ and
$d_\ell\in\cm_{1+n_i}$.
\setcounter{enumiv}{\value{enumi}}
\ee
Increasing $N$ even more, the fact that the sequences
$a_*$, $\wt a_{i,*}$, $b_*$, $\wt b_{i,*}$
and $c_*$ are all Cauchy allows us
to assume that 
\be
\setcounter{enumi}{\value{enumiv}}
\item
If $\ell>k$ are any integers $\geq N$,
then the five triangles
\[
\begin{array}{ccccc}
a_k^{}\la a_\ell^{}\la\delta_{k,\ell},& \qquad&
b_k^{}\la b_\ell^{}\la\wh\delta_{k,\ell}, &\quad &
c_k^{}\la c_\ell^{}\la\ov\delta_{k,\ell},\\
\wt a_{i,k}^{}\la \wt a_{i,\ell}^{}\la\wt\delta_{i,k,\ell},&\qquad&
\wt b_{i,k}^{}\la \wt b_{i,\ell}^{}\la\wt\delta'_{i,k,\ell}& &
\end{array}
\]
are such that
\[
\Tm\delta_{k,\ell},\ \delta_{k,\ell},\ 
\Tm\wh\delta_{k,\ell},\ \wh\delta_{k,\ell},\
\Tm\ov\delta_{k,\ell},\ 
\ov\delta_{k,\ell},\
\Tm\wt\delta_{i,k,\ell},\ 
\wt\delta_{i,k,\ell},\ 
\Tm\wt\delta'_{i,k,\ell},\ 
\wt\delta'_{i,k,\ell}
\]
all lie in
$\cm_{1+n_{i+2}^{}}$.
\setcounter{enumiv}{\value{enumi}}
\ee
Replacing each of the Cauchy sequences $x_*$ above,
with the subsequence $y_*=x_{(*+N-1)}$,
we may assume $N=1$.

Now: we have a map $\wt f_i:\wt a_i\la \wt b_i$ in the
category $\fs(\cs)\subset\fl(\cs)$, and we know that
$\wt a_i=\colim\,Y(\wt a_{i,*})$ and
that $\wt b_i=\colim\,Y(\wt b_{i,*})$.
Therefore the composite
$Y(\wt a_{i,1}^{})\la\wt a_i\stackrel{\wt f_i}\la\wt b_i$ must factor
through some $Y(\wt b_{i,\ell}^{})$; that is, we obtain
a morphism $\ph:\wt a_{i,1}^{}\la\wt b_{i,\ell}^{}$
rendering commutative
the square
\[\xymatrix@C+10pt{
Y(\wt a_{i,1}^{})\ar[r]^-{Y(\ph)}\ar[d] &Y(\wt b_{i,\ell}^{})\ar[d]\\
\wt a_i\ar[r]^-{\wt f_i} &\wt b_i
}\]
The diagram
\[\xymatrix@C+10pt{
  a_{1}^{}\ar[r]^-{f_{1}}\ar[d]_-{\alpha_{i,1}^{}} &b_{1}\ar[r] &b_\ell\ar[d]^-{\beta_{i,\ell}^{}}\\
  \wt a_{i,1}^{}
\ar[rr]^-\ph &&\wt b_{i,\ell}
}\]
is under no obligation to commute in $\cs$, but the map
$Y(\wt b_{i,\ell})\la\colim\,Y(\wt b_{i,*})=\wt b_i$
coequalizes the
two composites. Hence, at the cost of increasing $\ell$,
we do obtain a commutative diagram. Now we
replace the sequence
$\wt b_{i,*}$ by the subsequence $\wt b_{i,(*+\ell-1)}$
and the map of sequences $\beta_{i,*}^{}:b_*\la \wt b_{i,*}$ by
the composite
\[\xymatrix@C+10pt{
b_*\ar[r]^-{\beta_{i,*}^{}} &\wt b_{i,*}\ar[r] &\wt b_{i,(*+\ell-1)}\ .
}\]
We note that, for each integer $j>0$,
the triangle
$b_j\la \wt b_{i,j}\la d_j $ has $d_j\in\cm_{1+n_i}$,
while the triangle
$\wt b_{i,j}\la \wt b_{i,(j+\ell-1)}\la\wt\delta'_{i,j,(j+\ell-1)}$
has $\wt\delta'_{i,j,(j+\ell-1)}\in\cm_{1+n_{i+2}^{}}$.
Thus the completion of the composite
to a triangle 
$b_j\la \wt b_{i,(j+\ell-1)}\la d'_j$
has $d'_j\in\cm_{1+n_i}$.

We have done a great deal of replacing of
Cauchy sequences by subsequences,
but note that
\be
\setcounter{enumi}{\value{enumiv}}
\item
We still have a
Cauchy sequence
$a_*\stackrel{f_*}\la b_*\stackrel{g_*}\la c_*\stackrel{h_*}\la\T a_*$ of distinguished
triangles in $\cs$, such that
the colimit of its image under $Y$ is the given
strong triangle 
$A\stackrel f\la B\stackrel g\la C\stackrel h\la\T A$.
\setcounter{enumiv}{\value{enumi}}
\ee
This is because the sequences $a_*$, $b_*$ and $c_*$ were
replaced only by the straightforward, truncated
subsequences $a_{*+N-1}^{}$, $b_{*+N-1}^{}$ and $c_{*+N-1}^{}$.
The replacing of $\wt a_{i,*}$ and $\wt b_{i,*}$ by
subsequences was more complicated, but the sum total
is that
\be
\setcounter{enumi}{\value{enumiv}}
\item
We have maps of Cauchy sequences
$\alpha_{i,*}:a_*\la\wt a_{i,*}$ and
$\beta_{i,*}:b_*\la\wt b_{i,*}$
such that the colimit of the image
under $Y$ of these Cauchy sequences
are, respectively, $\alpha_i:A\la\wt a_i$ and
$\beta_i:B\la\wt b_i$.
\item
The hypotheses
of (i) hold for any integer $\ell\geq1$ (and not just
for $\ell\gg0$), and similarly 
the hypotheses of (ii) hold for any pair of
integers $\ell>k\geq1$.
\setcounter{enumiv}{\value{enumi}}
\ee
And furthermore we have
\be
\setcounter{enumi}{\value{enumiv}}
\item
A morphism $\ph:\wt a_{i,1}^{}\la \wt b_{i,1}$,
in the category $\cs$, such that
the two squares below both commute
\[\xymatrix@C+10pt{
a_1^{}\ar[r]^-{f_1}\ar[d]_{\alpha_{i,1}^{}} &b_1\ar[d]^{\beta_{i,1}^{}}
&  & Y(\wt a_{i,1}^{})\ar[r]^-{Y(\ph)}\ar[d]
&Y(\wt b_{i,1}^{})\ar[d]\\
\wt a_{i,1}^{}\ar[r]^-{\ph} &\wt b_{i,1}^{} & &
\wt a_i\ar[r]^-{\wt f_i} &\wt b_i
}\]
The square on the left commutes in $\cs$, while the
square on the right commutes in $\fl(\cs)$.
\setcounter{enumiv}{\value{enumi}}
\ee
Next let us focus on
the right-hand commutative square in (vi).
We have a morphism
$\wt f_i:\wt a_i\la\wt b_i$ in the category
$\fl(\cs)$, as well as Cauchy sequences
$\wt a_{i,*}$ and $\wt b_{i,*}$ with
$\wt a_i=\colim\,Y(\wt a_{i,*})$ and
$\wt b_i=\colim\,Y(\wt b_{i,*})$, and moreover
a commutative square
\[\xymatrix@C+10pt{
Y(\wt a_{i,1}^{})\ar[r]^-{Y(\ph)}\ar[d]
&Y(\wt b_{i,1}^{})\ar[d]\\
\wt a_i\ar[r]^-{\wt f_i} &\wt b_i
}\]
and \cite[Lemma~\ref{L28.19}(i)]{Neeman18A}
permits us to
\be
\setcounter{enumi}{\value{enumiv}}
\item
Choose
a Cauchy subsequence $\wt b'_{i,*}$ of $\wt b_{i,*}$,
with $\wt b'_{i,1}=\wt b_{i,1}^{}$, and
a morphism of Cauchy sequences
$\wt f_{i,*}:\wt a_{i,*}\la \wt b'_{i,*}$, such that
$\wt f_{i,1}=\ph$ with the morphism
$\ph$ as in (vi), and where $\wt f_i=\colim\,\wt f_{i,*}$.
\setcounter{enumiv}{\value{enumi}}
\ee
Now: in (ii), and with $N=1$ by our
choice of subsequences as remarked in (v),
for any $\ell>k\geq1$
the triangles
\[
\wt a_{i,k}^{}\la\wt a_{i,\ell}^{}\la\wt\delta_{i,k,\ell}
\qquad\text{ and }\qquad
\wt b'_{i,k}\la\wt b'_{i,\ell}\la\wt\delta'_{i,k,\ell}
\]
satisfy
\[
\{\Tm\wt\delta_{i,k,\ell},\ \wt\delta_{i,k,\ell},\
\Tm\wt\delta'_{i,k,\ell},\ \wt\delta'_{i,k,\ell}\}
\sub\cm_{1+n_{i+2}^{}}\ .
\]
Hence we see that
\[
\{\wt\delta_{i,k,\ell},\ \T\wt\delta_{i,k,\ell},\
\Tm\wt\delta'_{i,k,\ell},\ \wt\delta'_{i,k,\ell}\}
\sub\cm_{n_{i+2}^{}}\ .
\]
Now we apply 
\cite[Lemma~\ref{L28.19}~(ii) and (iii)]{Neeman18A}.
This
allows us to
\be
\setcounter{enumi}{\value{enumiv}}
\item
Extend the morphism of Cauchy sequences
$\wt f_{i,*}:\wt a_{i,*}\la \wt b'_{i,*}$,
obtained in (vii),
to a Cauchy sequence of triangles
$\wt a_{i,*}\stackrel{\wt f_{i,*}}\la \wt b'_{i,*}\stackrel{\wt g_{i,*}^{}}\la \wt c_{i,*}\stackrel{\wt h_{i,*}^{}}\la \T \wt a_{i,*}$
where we know, for any $\ell>k\geq1$, that the
triangle $\wt c_{i,k}^{}\la \wt c_{i,\ell}^{}\la\wt\delta''_{i,k,\ell}$ satisfies
$\{\Tm\wt\delta''_{i,k,\ell},\,\wt\delta''_{i,k,\ell}\}\subset\cm_{n_{i+2}^{}}$.
\setcounter{enumiv}{\value{enumi}}
\ee
Next, applying the functor $Y$ to the Cauchy sequence
of triangles 
$\wt a_{i,*}\stackrel{\wt f_{i,*}}\la \wt b'_{i,*}\stackrel{\wt g_{i,*}^{}}\la \wt c_{i,*}\stackrel{\wt h_{i,*}}\la \wt a_{i,*}$ that we obtained in
(viii), and taking colimits,
yields a strong triangle
$\wt a_i\stackrel{\wt f_i}\la \wt b_i\stackrel{\wt g_i^{}}\la\wt c_i\stackrel{\wt h_i}\la\T\wt a_i$
in the category $\fl(\cs)$. Since the objects
$\wt a_i$ and $\wt b_i$ both lie in $\fs(\cs)$,
from \cite[Remark~\ref{R200976}]{Neeman18A}
we learn that
this strong triangle is a distinguished
triangle in $\fs(\cs)$. More precisely: by
Lemma~\ref{L77.7}~(ii) and (iii) we
have that 
$\wt b_i\in\cl_{n_{i+1}^{}}^\perp\subset\T\cl_{n_{i+2}^{}}^\perp$ and
$\wt a_i\in\cl_{n_{i+2}^{}}^\perp$, and
the pre-triangle $\wt b_i\la\wt c_i\la\T\wt a_i$ gives that
$\wt c_i\in\T\cl_{n_{i+2}^{}}^\perp$. Summarizing,
we have produced
\be
\setcounter{enumi}{\value{enumiv}}
\item
A commutative diagram in $\fl(\cs)$
\[\xymatrix@C+10pt{
Y(\wt a_{i,1}^{})\ar[r]^-{Y(\wt f_{i,1})}\ar[d] &
Y(\wt b_{i,1}^{})\ar[d]\ar[r]^-{Y(\wt g_{i,1}^{})} &
Y(\wt c_{i,1}^{})\ar[d]\ar[r]^-{Y(\wt h_{i,1})} &
Y(\T\wt a_{i,1}^{})\ar[d]
\\
\wt a_i\ar[r]^-{\wt f_i} &\wt b_i\ar[r]^-{\wt g_i} &
\wt c_i\ar[r]^-{\wt h_i} &\wt \T a_i &
}\]
where the top row is the image
under $Y$ of the 
distinguished triangle
$\wt a_{i,1}^{}\stackrel{\wt f_{i,1}}\la \wt b'_{i,1}\stackrel{\wt g_{i,1}^{}}\la \wt c^{}_{i,1}\stackrel{\wt h_{i,1}}\la \wt a_{i,1}^{}$ in the category $\cs$, the bottom row
is a distinguished triangle in $\fs(\cs)$, and
where every object in the bottom
row belongs to $\cl_{n_{i+2}^{}}^\perp\cup\T\cl_{n_{i+2}^{}}^\perp\subset\cl_{1+n_{i+2}^{}}^\perp$. 
\setcounter{enumiv}{\value{enumi}}
\ee
Now, in  all of the work that we have done
since (vi), and which all focused on the
right-hand commutative square in (vi), the
left-hand commutative square of (vi) hasn't disappeared.
It can be rewritten as
\be
\setcounter{enumi}{\value{enumiv}}
\item
In the diagram below the rows are distinguished
triangles in $\cs$, and the square commutes:
\[\xymatrix@C+10pt{
a_1^{}\ar[r]^-{f_1}\ar[d]_{\alpha_{i,1}^{}} &b_1\ar[d]^{\beta_{i,1}^{}}
\ar[r]^-{g^{}_1}
&c_1^{}\ar[r]^-{h_1} &\T a_1 \\
\wt a_{i,1}^{}\ar[r]^-{\wt f_{i,1}} &\wt b_{i,1}^{}\ar[r]^-{\wt g_{i,1}^{}} &
\wt c_{i,1}^{}\ar[r]^-{\wt h_{i,1}}& \T \wt a_{i,1}^{}
}\]
\setcounter{enumiv}{\value{enumi}}
\ee
The reason for the commutativity of the square
is that the map $\wt f_{i,1}:\wt a_{i,1}^{}\la\wt b_{i,1}$ was
chosen, back in (vii),
to agree with $\ph:\wt a_{i,1}^{}\la\wt b_{i,1}$. Hence
the commutativity is just a rewriting of the
left-hand square of (vi).
\econ

Before stating the next Lemma we summarize the
highlights of what we have done so far.

\smr{S77.13}
We were given, 
back in Refinement~\ref{R77.5}(i), a Cauchy sequence
of distinguished triangles
$a_*\stackrel{f_*}\la b_*\stackrel{g_*}\la c_*\stackrel{h_*}\la\T a_*$, in the category $\cs$,
whose image under the functor $Y$ has, as its colimit,
the strong triangle
$A\stackrel f\la B\stackrel g\la C\stackrel h\la\T A$
that we fixed for the purpose of proving both
Proposition~\ref{P77.1} and 
Refinement~\ref{R77.5}.
And in Construction~\ref{C77.11}(viii)
we formed, in $\cs$, a Cauchy sequence of
triangles
$\wt a_{i,*}\stackrel{\wt f_{i,*}}\la \wt b'_{i,*}\stackrel{\wt g_{i,*}^{}}\la \wt c_{i,*}\stackrel{\wt h_{i,*}}\la \T\wt a_{i,*}$.
Applying the functor $Y$
and taking colimits, in the category $\MMod\cs$,
yielded the distinguished triangle
$\wt a_i\stackrel{\wt f_i}\la \wt b_i\stackrel{\wt g_i^{}}\la\wt c_i\stackrel{\wt h_i}\la\T\wt a_i$
in the category $\fs(\cs)$.
And combining this with the commutative square
of Lemma~\ref{L77.7} gives the diagram
\[\xymatrix@C+10pt{
A\ar[r]^-{f}\ar[d]^{\alpha_i} &
B \ar[d]^{\beta_i}\ar[r]^-{g} &
C\ar[r]^-{h} &
\T A
\\
\wt a_i\ar[r]^-{\wt f_i} &\wt b_i\ar[r]^-{\wt g_i^{}} &
\wt c_i\ar[r]^-{\wt h_i} &\wt \T a_i &
}\]
that we already met at the
beginning of Discussion~\ref{D77.9}.
\esmr

Summary~\ref{S77.13} has
prepared us to \emph{state} the next Lemma. The more
detailed information, contained in the highlighted
portions of Construction~\ref{C77.11},
will enter into the proof.

\lem{L77.15}
With the commutative diagram below be as
in Summary~\ref{S77.13}
\[\xymatrix@C+10pt{
A\ar[r]^-{f}\ar[d]^{\alpha_i} &
B \ar[d]^{\beta_i}\ar[r]^-{g} &
C\ar[r]^-{h} &
\T A
\\
\wt a_i\ar[r]^-{\wt f_i} &\wt b_i\ar[r]^-{\wt g_i^{}} &
\wt c_i\ar[r]^-{\wt h_i} &\wt \T a_i &
}\]
we may choose, in $\fl(\cs)$,
a morphism $\gamma_i:C\la\wt c_i$
rendering commutative the diagram
\[\xymatrix@C+10pt{
A\ar[r]^-{f}\ar[d]^{\alpha_i} &
B \ar[d]^{\beta_i}\ar[r]^-{g} &
C\ar[d]^{\gamma_i}\ar[r]^-{h} &
\T A\ar[d]^{\T\alpha_i}
\\
\wt a_i\ar[r]^-{\wt f_i} &\wt b_i\ar[r]^-{\wt g_i^{}} &
\wt c_i\ar[r]^-{\wt h_i} &\wt \T a_i &
}\]
And the properties that will be crucial,
in the rest of the proof of Proposition~\ref{P77.1}
and Refinement~\ref{R77.5}, are:
\be
\item
Every object in the bottom row belongs to
the category $\fs(\cs)\cap\cl_{1+n_{i+2}^{}}^\perp$.  
\item
All three morphisms
\[
\alpha_i:A\la\wt a_i\,,\qquad
\beta_i:B\la\wt b_i\,,\qquad
\gamma_i:C\la\wt c_i
\]
are of type-$(1+n_i)$, with respect to [respectively]
the pairs of Cauchy sequences
\[
(a_*,\wt a_{i,*}),\qquad
(b_*,\wt b_{i,*}),\qquad
(c_*,\wt c_{i,*})\ .
\]
\item
If we are given an integer $\ell>0$,
and assume that $c_n\in\cm_\ell$ for
all $n\gg0$, then there exists
an integer $N>0$ such that, for all $i>N$, we have
$\wt c_{i,n}\in\cm_\ell$ for $n\gg0$.
Therefore the $\wt c_i=\colim\,Y(\wt c_{i,*})$
must belong to $\cn_\ell$ for all $i>N$.
\setcounter{enumiv}{\value{enumi}}
\ee
\elem

\prf
Part (i) of the Lemma we already know, it was proved
in Construction~\ref{C77.11}(ix).
It's time to set to work, and prove the
existence of the morphism $\gamma_i$ and that
the resulting commutative
diagram satisfies properties (ii) and (iii)
of the Lemma.

In Construction~\ref{C77.11}(x)
we produced the diagram below 
\[\xymatrix@C+10pt{
a_1^{}\ar[r]^-{f_1}\ar[d]_{\alpha_{i,1}^{}} &b_1\ar[d]^{\beta_{i,1}^{}}
\ar[r]^-{g^{}_1}
&c_1^{}\ar[r]^-{h_1} &\T a_1 \\
\wt a_{i,1}^{}\ar[r]^-{\wt f_{i,1}} &\wt b_{i,1}^{}\ar[r]^-{\wt g_{i,1}^{}} &
\wt c_{i,1}^{}\ar[r]^-{\wt h_{i,1}}& \T \wt a_{i,1}^{}
}\]
The rows are distinguished
triangles in $\cs$, and the square commutes.
We may complete this commutative square to a
$3\times 3$ diagram where the rows and columns are
triangles
\[\xymatrix@C+10pt{
a_1^{}\ar[r]^-{f_1}\ar[d]_{\alpha_{i,1}^{}} &b_1\ar[d]^{\beta_{i,1}^{}}
\ar[r]^-{g^{}_1}
&c_1^{}\ar[d]^{\gamma_{i,1}^{}}\ar[r]^-{h_1} &\T a_1\ar[d]^{\T\alpha_{i,1}^{}} \\
\wt a_{i,1}^{}\ar[r]^-{\wt f_{i,1}}\ar[d] &\wt b_{i,1}^{}\ar[r]^-{\wt g_{i,1}^{}}\ar[d] &
\wt c_{i,1}^{}\ar[r]^-{\wt h_{i,1}}\ar[d]& \T \wt a_{i,1}^{}\ar[d] \\
\wh\delta_1\ar[r] &
\delta_1\ar[r] &
\ov\delta_1\ar[r] &
\T\wh\delta_1
}\]
and Construction~\ref{C77.11}(i), as refined in
Construction~\ref{C77.11}(v) after passing to subsequences,
tells us that $\delta_1\in\cm_{1+n_i}$ and
$\wh\delta_1\in\cm_{1+n_{i+1}^{}}$. Now $n_{i+1}^{}>n_i$
implies $\T\cm_{1+n_{i+1}^{}}\subset\cm_{1+n_i}$,
and therefore
$\T\wh\delta_1\in\cm_{1+n_i}$. The triangle
$\delta_1\la\ov\delta_1\la\T\delta_1$ now gives that
$\ov\delta_1\in\cm_{1+n_i}*\cm_{1+n_i}\subset\cm_{1+n_i}$.
And what we want to remember, for future
reference, is that
\be
\setcounter{enumi}{\value{enumiv}}
\item
We have produced, in the category $\cs$, a morphism
of triangles
\[\xymatrix@C+10pt{
a_1^{}\ar[r]^-{f_1}\ar[d]_{\alpha_{i,1}^{}} &b_1\ar[d]^{\beta_{i,1}^{}}
\ar[r]^-{g_1^{}}
&c_1^{}\ar[d]^{\gamma_{i,1}^{}}\ar[r]^-{h_1} &\T a_1\ar[d]^{\T\alpha_{i,1}^{}} \\
\wt a_{i,1}^{}\ar[r]^-{\wt f_{i,1}} &\wt b_{i,1}^{}\ar[r] &
\wt c_{i,1}^{}\ar[r]& \T \wt a_{i,1}^{} 
}\]
such that, in the triangle
$c_1^{}\stackrel{\gamma_{i,1}^{}}\la \wt c_{i,1}^{}\la\ov\delta_1$,
we have $\ov\delta_1\in\cm_{1+n_i}$.
\setcounter{enumiv}{\value{enumi}}
\ee
And now we combine (iv) above with
Construction~\ref{C77.11}(ix), to produce in $\fl(\cs)$ the
composite commutative diagram
\[\xymatrix@C+10pt{
&Y(a_1^{})\ar[r]^-{Y(f_1)}\ar[d]_{Y(\alpha_{i,1}^{})} &Y(b_1)\ar[d]^{Y(\beta_{i,1}^{})}
\ar[r]^-{Y(g^{}_1)}
&Y(c_1^{})\ar[d]^{Y(\gamma_{i,1}^{})}\ar[r]^-{Y(h_1)} &Y(\T a_1)\ar[d]^{Y(\T\alpha_{i,1}^{})} \\
(\dagger)&Y(\wt a_{i,1}^{})\ar[r]^-{Y(\wt f_{i,1})}\ar[d] &
Y(\wt b_{i,1}^{})\ar[d]\ar[r] &
Y(\wt c_{i,1}^{})\ar[d]\ar[r] &
Y(\T\wt a_{i,1}^{})\ar[d]
\\
&\wt a_i\ar[r]^-{\wt f_i} &\wt b_i\ar[r] &
\wt c_i\ar[r] &\wt \T a_i &
}\]
Now by Construction~\ref{C77.11}(ii), made
more precise after passing to subsequences
as in Construction~\ref{C77.11}(v), the Cauchy
sequences $a_*$, $b_*$ and $c_*$
all have the property that, for all pairs
of integers $\ell>k\geq1$,
the distinguished triangles
\[
\begin{array}{ccccc}
a_k^{}\la a_\ell^{}\la\delta_{k,\ell},& \qquad&
b_k^{}\la b_\ell^{}\la\wh\delta_{k,\ell}, &\quad &
c_k^{}\la c_\ell^{}\la\ov\delta_{k,\ell}
\end{array}
\]
are such that
\[
\Tm\delta_{k,\ell},\ \  \delta_{k,\ell},\ \  
\Tm\wh\delta_{k,\ell},\ \ \wh\delta_{k,\ell},\ \
\Tm\ov\delta_{k,\ell},\ \
\ov\delta_{k,\ell}
\]
all lie in
$\cm_{1+n_{i+2}^{}}$.
On the other hand we know from
Construction~\ref{C77.11}(ix) that the objects
$\wt a_i$, $\wt b_i$, $\wt c_i$ and $\T\wt a_i$ all lie
in $\fs(\cs)\cap\cl_{1+n_{i+2}^{}}^\perp$.
Now
\cite[Corollary~\ref{C28.105}]{Neeman18A}
permits us to factor the composite
produced in the
diagram $(\dagger)$ above, uniquely, through
\[\xymatrix@C+10pt{
&Y(a_1^{})\ar[r]^-{Y(f_1)}\ar[d] &Y(b_1)\ar[d]
\ar[r]^-{Y(g^{}_1)}
&Y(c_1^{})\ar[d]\ar[r]^-{Y(h_1)} &Y(\T a_1)\ar[d] \\
(\dagger\dagger)&A\ar[r]^-{f}\ar[d]^{\alpha'_i} &
B \ar[d]^{\beta'_i}\ar[r]^-{g} &
C\ar[d]^{\gamma_i}\ar[r]^-{h} &
\T A\ar[d]
\\
&\wt a_i\ar[r]^-{\wt f_i} &\wt b_i\ar[r]^-{\wt g^{}_i} &
\wt c_i\ar[r]^-{\wt h_i} &\wt \T a_i &
}\]
This defines for us a morphism $\gamma_i:C\la\wt c_i$.
The first step is to prove that $\alpha'_i=\alpha_i$ and
$\beta'_i=\beta_i$, with $\alpha_i$ and $\beta_i$
as in the hypotheses of the Lemma.
Once we prove that, we will have
shown that the morphism $\gamma_i:C\la\wt c_i$
above renders commutative the diagram in the
statement of the Lemma.

Recall, Lemma~\ref{L77.7}
produced for us morphisms
$\alpha_i:A\la\wt a_i$ and $\beta_i:B\la\wt b_i$,
as well as Cauchy sequences $a_*$, $\wt a_{i,*}$
$b_*$ in $\wt b_{i,*}$, with
isomorphisms in $\fl(\cs)$:
\[\xymatrix@C+0pt@R-10pt{
\colim\,Y(a_*) \ar@{=}[d] &
\colim\,Y(\wt a_{i,*}) \ar@{=}[d] &
\colim\,Y(b_*) \ar@{=}[d] &
\colim\,Y(\wt b_{i,*}) \ar@{=}[d] \\
A &\wt a_i &B &\wt b_i 
}\]
In Construction~\ref{C77.11},
we did some passing to subsequences,
and by the time we reached Construction~\ref{C77.11}(iv)
we had settled 
on maps of Cauchy sequences
$\alpha_{i,*}:a_*\la\wt a_{i,*}$ and
$\beta_{i,*}:b_*\la\wt b_{i,*}$
with $\alpha_i=\colim\,Y(\alpha_{i,*})$
and $\beta_i=\colim\,Y(\beta_{i,*})$.
These have the good properties
summarized in Construction~\ref{C77.11}(v).
In any case: the commutativity of the
two squares
\[\xymatrix@C+10pt{
Y(a_1^{})\ar[r]^-{Y(\alpha_{i,1}^{})}\ar[d] &Y(\wt a_{i,1}^{})\ar[d]
&  & Y(b_1^{})\ar[r]^-{Y(\beta_{i,1}^{})}\ar[d]
&Y(\wt b_{i,1}^{})\ar[d]\\
A\ar[r]^-{\alpha_{i}} &\wt a_{i} & &
B\ar[r]^-{\beta_i} &\wt b_i
}\]
comes because the vertical maps
correspond to the map from the first term
of a sequence
to the colimit. The {unique}
factorization of
\cite[Lemma~\ref{L28.102}]{Neeman18A},
 of the
composites 
\[\xymatrix@C+10pt{
Y(a_1^{})\ar[r]^-{Y(\alpha_{i,1}^{})} &Y(\wt a_{i,1}^{})\ar[r] 
& \wt a_i & Y(b_1^{})\ar[r]^-{Y(\beta_{i,1}^{})}
&Y(\wt b_{i,1}^{})\ar[r] &\wt b_i
}\]
as the composites
\[\xymatrix@C+10pt{
Y(a_1^{})\ar[r] &A\ar[r]^-{\exists!} 
& \wt a_i & Y(b_1^{})\ar[r]
&B\ar[r]^-{\exists!} &\wt b_i\ ,
}\]
which follows from the good properties
of the Cauchy sequences $a_*$ and $b_*$
spelled out in
Construction~\ref{C77.11}(ii) and 
\ref{C77.11}(iv),
has to be through the maps 
$\alpha_i:A\la\wt a_i$ and
$\beta_i:B\la\wt b_i$ that we fixed back in
way back in Lemma~\ref{L77.7}.
The bottom half of the diagram $(\dagger\dagger)$
is therefore
\[\xymatrix@C+10pt{
A\ar[r]^-{f}\ar[d]^{\alpha_i} &
B \ar[d]^{\beta_i}\ar[r]^-{g} &
C\ar[d]^{\gamma_i}\ar[r]^-{h} &
\T A\ar[d]^{\T\alpha_i}
\\
\wt a_i\ar[r]^-{\wt f_i} &\wt b_i\ar[r] &
\wt c_i\ar[r] &\wt \T a_i &
}\]
and this proves the assertions in the first 
paragraph of the current Lemma.

Now for the proof of (ii):
we were given, back in Lemma~\ref{L77.7}(iii),
that
$\alpha_i:A\la\wt a_i$ is of type-$(1+n_{i+1}^{})$
with
respect to $(a_*,\wt a_{i,*})$.
Since $n_{i+1}^{}>n_i$, it
trivially follows
that $\alpha_i:A\la\wt a_i$ is
also of type-$(1+n_i)$
with
respect to $(a_*,\wt a_{i,*})$.
And in Lemma~\ref{L77.7}(iv) we chose
$\beta_i:B\la\wt b_i$ to be of type-$(1+n_i)$ with
respect to $(b_*,\wt b_{i,*})$.
But $\gamma_i:C\la \wt c_i$ was only recently
constructed, and we need to estimate its type
with respect to $(c_*,\wt c_{i,*})$.

We propose to do this by applying Lemma~\ref{L3.307}.
We set $k=\ell=1$, and the commutative square
\[\xymatrix{
Y(c_1^{})\ar[r]^-{Y(\gamma_{i,1}^{})}\ar[d] &
   Y(\wt c_{i,1}^{})\ar[d] \\
C\ar[r]^-{\gamma_i} \ar[r] & \wt c_i
}\]
comes from the definition of
the map $\gamma_i:C\la\wt c_i$, in the passage
from the composite diagram $(\dagger)$ to its
unique factorization as the composite $(\dagger\dagger)$.
In (iv) above we showed that, in the
triangle
$c_1^{}\stackrel{\gamma_1^{}}\la\wt c_{i,1}^{}\la\ov\delta_1$
we have $\ov\delta_1\in\cm_{1+n_i}$.
From Construction~\ref{C77.11}(ii), made
more precise after passing to
subsequences as in Construction~\ref{C77.11}(v),
we know
that, for all $h>j\geq 1$, in the
triangle $c_j\la c_h^{}\la\ov\delta_{j,h}$,
we have $\ov\delta_{j,h}\in\cm_{1+n_{i+2}^{}}$.
Therefore certainly $\T\ov\delta_{j,h}\in\cm_{1+n_i}$.
And, again for pairs
of integers $h>j\geq 1$, we have the triangles
$\wt c_{i,j}\la\wt c_{i,h}^{}\la\wt\delta''_{i,j,h}$;
and in Construction~\ref{C77.11}(viii) we
proved that $\wt\delta''_{i,j,h}\in\cm_{n_{i+2}^{}}$.
As $n_{i+2}^{}>1+n_i$ we deduce
that $\wt\delta''_{i,j,h}\in\cm_{1+n_i}$.

All the conditions of Lemma~\ref{L3.307} are
satisfied, proving that
$\gamma_i:C\la \wt c_i$ is of
type-$(1+n_i)$
with respect to $(c_*,\wt c_{i,*})$. This ends the
proof of (ii).

And finally we come to the proof of (iii);
we are given some integer $\ell>0$,
and add the assumption that $c_j\in\cm_\ell$ for
all $j\gg0$.

Let $k$ be large enough so that
$1+n_k>\ell$.
Since the assertion of (iii)
concerns only the behavior
of  $\wt c_i$ and $\wt c_{i,*}$ 
for large $i$, we will
study only the $\wt c_i=\colim\,Y(\wt c_{i,*})$
meeting the condition $i>k$.
Appealing to Construction~\ref{C77.11}(viii)
again, we know that for pairs
of integers $h>j\geq 1$, the triangles
$\wt c_{i,j}\la\wt c_{i,h}^{}\la\wt\delta''_{i,j,h}$
satisfy $\wt\delta''_{i,j,h}\in\cm_{n_{i+2}^{}}\subset\cm_\ell$.
Therefore, as long as $i\geq k$,
it suffices to find a single $\wt c_{i,j}$ belonging
to $\cm_\ell$; the triangle would prove that,
for any $h>j$, we have $\wt c_{i,h}\in\cm_\ell*\cm_\ell\subset\cm_\ell$.

Next we apply
\cite[Lemma~\ref{L28.19}(i)]{Neeman18A} to
the morphism
$\gamma_i:C\la \wt c_i$ and the Cauchy
sequences $c_*$ and $\wt c_{i,*}$ with
$C=\colim\,Y(c_*)$ and with
$\wt c_i=\colim\,Y(\wt c_{i,*})$.
Replacing $\wt c_{i,*}$ by a subsequence,
we may choose a map of
Cauchy sequences $\gamma_{i,*}:c_*\la\wt c_{i,*}$,
with $\gamma_i=\colim\,Y(\gamma_{i,*})$.
And replacing $\wt c_{i,*}$ by a subsequence is
harmless---after all we have proved that for
each $i>k$, if a single $\wt c_{i,j}$ satisfies
$\wt c_{i,j}\in\cm_\ell$, then all but a finite
number do.

Recall that we have already proved (ii): we
know that the morphism
$\gamma_i:C\la \wt c_i$ is of
type-$(1+n_i)$
with respect to $(c_*,\wt c_{i,*})$.
Next we appeal to
\cite[Lemma~\ref{L28.19}(ii)]{Neeman18A}, and
complete the map of
Cauchy sequences $\gamma_{i,*}:c_*\la\wt c_{i,*}$,
created in the previous paragraph,
to a Cauchy sequence of triangles
$c_*\stackrel{\gamma_{i,*}}\la\wt c_{i,*}\la\ov\delta_*\la\T c_*$. We know that, for integers $m\gg0$, we have
$\ov\delta_m\in\cm_{1+n_i}\subset\cm_\ell$,
where the inclusion is once again because $i>k$.
But by hypothesis $c_m\in\cm_\ell$ for $m\gg0$,
and therefore $m\gg0$ implies
that $\wt c_{i,m}\in\cm_\ell*\cm_\ell\subset\cm_\ell$.
And $\wt c_i=\colim\,Y(\wt c_{i,*})$ must belong
to $\cn_\ell=\fs(\cs)\cap\cl_\ell$.

This completes the proof of (iii), and of the
entire Lemma.
\eprf

Before we move forward, it might be best to collect
on the highlights in one place---the technical
details that entered the proofs can now
safely be forgotten.

\smr{S77.17}
As in most of this section,
we fix in $\cs$ 
a Cauchy sequence of
distinguished triangles
$a_*\stackrel{f_*}\la b_*\stackrel{g_*}\la c_*\stackrel{h_*}\la\T a_*$, and let 
$A\stackrel f\la B\stackrel g\la C\stackrel h\la\T A$
be the strong triangle in $\fl(\cs)$ obtained
as the colimit.

In Lemma~\ref{L77.7}(ii) we 
produced (among other things), in the category
$\cs$ and for each
integer $i>0$, Cauchy sequences $\wt a_{i,*}$ and
$\wt b_{i,*}$. Set
$\wt a_i=\colim\,Y(\wt a_{i,*})$ and 
$\wt b_i=\colim\,Y(\wt b_{i,*})$.
In Construction~\ref{C77.11}(viii) we
produced, again in the category $\cs$, a
Cauchy sequence $\wt c_{i,*}$,
while in Construction~\ref{C77.11}(ix)
we put 
$\wt c_i=\colim\,Y(\wt c_{i,*})$.
And what will be relevant to us for the rest of
the section is that, in Lemma~\ref{L77.15},
we assembled it all into a commutative diagram
\[\xymatrix@C+10pt{
A\ar[r]^-{f}\ar[d]^{\alpha_i} &
B \ar[d]^{\beta_i}\ar[r]^-{g} &
C\ar[d]^{\gamma_i}\ar[r]^-{h} &
\T A\ar[d]^{\T\alpha_i}
\\
\wt a_i\ar[r]^-{\wt f_i} &\wt b_i\ar[r]^-{\wt g_i^{}} &
\wt c_i\ar[r]^-{\wt h_i} &\wt \T a_i &
}\]
And the relevant information about
this entire process, for the rest of
the current section, can be summarized in
\be
\item
The objects $\wt a_i$, $\wt b_i$, $\wt c_i$ and
$\T\wt a_i$ all lie in $\fs(\cs)\cap\cl_{1+n_{i+2}}^\perp$.  
\setcounter{enumiv}{\value{enumi}}
\ee
The reader can check that this was proved in
Construction~\ref{C77.11}(ix).
\be
\setcounter{enumi}{\value{enumiv}}
\item
The three morphisms  
\[
\alpha_i:A\la\wt a_i\,,\qquad
\beta_i:B\la\wt b_i\,,\qquad
\gamma_i:C\la\wt c_i
\]
are of type-$(1+n_i)$, with respect to [respectively]
the pairs of Cauchy sequences
\[
(a_*,\wt a_{i,*}),\qquad
(b_*,\wt b_{i,*}),\qquad
(c_*,\wt c_{i,*}).
\]
\setcounter{enumiv}{\value{enumi}}
\ee
For the proof the reader is referred to
Lemma~\ref{L77.15}(ii).
\be
\setcounter{enumi}{\value{enumiv}}
\item
If we assume that the object $B$ belongs to
the category $\fs(\cs)$, then for $i\gg0$
we may choose map $\beta_{i,*}:B\la \wt b_i$ to be 
the identity map $\id:B\la B$.
\setcounter{enumiv}{\value{enumi}}
\ee
For this assertion we refer the reader to
Lemma~\ref{L77.7}(v).
\be
\setcounter{enumi}{\value{enumiv}}
\item
If we are given an integer $\ell>0$,
and assume that the Cauchy sequence $c_*$
is such that $c_n\in\cm_\ell$ for $n\gg0$,
then $\wt c_i$ will belong $\cn_\ell$ for
$i\gg0$.
\setcounter{enumiv}{\value{enumi}}
\ee
For this final statement the reader is referred
to Lemma~\ref{L77.15}(iii).
\esmr

And now we are ready for the following.

\lem{L77.19}
Let the notation be as in Summary~\ref{S77.17}.
Then for any integer $i>0$, the commutative
diagram 
\[\xymatrix@C+10pt{
A\ar[r]^-{f}\ar[d]^{\alpha_i} &
B \ar[d]^{\beta_i}\ar[r]^-{g} &
C\ar[d]^{\gamma_i}\ar[r]^-{h} &
\T A\ar[d]^{\T\alpha_i}
\\
\wt a_i\ar[r]^-{\wt f_i} &\wt b_i\ar[r]^-{\wt g_i^{}} &
\wt c_i\ar[r]^-{\wt h_i} &\wt \T a_i 
}\]
factors, uniquely, as
\[\xymatrix@C+20pt{
A\ar[r]^-{f}\ar[d]^{\alpha_{i+3}^{}} &
B \ar[d]^{\beta_{i+3}^{}}\ar[r]^-{g} &
C\ar[d]^{\gamma_{i+3}^{}}\ar[r]^-{h} &
\T A\ar[d]^{\T\alpha_{i+3}^{}}
\\
\wt a_{i+3}^{}\ar[r]^-{\wt f_{i+3}^{}} \ar[d]^{\wt\alpha_i} &\wt b_{i+3}^{}\ar[r]^-{\wt g_{i+3}^{}} \ar[d]^{\wt\beta_i} &
\wt c_{i+3}^{}\ar[r]^-{\wt h_{i+3}^{}} \ar[d]^{\wt\gamma_i} &\wt \T a_{i+3}^{} \ar[d]^{\T\wt\alpha_i} \\
\wt a_i\ar[r]^-{\wt f_i} &\wt b_i\ar[r]^-{\wt g_i^{}} &
\wt c_i\ar[r]^-{\wt h_i} &\wt \T a_i 
}\]
where the morphisms
\[
\wt\alpha_i:\wt a_{i+3}^{}\la\wt a_i\,,\qquad
\wt\beta_i:\wt b_{i+3}^{}\la\wt b_i\,,\qquad
\wt\gamma_i:\wt c_{i+3}^{}\la\wt c_i
\]
are of type-$(1+n_i)$, with respect to [respectively]
the pairs of Cauchy sequences
\[
(\wt a_{i+3, *}^{},\wt a_{i,*}),\qquad
(\wt b_{i+3, *}^{},\wt b_{i,*}),\qquad
(\wt c_{i+3, *}^{},\wt c_{i,*}).
\]
\elem

\prf
From Summary~\ref{S77.17}(ii) we know that the
three morphisms
\[
\alpha_{i+3}^{}:A\la\wt a_{i+3}^{}\,,\qquad
\beta_{i+3}^{}:B\la\wt b_{i+3}^{}\,,\qquad
\gamma_{i+3}^{}:C\la\wt c_{i+3}^{}
\]
are of type-$(1+n_{i+3}^{})$, with respect to [respectively]
the pairs of Cauchy sequences
\[
(a_*,\wt a_{i+3,*}),\qquad
(b_*,\wt b_{i+3,*}),\qquad
(c_*,\wt c_{i+3,*})
\]
while Summary~\ref{S77.17}(i) informs us that
the objects
$\wt a_i$, $\wt b_i$, $\wt c_i$ and
all lie in $\fs(\cs)\cap\cl_{1+n_{i+2}}^\perp$.  
Corollary~\ref{C3.378} tells us that,
for any object $B\in\fs(\cs)\cap\cl_{1+n_{i+2}}^\perp$,
the functor $\Hom(-,B)$ takes any type-$(1+n_{i+3}^{})$
to an isomorphism. Applying this to the present situation
we have that each of the morphisms
\[
\alpha_i:A\la\wt a_i\,,\qquad
\beta_i:B\la\wt b_i\,,\qquad
\gamma_i:C\la\wt c_i
\]
factors, uniquely, as
\[
A\stackrel{\alpha_{i+3}^{}}\la \wt a_{i+3}^{}\stackrel{\wt\alpha_i}\la \wt a_i\,,\qquad
B\stackrel{\beta_{i+3}^{}}\la \wt b_{i+3}^{}\stackrel{\wt\beta_i}\la \wt b_i\,,\qquad
C\stackrel{\gamma_{i+3}^{}}\la \wt c_{i+3}^{}\stackrel{\wt\gamma_i}\la \wt c_i\,,
\]
where, by Lemma~\ref{L3.341}(ii), the morphisms
\[
\wt\alpha_i:\wt a_{i+3}^{}\la\wt a_i\,,\qquad
\wt\beta_i:\wt b_{i+3}^{}\la\wt b_i\,,\qquad
\wt\gamma_i:\wt c_{i+3}^{}\la\wt c_i
\]
must be of type-$(1+n_i)$, with respect to [respectively]
the pairs of Cauchy sequences
\[
(\wt a_{i+3, *}^{},\wt a_{i,*})\,,\qquad
(\wt b_{i+3, *}^{},\wt b_{i,*})\,,\qquad
(\wt c_{i+3, *}^{},\wt c_{i,*})\,.
\]
Anyway: with the three new morphisms at hand we can form
the diagram below
\[\xymatrix@C+20pt{
A\ar[r]^-{f}\ar[d]^{\alpha_{i+3}^{}} &
B \ar[d]^{\beta_{i+3}^{}}\ar[r]^-{g} &
C\ar[d]^{\gamma_{i+3}^{}}\ar[r]^-{h} &
\T A\ar[d]^{\T\alpha_{i+3}^{}}
\\
\wt a_{i+3}^{}\ar[r]^-{\wt f_{i+3}^{}} \ar[d]^{\wt\alpha_i} &\wt b_{i+3}^{}\ar[r]^-{\wt g_{i+3}^{}} \ar[d]^{\wt\beta_i} &
\wt c_{i+3}^{}\ar[r]^-{\wt h_{i+3}^{}} \ar[d]^{\wt\gamma_i} &\wt \T a_{i+3}^{} \ar[d]^{\T\wt\alpha_i} \\
\wt a_i\ar[r]^-{\wt f_i} &\wt b_i\ar[r]^-{\wt g_i^{}} &
\wt c_i\ar[r]^-{\wt h_i} &\wt \T a_i 
}\]
When we delete the middle row, the construction tells
us that the vertical composites deliver
the commutative diagram of Lemma~\ref{L77.15}, for
the integer $i>0$. If we delete the bottom row,
we are left with
the commutative diagram of 
Lemma~\ref{L77.15} for
the integer $i+3>0$.
Thus the composites in each of the squares
\[\xymatrix@C+20pt{
 \wt a_{i+3}^{}\ar[r]^-{\wt f_{i+3}^{}}\ar[d]^{\wt\alpha_i}& \wt b_{i+3}^{}
\ar[d]^{\wt\beta_i} &
\wt b_{i+3}^{} \ar[r]^-{\wt g_{i+3}^{}}
\ar[d]^{\wt\beta_i} &
\wt c_{i+3}^{}\ar[d]^{\wt\gamma_i} & 
\wt c_{i+3}^{}\ar[r]^-{\wt h_{i+3}^{}}\ar[d]^{\wt\gamma_i} & \T \wt a_{i+3}^{}\ar[d]^{\T \wt\alpha_i}
\\
\wt a_i\ar[r]^-{\wt f_i} & \wt b_i&
\wt b_i\ar[r]^-{\wt g_i}   \ar[r] &
\wt c_i &
\wt c_i\ar[r]^-{\wt h_i}  & \T \wt a_i
}\]
give a pair of maps rendering equal the composites
\[\xymatrix@C-3pt{
A\ar[r]^-{\alpha_{i+3}^{}} & \wt a_{i+3}^{}\ar@<0.5ex>[r] \ar@<-0.5ex>[r] & \wt b_i & 
B\ar[r]^-{\beta_{i+3}^{}} & \wt b_{i+3}^{}\ar@<0.5ex>[r] \ar@<-0.5ex>[r] & \wt c_i & 
C\ar[r]^-{\gamma_{i+3}^{}} & \wt c_{i+3}^{}\ar@<0.5ex>[r] \ar@<-0.5ex>[r] & \T \wt a_i\ .
}\]
But now we remember that the three morphisms
\[
\alpha_{i+3}^{}:A\la\wt a_{i+3}^{}\,,\qquad
\beta_{i+3}^{}:B\la\wt b_{i+3}^{}\,,\qquad
\gamma_{i+3}^{}:C\la\wt c_{i+3}^{}
\]
are of type-$(1+n_{i+3}^{})$, with respect to [respectively]
the pairs of Cauchy sequences
\[
(a_*,\wt a_{i+3,*}),\qquad
(b_*,\wt b_{i+3,*}),\qquad
(c_*,\wt c_{i+3,*}),
\]
and Summary~\ref{S77.17}(i) informs us that
the objects
$\wt b_i$, $\wt c_i$ and
$\T\wt a_i$ all lie in $\fs(\cs)\cap\cl_{1+n_{i+2}}^\perp$.  
Applying
Corollary~\ref{C3.378} again gives that, in each case,
the factorization of the composite through
a morphism of type-$(1+n_{i+3}^{})$ is unique,
and hence the three
squares must commute.
\eprf

We have done the hard work, and now the time has come
to put it all together.

\medskip

\nin
\emph{Proof of Proposition~\ref{P77.1} and of its
Refinement~\ref{R77.5}.}\ \
Lemma~\ref{L77.15} gives, for every integer $i>0$,
a commutative diagram
\[\xymatrix@C+10pt{
A\ar[r]^-{f}\ar[d]^{\alpha_{3i}^{}} &
B \ar[d]^{\beta_{3i}^{}}\ar[r]^-{g} &
C\ar[d]^{\gamma_{3i}^{}}\ar[r]^-{h} &
\T A\ar[d]^{\T\alpha_{3i}^{}}
\\
\wt a_{3i}^{}\ar[r]^-{\wt f_{3i}^{}} &\wt b_{3i}^{}\ar[r]^-{\wt g_{3i}^{}} &
\wt c_{3i}^{}\ar[r]^-{\wt h_{3i}^{}} &\wt \T a_{3i}^{} &
}\]
with the bottom row a distinguished triangle in
$\fs(\cs)$.
By Lemma~\ref{L77.19} we have, for every integer $i$,
a (unique) factorization of this commutative diagram
through the larger commutative diagram
\[\xymatrix@C+20pt{
A\ar[r]^-{f}\ar[d]^{\alpha_{3(i+1)}^{}} &
B \ar[d]^{\beta_{3(i+1)}^{}}\ar[r]^-{g} &
C\ar[d]^{\gamma_{3(i+1)}^{}}\ar[r]^-{h} &
\T A\ar[d]^{\T\alpha_{3(i+1)}^{}}
\\
\wt a_{3(i+1)}^{}\ar[r]^-{\wt f_{3(i+1)}^{}} \ar[d]^{\wt\alpha_{3i}^{}} &\wt b_{3(i+1)}^{}\ar[r]^-{\wt g_{3(i+1)}^{}} \ar[d]^{\wt\beta_{3i}^{}} &
\wt c_{3(i+1)}^{}\ar[r]^-{\wt h_{3(i+1)}^{}} \ar[d]^{\wt\gamma_{3i}^{}} &\wt \T a_{3(i+1)}^{} \ar[d]^{\T\wt\alpha_{3i}^{}} \\
\wt a_{3i}^{}\ar[r]^-{\wt f_{3i}^{}} &\wt b_{3i}^{}\ar[r]^-{\wt g_{3i}^{}} &
\wt c_{3i}^{}\ar[r]^-{\wt h_{3i}^{}} &\wt \T a_{3i}^{} 
}\]
In the category $\MMod\cs$ we deduce a map from
the strong triangle
$A\stackrel f\la B\stackrel g\la C\stackrel h\la\T A$
to the inverse limit of the sequence
$\wt a_{3*}^{}\stackrel{\wt f_{3*}^{}}\la \wt b_{3*}^{}\stackrel{\wt g_{3*}^{}}\la\wt c_{3*}^{}\stackrel{\wt h_{3*}^{}}\la\T\wt a_{3*}^{}$,
which is a sequence of distinguished triangles
in the category $\fs(\cs)$.
And in Summary~\ref{S77.17}(iii) we noted that, if the object
$B$ belongs to $\fs(\cs)$, then the morphisms
$\beta_{3i}^{}:B\la \wt b_{3i}^{}$ can be chosen
to be $\id:B\la B$ for all $i\gg0$, and hence the
inverse sequence $\wt b_{3i}^{}$ eventually stabilizes to
be $\cdots\la B\stackrel\id\la B\stackrel\id\la B$.
And in Summary~\ref{S77.17}(iv) we observed that,
if the sequence $c_*$ with $C=\colim\,Y(c_*)$
is such that $c_n\in\cm_\ell$  for $n\gg0$, then
we have that $\wt c_{3*}^{}$ satisfies
$\wt c_{3i}^{}\in\cn_\ell$ for $i\gg0$.
To complete the proof of Proposition~\ref{P77.1}, and of its
Refinement~\ref{R77.5}, we
need to show two things
\be
\item
The inverse sequence
$\wt a_{3*}^{}\stackrel{\wt f_{3*}^{}}\la \wt b_{3*}^{}\stackrel{\wt g_{3*}^{}}\la\wt c_{3*}^{}\stackrel{\wt h_{3*}^{}}\la\T\wt a_{3*}^{}$,
which we can view as a direct sequence
in $\fs(\cs)\op$, is Cauchy with respect to the
metric $\{\cn_i,\,i\in\nn\}$.
\item
The vertical maps in the
diagram below are isomorphisms
\[\xymatrix@C+40pt@R-2pt{
A\ar[d]^{\clim\alpha_{3*}^{}}\ar[r]^-f & B\ar[d]^{\clim\beta_{3*}^{}}\ar[r]^-g & C\ar[d]^{\clim\gamma_{3*}^{}}\ar[r]^-h &\T A\ar[d]^{\clim\T\alpha_{3*}^{}}\\
\clim \wt a_{3*}^{}\ar[r]_-{\clim \wt f_{3*}^{}}& \clim \wt b_{3*}^{}\ar[r]_-{\clim \wt g_{3*}^{}}& \clim \wt c_{3*}^{}\ar[r]_-{\clim \wt h_{3*}^{}}&\clim \T\wt a_{3*}^{}
}\]
\ee
First we prove (i). And we will confine
ourselves to proving that the
inverse sequence $\wt a_{3*}^{}$ is Cauchy,
leaving to the reader the proof
that $\wt b_{3*}^{}$ and $\wt c_{3*}^{}$ also are.

By Lemma~\ref{L77.19} we know
that, for every $i>0$, the morphism
$\wt\gamma_{3i}^{}:\wt a_{3(i+1)}^{}\la \wt a_{3i}^{}$
is of type-$(1+n_{3i}^{})$ with respect to
the Cauchy sequences $(\wt a_{3(i+1),*}^{}, \wt a_{3i,*}^{})$.
From Lemma~\ref{L3.309} we learn that
we may complete the morphism $\wt\gamma_{3i}^{}:\wt a_{3(i+1)}^{}\la \wt a_{3i}^{}$, in the
category $\fl(\cs)$, to a strong triangle
$E\la \wt a_{3(i+1)}^{}\stackrel{\wt\gamma_{3\scriptstyle{i}}^{}}\la \wt a_{3i}^{}\la \T E$, with $E$ and $\T E$ both in $\cl_{n_{3\scriptstyle{i}}^{}}$.
But in this strong triangle the objects
$\wt a_{3(i+1)}^{}$ and $\wt a_{3i}^{}$
both belong to $\fs(\cs)$, and by
\cite[Remark~\ref{R200976}]{Neeman18A}
this strong triangle must be a distinguished triangle
in $\fs(\cs)$. In particular $E$ belongs to
$\cn_{n_{3\scriptstyle{i}}^{}}=\cl_{n_{3\scriptstyle{i}}^{}}\cap\fs(\cs)$, and the
sequence $\wt a_{3*}^{}$ is indeed Cauchy.
This completes the proof of (i).

And finally we prove (ii). Once again we confine
ourselves to proving that the map
$\clim\alpha_{3*}^{}:A\la\clim\wt a_{3*}^{}$ is
an isomorphism, since the
proof that the other maps are isomorphisms
is essentially identical.

By (i) we know that the inverse sequence
$\wt a_{3*}^{}$ is Cauchy.
And Summary~\ref{S77.17}(ii) gives that, for each
integer $i>0$, the map $\alpha_{3i}^{}:A\la \wt a_{3i}^{}$
is of type-$(1+n_{3i}^{})$ with respect to
$(a_*,\wt a_{3i,*}^{})$. From
Lemma~\ref{L3.309} it follows that, for
every integer $i>0$,
there exists in $\fl(\cs)$
a strong triangle
$E_i\la A\stackrel{\alpha_{3i}^{}}\la \wt a_{3i}^{}\la\T E_i$,
with $E_i$ and $\T E_i$ both in
$\cl_{n_{3\scriptstyle{i}}^{}}$.
Therefore Lemma~\ref{L3.90909} applies,
proving that
the map $\clim\alpha_{3*}^{}:A\la\clim\wt a_{3*}^{}$ is
an isomorphism
[as is $\colim\,\wh Y(\alpha_{3*}^{}):\colim\,\wh Y(\wt a_{3*}^{})\la\wh Y(A)$]. This completes the proof of (ii),
as well as the proof of
Proposition~\ref{P77.1} and its Refinement~\ref{R77.5}.
\hfill{$\Box$}

\section{The excellence of $\fs(\cs)\op$, and the inclusion $\cs\op\subset\fs\big[\fs(\cs)\op\big]$}
\label{S97}

In Section~\ref{S77} we studied the relationship
between strong triangles in $\Mod\cs$ and
strong triangles in $\Mod{\fs(\cs)\op}$,
and now we are ready to develop the theory
further.

\pro{P97.1}
Let $\cs$ be a triangulated category, and assume that 
$\{\cm_i\mid i\in\nn\}$
is an excellent metric on $\cs$.
Then the metric $\{\cn\op_i,\,i\in\nn\}$
on the category $\fs(\cs)\op$ is also excellent.
\epro

\prf
By Lemma~\ref{L3.1}
we know that the metric $\{\cn\op_i,\,i\in\nn\}$
is good, and what needs showing is that
it satisfies the extra hypotheses of
Definition~\ref{D3.3}.
  
We begin with
Definition~\ref{D3.3}(i). For this,
start out with equalities that move the
information we have from $\big[\MMod\cs\big]\op$ to $\MMod{\fs(\cs)\op}$
\[
\begin{array}{ccccl}
\wh Y\big[\fs(\cs)\op\big]
 &=&\wh Y\left(\left[\fl(\cs)\cap\bigcup_{i\in\nn}\cl_i^\perp\right]\op\right)
     & &
  \text{by Remark~\ref{R3.2??}(ii),} \\
  &=&\fl\big[\fs(\cs)\op\big]\cap\bigcup_{i\in\nn}{^\perp\wh\cl_i}
   & &\text{by Proposition~\ref{P3.94949}.}\\
\end{array}
\]
And then continue with the equalities below, where
the passage from the first to the second line is
by intersecting both sides with
$\wh Y\big[\fs(\cs)\op\big]$.
\[
\begin{array}{ccccl}
\wh Y\big[\fs(\cs)\op\big]
  &=&\fl\big[\fs(\cs)\op\big]\cap\bigcup_{i\in\nn}{^\perp\wh\cl_i}
   & &\text{by previous equalities,}\\
  &=&\wh Y\big[\fs(\cs)\op\big]\cap\bigcup_{i\in\nn}{^\perp\wh\cl_i}
   & &\text{intersecting with $\wh Y\big[\fs(\cs)\op\big]$,}\\
  &=&\wh Y\big[\fs(\cs)\op\big]\cap\bigcup_{i\in\nn}{^\perp\wh Y(\cn\op_i)}
   & &\text{by Lemma~\ref{L3.2.5}(i).}\\
\end{array}
\]
As $\wh Y$ is fully faithful on $\fs(\cs)\op$ and
its subcategories, we deduce that
\[
\fs(\cs)\op\eq\bigcup_{i\in\nn}{^\perp\cn\op_i}\ ,
\]
completing the proof that Definition~\ref{D3.3}(i)
holds for $\fs(\cs)\op$ with
its metric $\{\cn\op_i,\,i\in\nn\}$.

Now for the proof of Definition~\ref{D3.3}(ii).
Let $m>0$ be any integer, and now take $n>m+1$ to be
the integer of Definition~\ref{D3.3}(iii) when
applied to the category $\cs$, with its
metric $\{\cm_i,\,i\in\nn\}$, and with the integer
$m+1>0$.

If 
$F$ is an object of
$\fs(\cs)$, then choose in $\cs$ a
Cauchy sequence $f_*$ with $F=\colim\,Y(f_*)$.
And now apply Lemma~\ref{L3.5}. It produces for
us, in the category $\fl(\cs)$,
a type-$(m+1)$ morphism $F\la D$ with respect
to Cauchy sequences $(f_*,d_*)$, with
$D\in\fs(\cs)\cap\cl_n^\perp$.
Lemma~\ref{L3.309} allows us to
extend $F\la D$ to a strong
triangle $E\la F\la D\la\T E$, with
$E$ and $\T E$ both in $\cl_m$.
As $F$ and $D$ both belong to $\fs(\cs)$,
\cite[Remark~\ref{R200976}]{Neeman18A}
tells us
that this strong triangle
is a distinguished triangle in $\fs(\cs)$.
In particular, $E$ belongs to $\fs(\cs)\cap\cl_m=\cn_m$.

But now, in the category $\fs(\cs)\op$, this distinguished
triangle rewrites
as $E\longleftarrow F\longleftarrow D$, with
$E\in\cn\op_m$ and with
\[
\begin{array}{ccccl}
\wh Y(D)&\in&\wh Y\Big(\big[\fs(\cs)\cap\cl_n^\perp \big]\op\Big) & & \\
&=& \wh Y\big[\fs(\cs)\op\big]\cap{^\perp\wh\cl_n}
 &\quad &\text{by Proposition~\ref{P3.94949}}\\
&=& \wh Y\big[\fs(\cs)\op\big]\cap{^\perp\wh Y(\cn\op_n)}
 &\quad &\text{by Lemma~\ref{L3.2.5}(i).}\\
\end{array}
\]
And because $\wh Y$ is fully faithful on $\fs(\cs)\op$
and its subcategories,
we learn that $D\in{^\perp\cn\op_n}$.

Thus the triangle $E\longleftarrow F\longleftarrow D$,
in the category $\fs(\cs)\op$, meets the criteria
of Definition~\ref{D3.3}(ii).

It remains to prove that Definition~\ref{D3.3}(iii)
holds for $\fs(\cs)\op$, together with its metric
$\{\cn\op_i,\, i\in\nn\}$. Once again we let $m>0$
be any integer, and choose $n>m+1$ to be the integer
corresponding to $m+1$ for the
category $\cs$; only this time we use the
correspondence in
Definition~\ref{D3.3}(ii).

Let $F$ be an object in $\fs(\cs)$, and choose in
$\cs$ a Cauchy sequence $f_*$ with $F=\colim\,Y(f_*)$.
Lemma~\ref{L3.4} allows us to choose a type-$(m+1)$
morphism $f:Y(E)\la F$, with respect to $(E,f_*)$,
and with $E\in{^\perp\cm_n}$. And Remark~\ref{R3.?!?}
gives that
\be
\item
The condition on $E\in\cs$, 
in the morphism $f:Y(E)\la F$, can be rephrased 
as $Y(E)\in{^\perp\cl_n}$.
\setcounter{enumiv}{\value{enumi}}
\ee
Now \cite[Lemma~\ref{L28.19}(i)]{Neeman18A}
permits us, up to replacing
$E$ and $f_*$ by subsequences, to choose
a map of Cauchy sequences $\ph:E\la f_*$ with
$f=\colim\,Y(\ph_*)$, and
\cite[Lemma~\ref{L28.19}(ii)]{Neeman18A}
permits us to extend
this to a Cauchy sequence
$E\stackrel{\ph_*}\la f_*\la d_*\la\T E$
of triangles in
$\cs$. And the assumption
that $f:Y(E)\la F$ is of type-$(m+1)$ with respect
to $(E,f_*)$ means that, in the Cauchy sequence
$d_*$, all but finitely many of
the terms lie in $\cm_{m+1}$.
Applying
$Y$ and taking colimits, we obtain the
strong triangle $Y(E)\la F\la D\la\T Y(E)$.

And now we apply 
of Proposition~\ref{P77.1}
to the
strong triangle $Y(E)\la F\la D\la\T Y(E)$.
More accurately: as in 
Refinement~\ref{R77.5}(i)
we choose our Cauchy sequence of
triangles to be
$E\stackrel{\ph_*}\la f_*\la d_*\la\T E$
above, and note that
\be
\setcounter{enumi}{\value{enumiv}}
\item
The object $F=\colim\,Y(f_*)$
lies in $\fs(\cs)$.
\item
The Cauchy sequence $d_*$ satisfies $d_\ell\in\cm_{m+1}$
for $\ell\gg0$.
\setcounter{enumiv}{\value{enumi}}
\ee
That is, the technical hyppotheses
of 
Refinement~\ref{R77.5}~(ii) and (iii) both
hold. Therefore, by the conclusion of
Refinement~\ref{R77.5},
we can produce in $\fs(\cs)$ a
Cauchy inverse
sequence $\wt e_*\stackrel{\wt\ph_*}\la\wt f_*\la\wt d_*\la\T\wt e_*$ of triangles,
with $\wt d_\ell\in\cn_{m+1}$ for $\ell\gg0$ and where
$\wh f_*$ becomes the constant inverse sequence
$\cdots\la F\stackrel\id\la F \stackrel\id\la F$ after
a finite number of steps. And
the inverse limit of this Cauchy sequence is
the 
strong triangle $Y(E)\la F\la D\la\T Y(E)$.

Now pass to a subsequence with $\wt f_*=F$ constant, and view
this in the category $\fs(\cs)\op$. We have
produced in $\fs(\cs)\op$ a type-$m$ map of Cauchy sequences
$\wt\ph\op_*:F\la\wt e_*$, where $F$ is the constant
sequence. But
\[
\colim\,\wh Y(\wt e_*)\eq \wh Y\Big(\clim\wt e_*\Big)
\eq\wh Y\big(Y(E)\big)
\]
belongs to $\wh Y\Big(\big[{^\perp\cl_n}\big]\op\Big)=\wh\cl_n^\perp$, where the equality is by
Proposition~\ref{P3.94949}. Taking the colimit
in $\MMod{\fs(\cs)\op}$,
of the map of Cauchy sequences $\wh Y(\wt\ph\op_*):\wh Y(F)\la\wh Y(\wt e_*)$,
we obtain that $\wh Y(\ph):\wh Y(F)\la\wh Y\big(Y(E)\big)$
is a type-$m$ morphism with respect to
the pair $(F,\wt e_*)$, and  where
$\wh Y\big(Y(E)\big)\in\wh\cl_n^\perp$.
This completes the proof
that Definition~\ref{D3.3}(iii) holds, and
hence concludes the entire proof of the
current Proposition.
\eprf

\exm{Emoreexcellentmetrics}
Let $\ct$ be a weakly approximable triangulated
category, and let $\cs=\ct^c$
be the subcategory of compact objects
in $\ct$. In Example~\ref{E3.39393901} we saw
that the metric 
$\{\cm_i,\,i\in\nn\}$, of
\cite[Example~\ref{E20.3}(i)]{Neeman18A},
is excellent on $\cs$.

From Proposition~\ref{P97.1} we learn that
the induced metric $\{\cn\op_i,\,i\in\nn\}$
is excellent on
$\fs(\cs)\op=(\ct^b_c)\op$.
We leave it to the reader to check
that this agrees with the
metric on $(\ct^b_c)\op$
given in 
\cite[Example~\ref{E20.3}(ii)]{Neeman18A}.
\eexm

Let $\cs$ be a triangulated category, and assume that 
$\{\cm_i\mid i\in\nn\}$
is an excellent metric on $\cs$.
Now let $\fs(\cs)\op$ be the triangulated
category of
\cite[Theorem~\ref{T28.128}]{Neeman18A},
with its excellent metric $\{\cn\op_i,\,i\in\nn\}$
of Proposition~\ref{P97.1}. We would like to understand
better the relation between these two
categories.

\rmk{R97.3}
Let's begin by only assuming that we have
a good (not necessarily excellent) metric
$\{\cm_i\mid i\in\nn\}$ on the triangulated category $\cs$.
The category $\fs(\cs)$ is triangulated
by \cite[Theorem~\ref{T28.128}]{Neeman18A}, and
has a good metric $\{\cn_i,\,i\in\nn\}$
by Lemma~\ref{L3.1}.
And the observation that we
want to stress here is that \emph{the
category $\fl(\cs)$, 
its subcategories $\cl_i$, and
its strong triangles
determine $\fs(\cs)$ and all
its structure.} After all:
Remark~\ref{R3.2??} says that, in the category
$\MMod\cs$, we have the equality
$\fs(\cs)=\fl(\cs)\cap\bigcup_{i\in\nn}\cl_i^\perp$,
and inside $\fl(\cs)$ this simplifies
to
\[\fs(\cs)\eq\bigcup_{i\in\nn}\cl_i^\perp\ .\]
And, by their definition,
\[
\cn_\ell\eq\fs(\cs)\cap\cl_\ell\eq\left[\bigcup_{i\in\nn}\cl_i^\perp
  \right]\cap\cl_\ell\eq\bigcup_{i\in\nn}\big[\cl_i^\perp\cap\cl_\ell\big]\ .
\]
And by
\cite[Remark~\ref{R200976}]{Neeman18A} we know that the
distinguished triangles on $\fs(\cs)$ are precisely
the strong triangles in $\fl(\cs)$ where two of the
objects (and therefore also the third) belong to
$\fs(\cs)$.
\ermk

So much for good metrics.

\dis{D97.5}
One feature, that sets excellent metrics apart, is that,
by Proposition~\ref{P3.94949}, the functor
$\wh Y$ induces an equivalence
$\fl(\cs)\op\cong\fl\big[\fs(\cs)\op\big]$
and, for each integer $\ell>0$,
it restricts to an equivalence
$\cl\op_\ell\cong\wh\cl_\ell$.
Now recall that, by Remark~\ref{R3.?!?2},
we have inclusion
\[
Y(\cs)\sub\bigcup_{i\in\nn}{^\perp\cl_i}
\]
and
\[
Y(\cm_n)\sub\bigcup_{i\in\nn}\big[{^\perp\cl_i}\cap\cl_n\big]\ .
\]
Since the left hand sides are contained in $\fl(\cs)$,
we can view these as inclusions of subcategories of $\fl(\cs)$.
By this we mean that all perpendiculars are to be viewed as
taken in $\fl(\cs)$, that is $^\perp\cl_n$ means $^\perp\cl_n\cap\fl(\cs)$.

If we apply the operator $(-)\op$, taking
subcategories of $\fl(\cs)$ to
subcategories of 
$\fl(\cs)\op\cong\fl\big[\fs(\cs)\op\big]$,
we obtain the inclusions
\[
Y(\cs)\op\sub\bigcup_{i\in\nn}{\wh\cl_i^\perp}
\]
and
\[
Y(\cm_n)\op\sub\bigcup_{i\in\nn}\big[{\wh\cl_i^\perp}\cap\wh\cl_n\big]\ .
\]

And now we come back
to Remark~\ref{R97.3}, applied to the triangulated
category $\fs(\cs)\op$ with its good metric
$\{\cn\op_i,\,i\in\nn\}$.
Remark~\ref{R97.3} gives formulas
for $\fs\big[\fs(\cs)\op\big]$,
and its excellent metric $\{\wh\cm_i,\,i\in\nn\}$,
as subcategories
of $\fl\big[\fs(\cs)\op\big]$. Explicitly it delivers the
formulas
\[
\fs\big[\fs(\cs)\op\big]\eq\bigcup_{i\in\nn}\wh\cl_i^\perp
\]
and
\[
\wh\cm_n\eq\bigcup_{i\in\nn}\big[{\wh\cl_i^\perp}\cap\wh\cl_n\big]\ .
\]
Comparing the formulas for $Y(\cs)\op$ and
for $Y(\cm_n)\op$, with the
formulas for  $\fs\big[\fs(\cs)\op\big]$ and for $\wh\cm_n$,
we obtain that, as full
subcategories of $\fl\big[\fs(\cs)\op\big]$,
we have inclusions
\[
Y(\cs)\op\subset\fs\big[\fs(\cs)\op\big]\qquad\text{ and }\qquad
Y(\cm_n)\op\subset\wh\cm_n\ .
\]
\edis

This allows us to define

\dfn{D97.7}
Let $\cs$ be a triangulated category, and assume that 
$\{\cm_i\mid i\in\nn\}$
is an excellent metric on $\cs$. Then 
the fully faithful functor $\Psi:\cs\op\la\fs\big[\fs(\cs)\op\big]$
is defined to be the composite
\[\xymatrix@C+30pt{
\cs\op\ar[r]^-{Y\op} & Y(\cs)\op\ar@{^{(}->}[r] & \fs\big[\fs(\cs)\op\big]
}\]
where $Y:\cs\la Y(\cs)$ is the equivalence of
$\cs$ with its image under the Yoneda embedding, and the inclusion
of $Y(\cs)\op$ into $\fs\big[\fs(\cs)\op\big]$ is as full
subcategories of $\fl\big[\fs(\cs)\op\big]$,
and comes from the formulas of Discussion~\ref{D97.5}.
\edfn

And now we are ready for the next main result.

\thm{T97.9}
Let $\cs$ be a triangulated category, and assume that 
$\{\cm_i\mid i\in\nn\}$
is an excellent metric on $\cs$. 
Then
\be
\item
The fully faithful
functor $\Psi:\cs\op\la\fs\big[\fs(\cs)\op\big]$,
of Definition~\ref{D97.7}, is triangulated.
\item
Every object of $\fs\big[\fs(\cs)\op\big]$ is
a direct summand of an object in the essential
image of $\Psi$.
\item
The good metric $\{\wh\cm_i,\,i\in\nn\}$
on the triangulated category $\fs\big[\fs(\cs)\op\big]$,
which is induced from the metric 
$\{\cm_i\mid i\in\nn\}$ on the triangulated
category $\cs$ by applying Lemma~\ref{L3.1} twice,
has the properties that
\begin{enumerate}
\item
$\Psi:\cs\op\la\fs\big[\fs(\cs)\op\big]$
takes
$\cm\op_i$ into $\wh\cm_i$.
\item
Every object of $\wh\cm_i$
is a direct summand of an object in $\Psi(\cm\op_i)$.
\end{enumerate}
\ee
\ethm

\prf
We begin with the proof of (i). Suppose $A\la B\la C\la\T A$
is a distinguished triangle in $\cs$. Then we can form
the constant Cauchy sequence of triangles
\[\xymatrix{
A\ar[r]\ar[d]^\id & B\ar[r]\ar[d]^\id & C\ar[r]\ar[d]^\id &\T A\ar[d]^\id \\
A\ar[r]\ar[d]^\id & B\ar[r]\ar[d]^\id & C\ar[r]\ar[d]^\id &\T A\ar[d]^\id \\
A\ar[r]\ar[d] & B\ar[r]\ar[d] & C\ar[r]\ar[d] &\T A\ar[d] \\
\vdots &\vdots &\vdots &\vdots 
}\]
Applying the functor $Y$ and taking colimits, we deduce
that $Y(A)\la Y(B)\la Y(C)\la\T Y(A)$
is a strong triangle in $\fl(\cs)$. But the
equivalence $\wh Y:\fl(\cs)\op\la\fl\big[\fs(\cs)\op\big]$
respects strong triangles by Proposition~\ref{P77.1},
and hence $\Tm\Psi(A)\la\Psi(C)\la\Psi(B)\la\Psi(A)$
is a strong triangle in $\fl\big[\fs(\cs)\op\big]$.
But as all the objects belong to
$\Psi(\cs\op)\subset\fs\big[\fs(\cs)\op\big]$,
\cite[Definition~\ref{D28.110}]{Neeman18A}
gives that it is
a distinguished triangle in $\fs\big[\fs(\cs)\op\big]$.

Next we observe that (iii)(a) has already been
proved. One of the inclusions we proved, at the end
of Discussion~\ref{D97.5}, was that
$Y(\cm_i)\op\subset\wh\cm_i$. And when
rephrased as in Definition~\ref{D97.7},
this can be rewritten as $\Psi(\cm\op_i)\subset\wh\cm_i$.

Now we will give a
unified proof of (ii) and (iii)(b);
the essential part of proof is
identical in the two cases. Suppose
therefore that $X$ is an object of
$\fs\big[\fs(\cs)\op\big]$; to prove
the more restrictive
(iii)(b) we assume that $X$ belongs
to $\wh\cm_\ell\subset\fs\big[\fs(\cs)\op\big]$.
Remark~\ref{R3.2??}
gives the formula
\[
\fs\big[\fs(\cs)\op\big]\eq
\fl\big[\fs(\cs)\op\big]\cap\bigcup_{i\in\nn}\wh\cl_i^\perp
\]
and intersecting with $\wh\cl_\ell$ gives that
\[
\wh\cm_\ell\eq
\wh\cl_\ell\cap\bigcup_{i\in\nn}\wh\cl_i^\perp
\]
Anyway:
there must exists an integer $m>0$ with
$X\in\fl\big[\fs(\cs)\op\big]\cap\wh\cl_m^\perp$,
and if we assume that $X\in\wh\cm_\ell$ then
$X\in\wh\cl_\ell\cap\cl_m^\perp$.
But Proposition~\ref{P3.94949} tells
us that the functor
$\wh Y$ induces an equivalence
$\fl(\cs)\op\cong\fl\big[\fs(\cs)\op\big]$
satisfying
$\wh Y(\cl\op_\ell)=\wh\cl_\ell$,
and we deduce that there must be
an object $F\in\fl(\cs)\cap{^\perp\cl_m}$ with
$\wh Y(F)=X$. If $X$ belongs to $\wh\cm_\ell$ then
$F$ belongs to $\cl_\ell\cap{^\perp\cl_m}$.

Choose in $\cs$ a Cauchy sequence $f_*$ with
$F=\colim\,Y(f_*)$, and if $F$ belongs to $\cl_\ell$
choose $f_*$ to be contained in 
$\cm_\ell$. And now we apply
Lemma~\ref{L3.4}: starting with the integer $m'\geq\max(\ell+1,m+1)$,
with $\ell$ and $m$ as above, we take the integer $n>m'$
as in Definition~\ref{D3.3}(ii). And then there
exists an object $E\in\cs$, a morphism
$\ph:Y(E)\la F$ in the category $\fl(\cs)$ of
type-$m'$ with respect to $(E,f_*)$, and, since we
have chosen $m'$ to be $\geq(\ell+1)$,
Corollary~\ref{C3.?!?!} guarantees that $E$
belongs to $\cm_\ell$.

Now by Lemma~\ref{L3.309} 
we may complete $\ph$, in the category $\fl(\cs)$,
to a strong triangle $\Tm D\la Y(E)\la F\la D$, with
both $\Tm D$ and $D$ in $\cl_{m'-1}\subset\cl_m$,
where the inclusion is because $m'\geq m+1$.
But as $F\in{^\perp\cl_m}$
and $D\in\cl_m$, the map $F\la D$ must vanish.

By Proposition~\ref{P77.1}, the functor $\wh Y$ takes
strong triangles to strong triangles. Therefore
the top row in the diagram below is a strong triangle in
$\fl\big[\fs(\cs)\op\big]$
\[\xymatrix@R-10pt{
\wh Y(D)\ar[r]^-0 &\wh Y(F)\ar[r]\ar@{=}[d] &
  \Psi(E)\ar[r] &\T \wh Y(D)\\
& X & &
}\]
and the vertical map is an isomorphism because
$F$ was chosen to be such that $\wh Y(F)\cong X$.
And now we recall Lemma~\ref{L3.303}: by assumption
$X$ is an object of $\fs\big[\fs(\cs)\op\big]$,
and hence the functor $\Hom(-,X)$ takes
strong triangles in $\fl\big[\fs(\cs)\op\big]$
to long exact sequences. Applying this to
the above, we deduce that the map
$\Hom\big(\Psi(E),X\big)\la\Hom(X,X)$ must be
surjective, in particular the identity map $\id:X\la X$
lies in the image. Therefore $X$ is a direct summand of
$\Psi(E)\in\Psi(\cs\op)$.
And in the case where $X$ is assumed to lie in
$\wh\cm_\ell$, our conclusion is that $X$ is a
direct summand of $\Psi(E)\in\Psi(\cm\op_\ell)$.
\eprf

\rmk{R97.11}
Of course: if the triangulated category $\cs$ is
idempotent-complete, and if the excellent
metric $\{\cm_i,\,i\in\nn\}$ satisfies the
property that each $\cm_i$ is idempotent-complete,
then the fully faithful, triangulated functor
$\Psi:\cs\op\la\fs\big[\fs(\cs)\op\big]$
must be an equivalence and $\Psi(\cm\op_i)=\wh\cm_i$.

All of our examples will be like that.

But in any case observe that, for
triangulated categories with
excellent metrics, the passage from $\cs$ to $\fs(\cs)\op$
does not change the categories $\fl(\cs)$ and
their subcategories
$\cl_i$. At least, it changes them only up to
the stupidity of moving to
opposite categories. Hence the real issue is only whether
$\cs$ is closed in $\fl(\cs)$
under direct summands, and $\cm_\ell$ is
closed in $\cl_\ell$ under
direct summands. And the formulas
\[
\fs(\cs)=\bigcup_{i\in\nn}\cl_i^\perp\ ,\qquad\qquad
\cn_\ell=\bigcup_{i\in\nn}\big[\cl_i^\perp\cap\cl_\ell\big]
\]
manifestly exhibit that these categories are closed
under direct summands in [respectively] $\fl(\cs)$ and
$\cl_\ell$.

Thus while there could be some idempotent completion
in passing from $\cs$ to 
$\Big(\fs\big[\fs(\cs)\op\big]\big)\op$, the
idempotent-completing stops there. For excellent
metrics, the process of
passing from the metric triangulated category $\cs$ to
the metric triangulated category $\fs(\cs)\op$ becomes
two-periodic after the first step.
\ermk

\section{Good extensions of very good metrics}
\label{S93}

Let $\cs$ be a triangulated category, and assume
$\{\cm_i,\,i\in\nn\}$ is a good metric on $\cs$.
In \cite[Definition~\ref{D21.1}]{Neeman18A}
we met the notion
of a good extension $F:\cs\la\ct$ with respect to
the metric, and the rest
of \cite[Section~\ref{S21}]{Neeman18A}
was largely devoted to showing
how, in the presence of a good extension, the constructions
of $\fl(\cs)$ and $\fs(\cs)$ can become computable.

So far, the current article
has made no mention of good
extensions. We have made
constructions which used only
a single triangulated category $\cs$ and its
metric $\{\cm_i,\,i\in\nn\}$. 
Since we are about to return to
good extensions, a brief reminder might be in place.

\rmd{R93.5}
Let $\cs$ be a triangulated category, and assume
$\{\cm_i\mid i\in\nn\}$ is good metric
on the category $\cs$. Suppose $F:\cs\la\ct$
is a fully faithful, triangulated functor.
\be
\item
The functor $\cy:\ct\la\MMod\cs$
is the restricted Yoneda map; it
takes an object $t\in\ct$ to
$\Hom\big(F(-),t\big)$.
We met this in
\cite[Notation~\ref{N21.-100}]{Neeman18A}.
\item
The fully faithful, triangulated
functor $F:\cs\la\ct$
is a \emph{good extension} with respect to the metric
if, for any Cauchy sequence
$a_*$ in the category $\cs$,
the natural map
$\colim\,Y(a_*)\la\cy\big(\hoco F(a_*)\big)$
is always an isomorphism.
See
\cite[Definition~\ref{D21.1}]{Neeman18A}.
\item
The subcategory $\fl'(\cs)$, of
\cite[Definition~\ref{D21.7}]{Neeman18A},
is by definition
the full subcategory of all isomorphs
of objects $\hoco F(a_*)$, with $a_*$ in
$\cs$ a Cauchy sequence.
\ee
\ermd

As we have said already,
this article has
so far abandoned good extensions.
We have made new constructions,
one of which, made 
back in Definition~\ref{D3.-3}(i), was to introduce
the full subcategories $\cl_i\subset\fl(\cs)$.
And these played a large role in what followed.
In the presence of a good extension $F:\cs\la\ct$
it makes sense to formulate the obvious analogue:

\dfn{D93.11}
Let $\cs$ be a triangulated category with a good metric
$\{\cm_i,\,i\in\nn\}$, and let $F:\cs\la\ct$ be a
good extension with respect to the metric.

The full subcategory $\cl'_i\subset\fl'(\cs)$
has for objects all the homotopy colimits
$\hoco F(a_*)$, where $a_*$ is a Cauchy sequence
contained in
$\cm_i$.
\edfn

It seems appropriate to start this section with a few,
elementary properties of the subcategories $\cl'_i$.

\lem{L93.-1}
Let the notation be as in Definition~\ref{D93.11}.
Then the essential image, of the subcategory
$\cl'_i\subset\ct$ under the functor
$\cy:\ct\la\MMod\cs$ of
\cite[Notation~\ref{N21.-100}]{Neeman18A},
is the subcategory $\cl_i\subset\fl(\cs)$ of
Definition~\ref{D3.-3}(i).
\elem

\prf
An object $X\in\ct$ belongs to $\cl'_i\subset\fl'(\cs)$
[respectively to $\cl_i\subset\fl(\cs)$] if there
exists a Cauchy sequence $a_*$, contained in $\cm_i$,
with $X=\hoco F(a_*)$ [respectively with $X=\colim\,Y(a_*)$].
The lemma now follows from the observation that, because
$F:\cs\la\ct$ is a good extension, the natural map
$\colim\,Y(a_*)\la\cy\big(\hoco F(a_*)\big)$ is
an isomorphism, see
\cite[Definition~\ref{D21.1}(ii)]{Neeman18A} or
Reminder~\ref{R93.5}(ii).
\eprf

\lem{L93.-3}
Let the notation be as in Definition~\ref{D93.11}.
Then $\cy^{-1}\big(\cl_i^\perp\big)=\big(\cl'_i\big)^\perp$.
\elem

\prf
First of all: since $Y(\cm_i)\subset\cl_i$ and
$F(\cm_i)\subset\cl'_i$ we have the inclusions
\[
\begin{array}{ccccc}
\cl_i^\perp&\sub &Y(\cm_i)^\perp &\sub& \fc(\cs)\\
\big(\cl'_i\big)^\perp&\sub &F(\cm_i)^\perp &\sub& \cy^{-1}\big(\fc(\cs)\big)
\end{array}
\]
where the inclusion $Y(\cm_i)^\perp \subset \fc(\cs)$
is by
\cite[Observation~\ref{O28.-1}]{Neeman18A} and
the inclusion
$F(\cm_i)^\perp \subset \cy^{-1}\big(\fc(\cs)\big)$
is by
\cite[Observation~\ref{O21.-1}]{Neeman18A}.
It immediately follows thay both $\cy^{-1}\big(\cl_i^\perp\big)$
and $\big(\cl'_i\big)^\perp$ are subcategories of
$\cy^{-1}\big(\fc(\cs)\big)$.

But now
\cite[Lemma~\ref{L21.13}]{Neeman18A}
says that, if $E\in\cl'_i$ and
$X\in\cy^{-1}\big(\fc(\cs)\big)$, then the natural
map $\Hom_\ct^{}(E,X)\la\Hom_{\MMod\cs}^{}\big(\cy(E),\cy(X)\big)$
is an isomorphism. Therefore $X\in\cy^{-1}\big(\fc(\cs)\big)$
will lie in $\big(\cl'_i\big)^\perp$ if and only if
$\cy(X)$ lies in $\cy(\cl'_i)^\perp=\cl_i^\perp$,
where the last equality in by Lemma~\ref{L93.-1}.
\eprf

\lem{L93.-5}
Let the notation be as in Definition~\ref{D93.11}.
Then $\big(\cl'_i\big)^\perp=F(\cm_i)^\perp$,
and
\[
\fl'(\cs)\cap\cy^{-1}\big(\fc(\cs)\big)\eq
\bigcup_{i\in\nn}\fl'(\cs)\cap\big(\cl'_i\big)^\perp\ .
\]
\elem

\prf
Lemma~\ref{L93.-3} gives
the equality
$\cy^{-1}\big(\cl_i^\perp\big)=\big(\cl'_i\big)^\perp$,
while in
\cite[Observation~\ref{O21.-1}]{Neeman18A} we saw
that $\cy^{-1}\big(Y(\cm_i)^\perp\big)=F(\cm_i)^\perp$.
The equality
\[\big(\cl'_i\big)^\perp=F(\cm_i)^\perp\]
follows from the equalities
above by applying $\cy^{-1}$ to the equality
$Y(\cm_i)^\perp=\cl_i^\perp$ of Lemma~\ref{L3.2.5}(ii).

Now: if we intersect the equality above with
$\fl'(\cs)$, and then take the union over $i\in\nn$,
we obtain the first equality in
\[
\bigcup_{i\in\nn}\fl'(\cs)\cap\big(\cl'_i\big)^\perp
\eq
\bigcup_{i\in\nn}\fl'(\cs)\cap F(\cm_i)^\perp
\eq
\fl'(\cs)\cap\cy^{-1}\big(\fc(\cs)\big)
\]
while the second equality follows from  
\cite[Observation~\ref{O21.-1}]{Neeman18A}.
\eprf

\rmd{R94.-7}
We remind the reader that, in
\cite[Corollary~\ref{C21.15}]{Neeman18A},
we proved that the functor $\cy:\ct\la\MMod\cs$ restricts
to an equivalence of categories
$\fl'(\cs)\cap\cy^{-1}\big(\fc(\cs)\big)\la\fs(\cs)$.
In the next Lemma we allow ourselves,
in a notational crime, to use the
symbol $\cn_n$ both for the
subcategory $\cn_n\subset\fs(\cs)=\fl(\cs)\cap\fc(\cs)$,
and for its inverse image under the equivalence
$\cy:\fl'(\cs)\cap\cy^{-1}\big(\fc(\cs)\big)\la\fs(\cs)$.
\ermd

We deduce the following.

\lem{L94.-9}
Let the notation be as in Definition~\ref{D93.11}.
The equivalence
$\cy:\fl'(\cs)\cap\cy^{-1}\big(\fc(\cs)\big)\la\fs(\cs)$,
of Reminder~\ref{R94.-7}, identifies the subcategory
$\cn_n\subset\fs(\cs)$ of Definition~\ref{D3.-3}(ii)
with the full subcategory
\[
\bigcup_{i\in\nn}\cl'_n\cap\big(\cl'_i\big)^\perp\sub\fl'(\cs)\ .
\]
\elem

\prf
Lemma~\ref{L93.-1} tells us that $\cy(\cl'_n)=\cl_n$.
Therefore
\[
\cn_n\eq\cl_n\cap\fc(\cs)\eq \cy(\cl'_n)\cap\fc(\cs)\ .
\]
Now consider the equalities
\[
\cy\Big[\cl'_n\cap\cy^{-1}\big(\fc(\cs)\big)\Big]
\eq
\cy(\cl'_n)\cap\fc(\cs)\eq\cn_n
\]
where the first is obvious, and the second comes from
the displayed
equalities before.
Coupling this with the fact that
the map
$\cy:\fl'(\cs)\cap\cy^{-1}\big(\fc(\cs)\big)\la\fs(\cs)$
is an equivalence, see Reminder~\ref{R94.-7},
we can rewrite this (criminally) as
$\cn_n=\cl'_n\cap\cy^{-1}\big(\fc(\cs)\big)$
And now Lemma~\ref{L93.-5} allows us expand
$\cl'_n\cap\cy^{-1}\big(\fc(\cs)\big)$
as the right-hand-side expression in
\[
\cn_n\eq
\bigcup_{i\in\nn}\cl'_n\cap\big(\cl'_i\big)^\perp\ .
\]
\eprf

\rmk{R93.-11}
Let us continue to
commit the notational crime of identifying
$\fl'(\cs)\cap\cy^{-1}\big(\fc(\cs)\big)\subset\fl'(\cs)$
with its image under the functor $\cy$, which
is $\fs(\cs)\subset\fl(\cs)$. After all, on the subcategories
the map restricts to an equivalence.

After agreeing to commit this notational crime,
we view $\fs(\cs)$ as full subcategory of
$\fl'(\cs)$ and, if the perpendiculars are
understood in the category $\fl'(\cs)$,
then the formula of Lemma~\ref{L93.-5} simplifies to
\[
\fs(\cs)\eq\bigcup_{i\in\nn}\big(\cl'_i\big)^\perp\ ,
\]
while the formula of Lemma~\ref{L94.-9} remains
\[
\cn_n\eq\cl'_n\cap\bigcup_{i\in\nn}\big(\cl'_i\big)^\perp\ .
\]
So far in this section our metrics have only been good,
there has, so far, been no mention of excellence of
metrics.
\ermk

Recall: starting in Section~\ref{S70} we have
mostly confined our
attention to excellent metrics---as the name suggests
they have special properties that render them more
useful than the average good metric. For the highlights
the reader is referred to Proposition~\ref{P97.1},
Theorem~\ref{T97.9} and Remark~\ref{R97.11}.

Of course, it isn't surprising that,
if the metric $\{\cm_i,\,i\in\nn\}$ on the
triangulated category $\cs$ happens to be excellent,
then good 
extensions $F:\cs\la\ct$ will also have unexpectedly
nice properties. But for the key property we need,
there is no need to assume the full excellence
of the metric. Hence we introduce
an intermediate concept.

\dfn{D93.1}
Let $\cs$ be a triangulated category with a
good metric $\{\cm_i\mid i\in\nn\}$.
The metric is called \emph{very good} if
the following holds:
\be
\item
$\cs=\cup_{i\in\nn}\,{^\perp\cm_i}$.
\item
For every integer $m\in\nn$ there
exists an integer $n>m$ such that any object
$F\in\cs$ admits, in $\cs$,
a triangle $E\la F\la D\la\T E$ with $E\in {^\perp\cm_n}$
and with $D\in\cm_m$.
\ee
\edfn

\rmk{R93.3}
If the reader compares Definitions~\ref{D93.1}
and \ref{D3.3}, it becomes obvious that
every excellent metric is very good. A very good
metric is excellent if, in addition to the
hypotheses postulated in Definitions~\ref{D93.1},
it satisfies the extra condition of
Definition~\ref{D3.3}(iii).

As an aside, we note that Definition~\ref{D3.3}(iii)
is complicated to state. So much so
that we needed to spend the whole of
Section~\ref{S3} preparing the ground.
\ermk

The following little observation will be helpful.

\lem{L93.7}
Let $\cs$ be a triangulated category, assume
$\{\cm_i\mid i\in\nn\}$ is a good metric
on the category $\cs$,
and suppose that $F:\cs\la\ct$
is a good extension with respect to the metric.

Suppose  $x_*$ is a Cauchy sequence in $\cs$, and 
$X\in\ct$ is the object 
$X=\hoco F(x_*)$. For every integer $n>0$, there
exists an integer $N>0$ such that, if $i\geq N$, then
in the triangle $F(x_i)\la X\la d_{i}$ we
have
$d_i\in F\big({^\perp\cm_n}\big)^\perp$.

\elem

\prf
Since $x_*$ is a Cauchy sequence in $\cs$,
for every integer $n>0$ there
exists an integer $N>0$ such that, if $j>i\geq N$, then
in the triangle $x_i\la x_j\la d_{i,j}\la\T x_i$ we
have $d_{i,j}\in\cm_n$. Therefore,
for any object $y\in{^\perp\cm_n}$, we have that 
$\Hom_\cs^{}(y,-)$ takes
the map $x_i\la x_j$ to an epimorphism and
the map $\T x_i\la\T x_j$ to a monomorphism.
As $F$ is fully faithful, this means that
$\Hom_\ct^{}\big(F(y),-\big)$
take $F(x_i)\la F(x_j)$
to an epimorphism and $\T F(x_i)\la\T F(x_j)$ to a
monomorphism.
But now $\Hom\big(F(y),-\big)$ takes $X=\hoco F(x_*)$
to $\colim\,\Hom\big(F(y),F(x_j)\big)$, and we conclude
that $\Hom_\ct^{}\big(F(y),-\big)$ must take the map
$F(x_i)\la X$ to an epimorphism and $\T F(x_i)\la \T X$ to
a monomorphism.
The long exact sequence,
obtained by applying the
functor $\Hom_\ct^{}\big(F(y),-\big)$ to
the sequence
\[\xymatrix@C+10pt{
F(x_i)\ar[r] & X\ar[r] & d_i\ar[r] & \T F(x_i)\ar[r] & \T X
}\]
gives that $\Hom_\ct^{}\big(F(y),d_i\big)=0$.
As this is true for all $y\in{^\perp\cm_n}$
and all 
$i\geq N$, the Lemma is proved.
\eprf

We are now ready for the main result.

\pro{P93.9}
Let $\cs$ be a triangulated category, and assume
$\{\cm_i\mid i\in\nn\}$ is a very good metric
on the category $\cs$.
Assume further that, for every integer $i\in\nn$, the
subcategory $\cm_i\subset\cs$ is closed in $\cs$ under
direct summands.
And finally suppose that $F:\cs\la\ct$
is a good extension with respect to the metric.

Then the subcategory $\fl'(\cs)$, of
\cite[Definition~\ref{D21.7}]{Neeman18A},
is a triangulated
subcategory of $\ct$. Moreover:
suppose we are given
\be
\item
A morphism $\ph:A\la B$ in the category $\fl'(\cs)$.
\item
Cauchy sequences $a_*$ and $b_*$
in the category $\cs$, and isomorphisms
$A\cong\hoco F(a_*)$ and $B\cong\hoco F(b_*)$.
These must exist by Reminder~\ref{R93.5}(iii).

For future reference we write $\alpha:F(a_*)\la A$
for the map in $\ct$ from the Cauchy sequence to
$A$, and similarly for $\beta:F(b_*)\la B$.
\item
Combining Reminder~\ref{R93.5}(iii)
with
\cite[Lemma~\ref{L28.19}(i)]{Neeman18A}
we may, after passing to subsequences
of $a_*$ and $b_*$, choose
a map of Cauchy sequences $\ph_*:a_*\la b_*$
with $\colim\,Y(\ph_*)=\cy(\ph)$.
\setcounter{enumiv}{\value{enumi}}
\ee
Let $A\stackrel\ph\la B\la C\la\T A$ be a triangle
in $\ct$.
Then, after passing to subsequences,
we may complete $\ph_*$ in $\cs$ to a Cauchy sequence
of triangles $a_*\stackrel{\ph_*}\la b_*\la c_*\la\T a_*$,
and do it in such a way that
\be
\setcounter{enumi}{\value{enumiv}}
\item
The commutative square
\[\xymatrix@C+20pt{
F(a_*) \ar[d]^\alpha\ar[r]^-{F(\ph_*)} & F(b_*)\ar[d]^\beta \\
A\ar[r]^-\ph & B
}\]
may be extended, in $\ct$, to
a commutative diagram
\[\xymatrix@C+20pt{
F(a_*) \ar[d]^\alpha\ar[r]^-{F(\ph_*)} & F(b_*)\ar[r]\ar[d]^\beta &
 F(c_*)\ar[d]^\gamma \ar[r] &\T F(a_*)\ar[d]^{\T\alpha} \\
A\ar[r]^-\ph & B\ar[r] & C\ar[r] &\T A
}\]
\item
There is an isomorphism
$\hoco F(c_*)\la C$ in the category $\ct$,
which exhibits the pre-triangle
$\cy(A)\stackrel{\cy(\ph)}\la\cy(B)\la\cy(C)\la\T\cy(A)$
as naturally isomorphic to
the strong triangle
obtained by taking
the colimit in $\MMod\cs$ of the sequence 
$Y(a_*)\stackrel{Y(\ph_*)}\la Y(b_*)\la Y(c_*)\la\T Y(a_*)$.
\setcounter{enumiv}{\value{enumi}}
\ee
\epro

\prf
We prove the ``moreover'' assertion,
which means we may assume
given all the information
chosen in (i), (ii) and (iii).
By (iii) we are given a map
of sequences $\ph_*:a_*\la b_*$,
and the first step is to choose our subsequence.
We are about to produce an increasing function
$\rho:\nn\la\nn$, and replace
the cauchy sequences 
$a_1^{}\la a_2^{}\la a_3^{}\la\cdots$
and $b_1^{}\la b_2^{}\la b_3^{}\la\cdots$
by the Cauchy sequences
$a_{\rho(1)}\la a_{\rho(2)}\la a_{\rho(3)}\la\cdots$
and $b_{\rho(1)}\la b_{\rho(2)}\la b_{\rho(3)}\la\cdots$

We do this inductively, starting with $\rho(1)=1$.
Now for the inductive step. Suppose $\rho(k)$ has been
constructed; our metric $\{\cm_i\mid i\in\nn\}$
is assumed to be
very good, and Definition~\ref{D93.1}(i) allows
us to
\be
\setcounter{enumi}{\value{enumiv}}
\item
Choose an 
integer $m\in\nn$ such that $^\perp\cm_m$ contains
$a_{\rho(k)}$, $b_{\rho(k)}$, and $\T a_{\rho(k)}$. 
\item
Next: with $m$ as above, set $n$ to be the integer 
$n=m+\rho(k)+1$. Lemma~\ref{L93.7} allows us to choose
an integer $N>n$ such that, for any $i\geq N$, in the triangles
\[
F(a_i)\stackrel{\alpha_i}\la A\la d_i\qquad\text{ and }\qquad
F(b_i)\stackrel{\beta_i}\la B\la\ov d_i
\]
the objects $d_i$, $\ov d_i$ and $\T d_i$ all lie
in $F\big({^\perp\cm_n}\big)^\perp$. Put $\rho(k+1)=N$, with
$N$ chosen as above.
\setcounter{enumiv}{\value{enumi}}
\ee
This completes the inductive definition of the
function $\rho:\nn\la\nn$.

And now we pass to the subsequence, replacing
$a_1^{}\la a_2^{}\la a_3^{}\la\cdots$
and $b_1^{}\la b_2^{}\la b_3^{}\la\cdots$
by the Cauchy subsequences
$a_{\rho(1)}\la a_{\rho(2)}\la a_{\rho(3)}\la\cdots$
and $b_{\rho(1)}\la b_{\rho(2)}\la b_{\rho(3)}\la\cdots$\ \
With our new sequences $a_*$ and $b_*$, the key
features of
(vi) and (vii) are transformed into
\be
\setcounter{enumi}{\value{enumiv}}
\item
For each $k\in\nn$ there is an 
integer $m=m(k)\in\nn$ such that $^\perp\cm_m$ contains
$a_k^{}$, $b_k$, and $\T a_k^{}$,
and an integer $n=n(k)>\max\big(m(k),k+1\big)$
such that, for all
$i>k$,
in the triangles
\[
F(a_i)\stackrel{\alpha_i}\la A\la d_i\qquad\text{ and }\qquad
F(b_i)\stackrel{\beta_i}\la B\la\ov d_i
\]
the objects $d_i$, $\ov d_i$ and $\T d_i$ all lie
in $F\big({^\perp\cm_n}\big)^\perp$.
\setcounter{enumiv}{\value{enumi}}
\ee

Next: by (ii) and (iii),
for every integer $k\in\nn$, we are given
in $\ct$ a commutative
square
\[\xymatrix@C+20pt{
F(a_k^{}) \ar[d]^{\alpha_k^{}}\ar[r]^-{F(\ph_k)} & F(b_k^{})\ar[d]^{\beta_k^{}} \\
A\ar[r]^-\ph & B
}\]
This square may be extended, in the
triangulated category $\ct$, to
a $3\times3$ diagram whose rows and columns are
triangles
\[\xymatrix@C+20pt{
&F(a_k^{}) \ar[d]^{\alpha_k^{}}\ar[r]^-{F(\ph_k^{})} & F(b_k^{})\ar[r]\ar[d]^{\beta_k^{}} &
F(c_k^{})\ar[d]^{\gamma_k^{}} \ar[r] &\T F(a_k^{})\ar[d]^{\T\alpha_k^{}} &\\
(\dagger)_k&A\ar[r]^-\ph\ar[d] & B\ar[d]\ar[r] &
  C\ar[r]\ar[d] &\T A\ar[d]&\\
& d_k^{}\ar[r] & \ov d_k^{}\ar[r] &  
  \wh d_k^{}\ar[r] & \T d_k^{} &
}\]
Choose and fix such a $3\times3$ diagram $(\dagger)_k$
for each $k\in\nn$.
Now recall: (viii) chooses for us
integers $m(k)$ and $n(k)$ such that
\[
\{a_k^{},\,b_k,\,\T a_k^{}\}\subset{^\perp\cm_{m(k)}}\qquad
\text{ and }\qquad
\{d_{k+1}^{},\,\ov d_{k+1},\,\T d_{k+1}^{}\}\subset
F\big({^\perp\cm_{n(k)}}\big)^\perp
\]
The triangles
\[
b_{k}^{}\la c_{k}^{}\la\T  a_{k}^{}
\qquad
\text{ and }\qquad
\ov d_{k+1}^{}\la \wh d_{k+1}^{}\la\T\ov  d_{k+1}^{}\ ,
\]
from the top row of $(\dagger)_k^{}$
and the bottom row of $(\dagger)_{k+1}^{}$,
combine to tell us that
\[
\begin{array}{c}
c_{k}^{}\quad\in\quad{^\perp\cm_{m(k)}}*{^\perp\cm_{m(k)}}
  \sub{^\perp\cm_{m(k)}}
\\
\wh d_{k+1}^{}\quad\in\quad F\big({^\perp\cm_{n(k)}}\big)^\perp*F\big({^\perp\cm_{n(k)}}\big)^\perp\sub F\big({^\perp\cm_{n(k)}}\big)^\perp
\end{array}
\]
Combining this information, with what we already know
from (viii),
we conclude
\be
\setcounter{enumi}{\value{enumiv}}
\item
In the diagram 
$(\dagger)_{k}^{}$, the objects
$a_{k}^{}$, $b_{k}^{}$, $c_{k}^{}$ and $\T a_{k}^{}$
of the top row
all belong to $^\perp\cm_{m(k)}\subset\cs$.
\item
In the diagram $(\dagger)_{k+1}^{}$, the objects
$d_{k+1}^{}$, $\ov d_{k+1}^{}$, $\wh d_{k+1}^{}$ and $\T d_{k+1}^{}$
of the bottom row 
all belong to $F\big({^\perp\cm_{n(k)}}\big)^\perp\subset\ct$.
\setcounter{enumiv}{\value{enumi}}
\ee
Now: for each of
the four columns
$\Tm ({^?d}_{k+1}^{})\la x_{k+1}^{}\la X \la {^?d}_{k+1}^{}$,
where ${^?d}_{k+1}^{}$
stands (as appropriate) for any of
$d_{k+1}^{}$, $\ov d_{k+1}^{}$, $\wh d_{k+1}^{}$ and $\T d_{k+1}^{}$
in the diagram $(\dagger)_{k+1}^{}$, we have
${^?d}_{k+1}^{}\in F\big({^\perp\cm_{n(k)}}\big)^\perp$.
Now since
$\cm_{n(k)-1}$ contains $\cm_{n(k)}\cup\Tm\cm_{n(k)}$ we have that
$F\big({^\perp\cm_{n(k)-1}}\big)^\perp$
contains both $F\big({^\perp\cm_{n(k)}}\big)^\perp$
and $\Tm F\big({^\perp\cm_{n(k)}}\big)^\perp$, and in particular
both $\Tm ({^?d}_{k+1}^{})$ and ${^?d}_{k+1}^{}$ lie in
$F\big({^\perp\cm_{n(k)-1}}\big)^\perp$.
Thus
\be
\setcounter{enumi}{\value{enumiv}}
\item
Let the integer $n(k)>m(k)$ be as in  
(viii).
For any object $y\in{^\perp\cm_{n(k)-1}}$, the functor
$\Hom\big(F(y),-\big)$ takes the four vertical morphisms
\[
\begin{array}{ccc}
\alpha_{k+1}^{}:F(a_{k+1}^{}) \la A\,,&\qquad&
\beta_{k+1}^{}:F(b_{k+1}^{}) \la B\,,\\
\gamma_{k+1}^{}:F(c_{k+1}^{}) \la C\,,&\qquad&
\T\alpha_{k+1}^{}:\T F(a_{k+1}^{}) \la \T A
\end{array}
\]
of the diagram $(\dagger)_{k+1}^{}$ to isomorphisms.
\setcounter{enumiv}{\value{enumi}}
\ee
Now by (viii) the integers $m(k)<n(k)$
satisfy $m(k)\leq n(k)-1$, therefore
$\cm_{n(k)-1}\subset\cm_{m(k)}$ and hence
${^\perp\cm_{m(k)}}\subset{^\perp\cm_{n(k)-1}}$. On the other hand
(ix)
tells us that the three objects
$a_{k}^{}$, $b_{k}^{}$, $c_{k}^{}$
of the top row of $(\dagger)_{k}^{}$
all belong to $^\perp\cm_m\subset{^\perp\cm_{n-1}}$.
Putting this together with (xi) gives
\be
\setcounter{enumi}{\value{enumiv}}
\item
If 
$y\in\cs$ is one of the objects
$a_{k}^{}$, $b_{k}^{}$, $c_{k}^{}$
in the top row of $(\dagger)_{k}^{}$,
and $F(x)\la X$ is one of the four vertical morphisms
of the diagram $(\dagger)_{k+1}^{}$
\[
\begin{array}{ccc}
\alpha_{k+1}^{}:F(a_{k+1}^{}) \la A\,,&\qquad&
\beta_{k+1}^{}:F(b_{k+1}^{}) \la B\,,\\
\gamma_{k+1}^{}:F(c_{k+1}^{}) \la C\,,&\qquad&
\T\alpha_{k+1}^{}:\T F(a_{k+1}^{}) \la \T A\,,
\end{array}
\]
then any map $F(y)\la X$ factors as
$F(y)\stackrel{F(\psi)}\la F(x)\la X$
for a unique morphism $\psi:y\la x$
in the category $\cs$.
\setcounter{enumiv}{\value{enumi}}
\ee
Therefore the three maps
\[
\alpha_{k}^{}:F(a_{k}^{}) \la A\,,\qquad
\beta_{k}^{}:F(b_{k}^{}) \la B\,,\qquad
\gamma_{k}^{}:F(c_{k}^{}) \la C
\]
factor as
\[
F(a_{k}^{}) \stackrel{F(\wt\alpha_{k}^{})}\la F(a_{k+1}^{})\stackrel{\alpha_{k+1}^{}}\la A\,, \qquad\qquad 
F(b_{k}^{}) \stackrel{F(\wt\beta_{k}^{})}\la F(b_{k+1}^{})\stackrel{\beta_{k+1}^{}}\la B \,,
\]
\[
F(c_{k}^{}) \stackrel{F(\wt\gamma_{k}^{})}\la F(c_{k+1}^{})\stackrel{\gamma_{k+1}^{}}\la C 
\]
for a unique triple of morphisms 
\[
\wt\alpha_{k}^{}:a_{k}^{} \la a_{k+1}^{}\,,\qquad
\wt\beta_{k}^{}:b_{k}^{} \la b_{k+1}^{}\,,\qquad
\wt\gamma_{k}^{}:c_{k}^{} \la c_{k+1}^{}
\]
in the category $\cs$.

Next let us study the diagram below for commutativity
\[\xymatrix@C+10pt{
&F(a_{k}^{}) \ar[d]^{F(\wt\alpha_{k}^{})}\ar[r]^-{F(\ph_{k}^{})} & F(b_{k}^{})\ar[r]\ar[d]^{F(\wt\beta_{k}^{})} &
F(c_{k}^{})\ar[d]^{F(\wt\gamma_{k}^{})} \ar[r] &\T F(a_{k}^{})\ar[d]^{\T F(\wt\alpha_{k}^{})} &\\
(\dagger\dagger)_{k}&F(a_{k+1}^{}) \ar[d]^{\alpha_{k+1}^{}}\ar[r]^-{F(\ph_{k+1}^{})} & F(b_{k+1}^{})\ar[r]\ar[d]^{\beta_{k+1}^{}} &
F(c_{k+1}^{})\ar[d]^{\gamma_{k+1}^{}} \ar[r] &\T F(a_{k+1}^{})\ar[d]^{\T\alpha_{k+1}^{}} &\\
&A\ar[r]^-\ph & B\ar[r] &
  C\ar[r] &\T A&\\
}\]
If we delete the middle row of $(\dagger\dagger)_{k}$
then the resulting
diagram is the top
part of the commutative diagram $(\dagger)_k^{}$,
and if we delete the top row of $(\dagger\dagger)_{k}$
then
we are left with top part of the
commutative diagram $(\dagger)_{k+1}^{}$.
It follows that the composites in
each of the three squares
\[\xymatrix@C+15pt{
a_{k}^{} \ar[d]^{\wt\alpha_{k}^{}}\ar[r]^-{\ph_{k}^{}} & b_{k}^{}\ar[d]^{\wt\beta_{k}^{}} &
b_{k}^{}\ar[r]\ar[d]^{\wt\beta_{k}^{}} &
c_{k}^{}\ar[d]^{\wt\gamma_{k}^{}} &
c_{k}^{}\ar[d]^{\wt\gamma_{k}^{}} \ar[r] &\T a_{k}^{}\ar[d]^{\T \wt\alpha_{k}^{}}\\
a_{k+1}^{} \ar[r]^-{\ph_{k+1}^{}} & b_{k+1}^{} &
 b_{k+1}^{}\ar[r] &
c_{k+1}^{} &
c_{k+1}^{} \ar[r] &\T a_{k+1}^{}
}\]
delivers a pair of maps, rendering equal the composites
\[\xymatrix@C+3pt@R-15pt{
 F(a_k^{})\ar@<0.5ex>[r] \ar@<-0.5ex>[r] & F(b_{k+1}^{})\ar[r]^-{\beta_{k+1}^{}} & B  & 
 F(b_k^{})\ar@<0.5ex>[r] \ar@<-0.5ex>[r] &  F(c_{k+1}^{})\ar[r]^-{\gamma_{k+1}^{}}  & C \\
& F(c_k^{})\ar@<0.5ex>[r] \ar@<-0.5ex>[r] &\T F(a_{k+1}^{})
 \ar[r]^-{\T\alpha_{k+1}^{}}  & \T A 
}\]
By the uniqueness assertion of (xii) the three squares in
$\cs$ commute, and we deduce that the diagram
\[\xymatrix@C+20pt{
a_{k}^{} \ar[d]^{\wt\alpha_{k}^{}}\ar[r]^-{\ph_{k}^{}} & b_{k}^{}\ar[r]\ar[d]^{\wt\beta_{k}^{}} &
c_{k}^{}\ar[d]^{\wt\gamma_{k}^{}} \ar[r] &\T a_{k}^{}\ar[d]^{\T \wt\alpha_{k}^{}} \\
a_{k+1}^{} \ar[r]^-{\ph_{k+1}^{}} & b_{k+1}^{}\ar[r] &
c_{k+1}^{} \ar[r] &\T a_{k+1}^{} \\
}\]
commutes in $\cs$. 

Let us have a brief recap:
we have, after passing
to a subsequence, completed
$\ph_*:a_*\la b_*$ in $\cs$ to a sequence
of triangles $a_*\stackrel{\ph_*}\la b_*\la c_*\la\T a_*$,
and we have done it in such a way that
the commutative square
\[\xymatrix@C+20pt{
F(a_*) \ar[d]^\alpha\ar[r]^-{F(\ph_*)} & F(b_*)\ar[d]^\beta \\
A\ar[r]^-\ph & B
}\]
was extended, in $\ct$, to
a commutative diagram
\[\xymatrix@C+20pt{
F(a_*) \ar[d]^\alpha\ar[r]^-{F(\ph_*)} & F(b_*)\ar[r]\ar[d]^\beta &
 F(c_*)\ar[d]^\gamma \ar[r] &\T F(a_*)\ar[d]^{\T\alpha} \\
A\ar[r]^-\ph & B\ar[r] & C\ar[r] &\T A
}\]
To prove the assertion
(iv) of the current Proposition,
it remains to show that
\be
\setcounter{enumi}{\value{enumiv}}
\item
The sequence $c_*$ is Cauchy in $\cs$.
\setcounter{enumiv}{\value{enumi}}
\ee
After that, we will have to worry about proving (v) in
the statement of the Proposition.

First we prove (xiii). In the
triangulated category $\cs$ we extend the
morphism $\wt\gamma_k^{}:c_k^{}\la c_{k+1}^{}$
to a triangle
$c_k^{}\stackrel{\wt\gamma_k^{}}\la c_{k+1}^{}\la \delta_k^{}$,
and then in the category $\ct$
we extend the commutative square
\[\xymatrix@C+20pt{
F(c_k^{}) \ar@{=}[d] \ar[r]^-{F(\wt\gamma_k^{})} &
  F(c_{k+1}^{})\ar[d]^{\gamma_{k+1}^{}}\\
F(c_k^{}) \ar[r]^-{\gamma_{k}^{}} & C
}\]
to an octahedron
\[\xymatrix@C+20pt{
F(c_k^{}) \ar@{=}[d] \ar[r]^-{F(\wt\gamma_k^{})} &
 F(c_{k+1}^{})\ar[d]^{\gamma_{k+1}^{}}\ar[r] &
 F(\delta_k^{})\ar[d]\\
F(c_k^{}) \ar[r]^-{\gamma_{k}^{}} & C\ar[d]\ar[r] &
 \wh d_k^{}\ar[d] \\
&\wh d_{k+1}^{}\ar@{=}[r] & \wh d_{k+1}^{}
}\]
By (x) the object $\wh d_{k+1}^{}$ belongs to
$F\big({^\perp\cm_{n(k)}}\big)^\perp$,
and in (viii) we noted that $n(k)\geq k+2$.
Therefore $\wh d_{k+1}^{}\in F\big({^\perp\cm_{k+2}}\big)^\perp$
as long as $k\geq 1$. Thus for $k\geq2$ we have that
$\wh d_{k}^{}\in F\big({^\perp\cm_{k+1}}\big)^\perp$,

But now the inclusion
$\Tm\cm_{k+2}\subset\cm_{k+1}$ implies that
$\Tm F\big({^\perp\cm_{k+2}}\big)^\perp\subset
F\big({^\perp\cm_{k+1}}\big)^\perp$,
and hence
$\Tm\wh d_{k+1}^{},\wh d_{k}^{}$ both
lie in $F\big({^\perp\cm_{k+1}}\big)^\perp$.
The triangle $\Tm\wh d_{k+1}^{}\la F(\delta_k^{})\la\wh d_{k}^{}$
tells us that
$F(\delta_k^{})\in F\big({^\perp\cm_{k+1}}\big)^\perp*F\big({^\perp\cm_{k+1}}\big)^\perp\subset F\big({^\perp\cm_{k+1}}\big)^\perp$, 
and as $F$ is fully faithful we deduce that
$\delta_k^{}\in \big({^\perp\cm_{k+1}}\big)^\perp$.

But now we remember that our metric $\{\cm_i,\,i\in\nn\}$ is
assumed to be very good, and Definition~\ref{D93.1}(ii)
says that, for any integer $m>0$, there exists an integer
$n>m$ such that every object $F\in\cs$ admits a triangle
$E\la F\la D$ with $E\in{^\perp\cm_n}$ and with
$D\in\cm_m$. Applying this to the objects
$F=\delta_k^{}$, with $k+1>n$, we produce triangles
$E_k\la\delta_k^{}\la D_k$ with $E_k\in{^\perp\cm_n}$
and with $D_k\in\cm_m$. However: we know that
$\delta_k^{}\in\big({^\perp\cm_{k+1}}\big)^\perp\subset
\big({^\perp\cm_{n}}\big)^\perp$, and hence the map
from $E_k\in{^\perp\cm_n}$ to
$\delta_k^{}\in\big({^\perp\cm_{n}}\big)^\perp$ must vanish.
Thus $\delta_k^{}$ has to be a direct summand of
$D_k\in\cm_m$, and the hypotheses of the Proposition
give that $\cm_m$ is closed under direct summands.
We conclude that $\delta_k^{}\in\cm_m$ for all
$k+1>n$.
And now
\cite[Lemma~\ref{L20.8}]{Neeman18A}
informs us that the sequence
$c_*$ is Cauchy, proving (xiii) and therefore
completing the proof of (iv).

We have proved (iv)
meaning we now have in $\cs$
a Cauchy sequence of triangles
$a_*\stackrel{\ph_*}\la b_*\la c_*\la\T a_*$,
and in $\ct$ 
a commutative diagram
\[\xymatrix@C+20pt{
F(a_*) \ar[d]^\alpha\ar[r]^-{F(\ph_*)} & F(b_*)\ar[r]\ar[d]^\beta &
 F(c_*)\ar[d]^\gamma \ar[r] &\T F(a_*)\ar[d]^{\T\alpha} \\
A\ar[r]^-\ph & B\ar[r] & C\ar[r] &\T A
}\]
Applying the functor $\cy:\ct\la\MMod\cs$
of
\cite[Notation~\ref{N28.1}]{Neeman18A},
and remembering that $\cy\circ F=Y$,
gives us in the category $\MMod\cs$
a commutative diagram
\[\xymatrix@C+20pt{
Y(a_*) \ar[d]^{\cy(\alpha)}\ar[r]^-{Y(\ph_*)} & Y(b_*)\ar[r]\ar[d]^{\cy(\beta)} &
 Y(c_*)\ar[d]^{\cy(\gamma)} \ar[r] &\T Y(a_*)\ar[d]^{\T\cy(\alpha)} \\
\cy(A)\ar[r]^-{\cy(\ph)} & \cy(B)\ar[r] & \cy(C)\ar[r] &\T \cy(A)
}\]
Of course: in the abelian category we can take
the colimit of the sequence in the top row, and we deduce
a map
\[\xymatrix@C+20pt{
\cy(A) \ar@{=}[d]\ar[r]^-{\cy(\ph)} & \cy(B)\ar[r]\ar@{=}[d] &
 \colim\,Y(c_*)\ar[d]^{\wt\gamma} \ar[r] &\T \cy(A)\ar@{=}[d]\ar[r]^-{\T\cy(\ph)} &\T \cy(B)\ar@{=}[d] \\
 \cy(A)\ar[r]^-{\cy(\ph)} & \cy(B)\ar[r] & \cy(C)\ar[r]
 &\T \cy(A)\ar[r]^-{\T\cy(\ph)} &\T \cy(B)
}\]
Both rows are exact sequences in $\MMod\cs$, and
the 5-lemma allows us to deduce that
the map $\wh\gamma$ is an isomorphism.

But now: we have proved that the sequence $c_*$
is Cauchy, and because $F:\cs\la\ct$ is a
good extension we know that the natural map
$\colim\,Y(c_*)\la\cy\big(\hoco F(c_*)\big)$ is an
isomorphism. The isomorphism $\wt\gamma$ of the diagram
above rewrites as
$\wt\gamma:\cy\big(\hoco F(c_*)\big)\la\cy(C)$,
and
\cite[Lemma~\ref{L21.11}]{Neeman18A}
allows us to find a (non-unique) morphism
$\gamma':\hoco F(c_*)\la C$,
in the category $\ct$, with $\cy(\gamma')=\wt\gamma$.

It remains to prove that $\gamma'$ is an isomorphism.
Complete $\gamma'$ in the category $\ct$ to a triangle
$\hoco F(c_*)\stackrel{\gamma'}\la C\la D$, and as the functor
$\cy:\ct\la\MMod\cs$ is homological we
deduce that $\cy(D)=0$. By the definition of
the functor $\cy$,
back in
\cite[Notation~\ref{N21.-100}]{Neeman18A},
this means that
$D\in\cs^\perp$, and 
\cite[Lemma~\ref{L21.19182}]{Neeman18A}
allows us to deduce that
$D\in{\loc\cs}^\perp$.

On the other hand the objects $A=\hoco F(a_*)$
and $B=\hoco F(b_*)$ clearly belong to
$\loc\cs$, and the triangle $A\la B\la C\la\T A$
tells us that $C\in\loc\cs$. But then,
in the triangle
$\hoco F(c_*)\stackrel{\gamma'}\la C\la D$
we have that both $C$ and $\hoco F(c_*)$ belong
to $\loc\cs$, and we deduce that
$D\in\loc\cs\cap{\loc\cs}^\perp$ must vanish.
Therefore $\gamma':\hoco F(c_*)\la C$ is indeed
an isomorphism in $\ct$.
\eprf

In the beginning of the section, before we introduced
the notion of very good metrics, we studied the
categories $\cl'_i\subset\fl'(\cs)$ of
Definition~\ref{D93.11}. Needless to say the various
formulas we proved back then, for good metrics, also
hold for very good metrics. For the key formulas see
Remark~\ref{R93.-11}. It is nonetheless worth noting:

\rmk{R93.012348}
Let the assumptions be as in Proposition~\ref{P93.9}.
If we agree to continue to commit the notational
crime of Remark~\ref{R93.-11}, that is view
$\fs(\cs)$ as a full subcategory of $\fl'(\cs)$,
then the inclusion $\fs(\cs)\la\fl'(\cs)$ is a
triangulated functor.

After all: Proposition~\ref{P93.9} says the
$\fl'(\cs)$ is a triangulated subcategory
of $\ct$, while
\cite[Theorem~\ref{T20.17}]{Neeman18A}
gives the non-criminal
version of the statement that $\fs(\cs)\subset\ct$ is
also a triangulated subcategory.
\ermk

The following is an easy consequence of
Proposition~\ref{P93.9}.

\cor{C93.13}
Let the assumptions be as in Proposition~\ref{P93.9}.
We remind the reader:
$\cs$ is a triangulated category, and 
$\{\cm_i\mid i\in\nn\}$ is a very good metric
on the category $\cs$.
We assume further that, for every integer $i\in\nn$, the
subcategory $\cm_i\subset\cs$ is closed in $\cs$ under
direct summands.
And finally we suppose that $F:\cs\la\ct$
is a good extension with respect to the metric.

Then, for every integer $i>0$, we have that
$\cl'_i*\cl'_i\subset\cl'_i$.
\ecor

\prf
Suppose $A\la B\la C\stackrel\ph\la\T A$ is a triangle in $\ct$,
with $A$ and $C$ both in $\cl'_i$. We apply
Proposition~\ref{P93.9}
to the morphism $\Tm\ph:\Tm C\la A$.

First of all, because $A$ and $C$ both lie in $\cl'_i$
we may choose Cauchy sequences $a_*$ and $c_*$, both
contained in $\cm_i$, with $A=\hoco F(a_*)$ and
with $C=\hoco F(c_*)$. After passing to
subsequences, Proposition~\ref{P93.9}
creates for us a Cauchy sequence of
triangles $\Tm c_*\stackrel{\Tm\ph_*}\la a_*\la b_*\la c_*$
such that $B\cong\hoco F(b_*)$.
But for each integer $n>0$, the triangle
$a_n\la b_n\la c_n$ tells
us that $b_n\in\cm_i*\cm_i\subset\cm_i$, proving
that $B$ belongs to $\cl'_i$.
\eprf

Proposition~\ref{P93.9} shows that, for a triangulated
category $\cs$ with a very good metric
$\{\cm_i,\,i\in\nn\}$, and assuming
that the subcategories $\cm_i\subset\cs$ are closed in $\cs$
under direct summands, any good
extension $F:\cs\la\ct$ has the property that
the intermediate category $\fl'(\cs)$ is also
triangulated.
It is natural to wonder if the
embedding $\cs\la\fl'(\cs)$ is also a good extension
with respect to the metric.

This does happen, and as we will see there are
useful examples out there.

\cor{C93.95}
Let $\cs$ be a triangulated category, and assume
$\{\cm_i,\,i\in\nn\}$ is a good metric on
$\cs$. Assume further that the category
$\cs$ has the property that, for any object
$E\in\cs$, the coproduct of $\coprod_{n\in\nn}E$ exists
in $\cs$. By this we mean: the coproduct of
countably many copies, of the single object
$E\in\cs$ exists in $\cs$.

Suppose further that
\be
\item
For every integer $i>0$ there exists an integer
$r=r(i)>0$ such that,
for every object $E\in\cm_{i+r}$,
the coproduct $\coprod_{n\in\nn}E$ belongs
to $\cm_i$.
\setcounter{enumiv}{\value{enumi}}
\ee
Then, for any Cauchy sequence $E_*$ in $\cs$,
the coproduct $\coprod_{i=1}^\infty Y(E_i)$, which
we may form in $\MMod\cs$, actually belongs
to $\fl(\cs)\subset\MMod\cs$.

And now for a result about $\fl'(\cs)$.
Assume in addition
\be
\setcounter{enumi}{\value{enumiv}}
\item
The metric $\{\cm_i,\,i\in\nn\}$ on the category
$\cs$ is very good (not only good),
and the subcategories $\cm_i\subset\cs$
are closed in $\cs$ under direct summands.
\item
We are given a good extension $F:\cs\la\ct$
with respect to the metric, and $F$ respects
the coproducts in $\cs$
of the form
$\coprod_{n\in\nn}E$,
whose existence is guaranteed in (i).
\ee
Then, for any Cauchy sequence $E_*$ in $\cs$,
the coproduct $\coprod_{i=1}^\infty F(E_i)$, which
we may form in $\ct$, actually belongs
to $\fl'(\cs)\subset\ct$. Moreover: the functor
$\cy:\fl'(\cs)\la\fl(\cs)$ respects these coproducts.

And finally, the embedding $\cs\la\fl'(\cs)$ is
a good extension with respect to the metric.
\ecor

\prf
Let $E_*$ be a Cauchy sequence in $\cs$, and out
of it form the sequence $\wt E_n$ whose
$n\mth$ term is
$\oplus_{i=1}^{n-1}E_i\oplus\coprod_{i=n}^\infty E_n$.
And let the map $\wt E_n\la \wt E_{n+1}$ be the morphism
\[\xymatrix@C+30pt{
\ds\left(\bigoplus_{i=1}^{n}E_i\right)
\oplus\left(\coprod_{i=n+1}^\infty E_n\right)
\ar[r]^-{\id\oplus f} &
\ds\left(\bigoplus_{i=1}^{n}E_i\right)
\oplus\left(\coprod_{i=n+1}^\infty E_{n+1}\right)
}\]
where $f$ stands for the countable
coproduct of the map $E_n\la E_{n+1}$.
By (i) the sequence $\wt E_*$ is Cauchy in $\cs$,
and hence $\colim\, Y\big(\wt E_*\big)$ belongs to
$\fl(\cs)\subset\MMod\cs$.

But $Y\big(\wt E_*\big)$ is a direct sum
of sequences of the form
$Y(E_1)\la Y(E_2)\la Y(E_3)\la\cdots\la Y(E_n)$,
each of which stabilizes to some $Y(E_n)$
after a finite number of steps.
Hence $\colim\,Y\big(\wt E_*\big)=\coprod_{i=1}^\infty Y(E_i)$.

Now for the rest of the proof, where we add the
hypotheses (ii) and (iii).
By (iii) the functor
$F$ respects those coproducts in $\cs$
whose existence
is guaranteed in (i), and we have
that $F\big(\wt E_n\big)
=\oplus_{i=1}^{n-1}F(E_i)\oplus\coprod_{i=n}^\infty F(E_n)$. And
the map $F\big(\wt E_n\big)\la F\big(\wt E_{n+1}\big)$
is equal to
\[\xymatrix@C+30pt{
\ds\left(\bigoplus_{i=1}^{n}F(E_i)\right)
\oplus\left(\coprod_{i=n+1}^\infty F(E_n)\right)
\ar[r]^-{\id\oplus F(f)} &
\ds\left(\bigoplus_{i=1}^{n}F(E_i)\right)
\oplus\left(\coprod_{i=n+1}^\infty F(E_{n+1})\right)\ .
}\]
Since this is a coproduct of sequences that stabilize,
its homotopy colimit is a genuine colimit. We have
exhibited
$\coprod_{i=1}^\infty F(E_i)$ as $\hoco F\big(\wt E_*\big)$
for the Cauchy sequence $\wt E_*$ in $\cs$.
This shows that the object $\coprod_{i=1}^\infty F(E_i)$,
formed in the category $\ct$, belongs to
the subcategory $\fl'(\cs)$.

Inspection of the last paragraphs tells
us that the process by which we formed these coproducts,
in the categories $\fl'(\cs)$ and in $\fl(\cs)$
respectively, is compatible with the functor
$\cy:\fl'(\cs)\la\fl(\cs)$.

By
Proposition~\ref{P93.9} the subcategory
$\fl'(\cs)\subset\ct$ is triangulated, and hence
the homotopy colimit $\hoco F(E_*)$, formed
in $\ct$, actually lies
in $\fl'(\cs)$. The fact that this homotopy
colimit satisfies the requirement
of
\cite[Definition~\ref{D21.1}(ii)]{Neeman18A} 
immediately follows.
\eprf

\section{Good extensions of excellent metrics}
\label{S94}

Starting from Section~\ref{S70},
we have spent much effort studying excellent metrics.
That is we fixed a triangulated category
$\cs$ together with an excellent  
metric $\{\cm_i,\,i\in\nn\}$, and
made an elaborate study of various
completions that can be constructed out of it.
Excellent metrics have
striking properties that set them apart from
run-of-the-mill good metrics.

And then, in Section~\ref{S93},
we went back to take another look at good extensions
with respect to metrics:
that is the starting point is no longer
a single triangulated category $\cs$
together with a choice of metric,
but we add to the baggage a fully faithful
functor $F:\cs\la\ct$ satisfying useful
extra properties. Those were introduced
back in
\cite[Section~\ref{S21}]{Neeman18A},
where it was
shown that they can be useful in studying
completions for good metrics.
In Section~\ref{S93} we took another look
at good extensions, this time assuming
the metric is ``very good''---that is better than
merely ``good'', but not quite as
well behaved as ``excellent''.
Being very good sufficed
to show that a number of the
constructions of
\cite[Section~\ref{S21}]{Neeman18A}
become superior to average.

\plm{P94.1}
In this section we address two questions:
\be
\item
Let $\cs$ be a triangulated category, and  
assume $\{\cm_i,\,i\in\nn\}$ is a very good
metric on $\cs$. Assume further that 
$F:\cs\la\ct$ is a good extension with respect
to the metric.

Can one detect, in $\fl'(\cs)$, that
the metric $\{\cm_i,\,i\in\nn\}$ is in fact
excellent?
\item
In Proposition~\ref{P3.94949} we showed that,
for excellent metrics $\{\cm_i,\,i\in\nn\}$,
the functor $\wh Y:\big(\MMod\cs\big)\op\la\MMod{\cs\op}$
restricts to an equivalence
$\fl(\cs)\op\cong\fl\big(\fs(\cs)\big)\op$,
and for every integer $\ell>0$
this further restricts to equivalences
$\cl\op_\ell\cong\wh\cl_\ell$.

Are there analoguous statements for the subcategories
$\fl'(\cs)\subset\ct$ of
\cite[Definition~\ref{D21.7}]{Neeman18A}
and its subcategories $\cl'_\ell$ of
Definition~\ref{D93.11}?
\ee
\eplm

Before we launch into this, a little elaboration
of Problem~\ref{P94.1}(i) might be in order.

Suppose we are given a
triangulated category
$\cs$ together with 
good metric $\{\cm_i,\,i\in\nn\}$.
For this metric to be very good means that
we impose, on the metric $\{\cm_i,\,i\in\nn\}$,
two additional restrictions, see
Definition~\ref{D93.1}. And the reader can
check that these properties are statement
purely about the objects and triangles in $\cs$.
Next: for the very good metric to be
excellent, it must satisfy
the extra hypothesis of
Definition~\ref{D3.3}(iii), and this is
a statement about the existence of
morphisms in $\fl(\cs)$ of the form
$E\la D$, with $E\in\cs$ and $D\in\fs(\cs)$.
Merely formulating
Definition~\ref{D3.3}(iii)
requires introducing auxiliary constructions---OK,
the auxiliary constructions in question never
use anything beyond the triangulated
category $\cs$, its good metric $\{\cm_i,\,i\in\nn\}$,
and the Yoneda embedding $Y:\cs\la\MMod\cs$.
But it is still natural to ask if,
in the presence of a good extension $F:\cs\la\ct$,
there is a way to reformulate
Definition~\ref{D3.3}(iii)
as a statement about the existence of
morphisms $E\la D$ in $\fl'(\cs)\subset\ct$.

Not surprisingly the answer to Problem~\ref{P94.1}(i)
is a simple and almost
unqualified Yes. We formulate this in

\pro{P94.3}
Let $\cs$ be a triangulated category, and  
assume $\{\cm_i,\,i\in\nn\}$ is a very good
metric on $\cs$. Suppose that 
$F:\cs\la\ct$ is a good extension with respect
to the metric. Assume further that each
of the subcategories $\cm_i$ is closed in $\cs$
under direct summands.

The metric $\{\cm_i,\,i\in\nn\}$
is excellent if and only if,
for every integer $m>0$, there exists an
integer $n>0$ such that every object $B\in\cs$
admits, in the category $\fl'(\cs)\subset\ct$,
a distinguished triangle $A\la F(B)\la C$ with
$C\in\fs(\cs)\cap\big(\cl'_n\big)^\perp$ and with
$A\in\cl'_m$.
\epro

Let us note that, in the statement of the Proposition
as in its proof, we commit the notational abuse
of Remark~\ref{R93.-11}: the category
$\fs(\cs)\subset\fl(\cs)$
is being identified with the subcategory
$\fl'(\cs)\cap\cy^{-1}\big(\fc(\cs)\big)\subset\fl'(\cs)$.
They are of course equivalent: the functor
$\cy:\fl'(\cs)\la\fl(\cs)$ restricts to a triangulated
equivalence on the subcategories.

\prf
Let us begin by proving the direction $\Longrightarrow$,
meaning we assume the metric on $\{\cm_i,\,i\in\nn\}$
is excellent and produce from $m$ the integer $n$,
satisfying the existence statement in the Proposition.

The recipe for the integer $n$ is easy: we begin with an
integer $m>0$. Let $n>m+1$ be the integer
that we obtain, by applying Definition~\ref{D3.3}(iii) to
the integer $m+1$.
By this we mean:
for any object $B\in\cs$, there exists in $\cs$
a Cauchy sequence $c_*$, such that
$C'=\colim\, Y(c_*)$ belongs to $\fs(\cs)\cap\cl_{n}^\perp$, 
and in the category $\fl(\cs)$ a morphism
$\ph':Y(B)\la C'$ which is type-$(m+1)$
with respect to $(B, c_*)$. By
\cite[Lemma~\ref{L28.19}(i)]{Neeman18A} we
may, after replacing $(B,c_*)$ by subsequences,
choose a map of Cauchy sequences
$\ph_*:B\la c_*$ with $\ph'=\colim\,Y(\ph_*)$.

Now let $C=\hoco F(c_*)$, and let
$\ph:F(B)\la C$ be the composite
$F(B)\la F(c_1^{})\la\hoco F(c_*)$.
Applying the functor $\cy$ gives
that $\cy(\ph)$ is the composite
$Y(B)\la Y(c_1^{})\la\colim\, Y(c_*)=C'$,
which is the map $\ph':Y(B)\la C'$ above.
For future reference note that
we therefore know that $\cy(\ph):Y(B)\la\cy(C)$
is of type-$(m+1)$ with respect to $(B,c_*)$,
and that $\cy(C)\in\fs(\cs)\cap\cl_{n}^\perp$.

Next: the object $C=\hoco F(c_*)$ clearly belongs
to $\fl'(\cs)$. But it also lies in
$\cy^{-1}\big(\cl_{n}^\perp\big)=\big(\cl'_{n}\big)^\perp$
where the equality is by Lemma~\ref{L93.-3}.
Thus
\[
C\quad\in\quad \fl'(\cs)\cap \big(\cl'_{n}\big)^\perp
\eq \fs(\cs)\cap \big(\cl'_{n}\big)^\perp
\]
where the equality is in the abuse of notation
that we adopted in Remark~\ref{R93.-11}, and
follows from the first display formula in
Remark~\ref{R93.-11}.

This puts us in the position where
Proposition~\ref{P93.9}~(i), (ii) and(iii) are
safisfied.
Recalling that $\cy(\ph)$ is
a type-$(m+1)$ morphism with respect
to $(B,c_*)$,
the map of Cauchy sequences $\ph_*:B\la c_*$ that
we chose must be of type-$(m+1)$. 
Now in the category $\ct$ complete $\ph$ to a triangle
$A\la F(B)\stackrel\ph\la C\la\T A$.
To finish the proof of the direction
$\Longrightarrow$, it suffices to
show that $A\in\cl'_m$. 

In the category $\cs$ we have produced a
type-$(m+1)$ morphism
of Cauchy sequences $\ph_*:B\la c_*$, and in
the category $\ct$ a commutative square
\[\xymatrix@C+20pt{
F(B) \ar@{=}[d]\ar[r]^-{F(\ph_*)} & F(c_*)\ar[d]^\gamma \\
F(B)\ar[r]^-\ph & C
}\]
Proposition~\ref{P93.9}~(iv) and (v)
permit us, after replacing $\ph_*:B\la c_*$ by
a subsequence, to 
extend $\ph_*:B\la c_*$ in the category $\cs$
to a Cauchy sequence of
triangles $a_*\la B\stackrel{\ph_*}\la c_*\la\T a_*$,
and do it in such a way that we obtain,
in $\ct$,
a commutative diagram
\[\xymatrix@C+20pt{
F(a_*) \ar[d]^\alpha\ar[r] & F(B)\ar[r]^-{F(\ph_*)}\ar@{=}[d] &
 F(c_*)\ar[d]^\gamma \ar[r] &\T F(a_*)\ar[d]^{\T\alpha} \\
A\ar[r] & F(B)\ar[r]^-\ph & C\ar[r] &\T A
}\]
where we furthermore know that
there is an isomorphism
$\hoco F(a_*)\la A$ in the category $\ct$.
But because the morphism $\ph_*:B\la c_*$ is
of type-$(m+1)$, Definition~\ref{D3.341}
tells us that the Cauchy sequence $\T a_*$
must satisfy $\T a_i\in\cm_{m+1}$ for all $i\gg0$.
But then $a_i\in\cm_m$ for all $i\gg0$,
and $A=\hoco F(a_*)$ must belong to $\cl'_m$.

We have proved the direction $\Longrightarrow$,
and it remains to prove the direction
$\Longleftarrow$.

Assume therefore that the existence statement of
the current Proposition is satisfied, and
start with an integer $m>0$. Out of
it we need to produce an integer $n>m$ satisfying
the existence statement
in Definition~\ref{D3.3}(iii). Take the integer
$m>0$, then  plug the integer $m+1$
into the machine of Proposition~\ref{P94.3},
and an integer $n>m+1$ pops out. And what the assertion
of the Proposition states is that, for any object
$B\in\cs$, there exists in $\ct$ a triangle $A\stackrel\ph\la F(B)\stackrel\psi\la C$
with $A\in\cl'_{m+1}$ and with
$C\in\fs(\cs)\cap\big(\cl'_n\big)^\perp$.

Because $A$ belongs to $\cl'_{m+1}$, there must
exist in $\cm_{m+1}$ a Cauchy sequence $a_*$
with $A=\hoco F(a_*)$. Now the
functor $F:\cs\la\ct$ is fully faithful,
hence for each integer $i\in\nn$
the composite $F(a_i)\la A\stackrel\ph\la F(B)$ must
be $F(\ph_i)$ for a unique morphism $\ph_i:a_i\la B$
in the category $\cs$. We have produced, in $\cs$,
a morphism of Cauchy sequences $\ph_*:a_*\la B$, and in
$\ct$ we have a commutative square
\[\xymatrix@C+20pt{
F(a_*) \ar[d]^\alpha\ar[r]^-{F(\ph_*)} & F(B)\ar@{=}[d] \\
A\ar[r]^-\ph & F(B)
}\]
Proposition~\ref{P93.9}~(iv) and (v)
permit us, after replacing $\ph_*:a_*\la B$ by
a subsequence, to 
extend $\ph_*:a_*\la B$ in the category $\cs$
to a Cauchy sequence of
triangles $a_*\stackrel{\ph_*}\la B\stackrel{\psi_*}\la c_*\la\T a_*$,
and do it in such a way that we obtain,
in $\ct$,
a commutative diagram
\[\xymatrix@C+20pt{
F(a_*) \ar[d]^\alpha\ar[r]^-{F(\ph_*)} & F(B)\ar[r]^{F(\psi_*)}\ar@{=}[d] &
 F(c_*)\ar[d]^\gamma \ar[r] &\T F(a_*)\ar[d]^{\T\alpha} \\
A\ar[r]^-\ph & F(B)\ar[r]^-\psi & C\ar[r] &\T A
}\]
where the
functor $\cy:\ct\la\MMod\cs$ takes
the distinguished triangle
$A\stackrel\ph\la F(B)\stackrel\psi\la C\la\T A$
to the strong triangle, which is the colimit
of $Y$ applied to the Cauchy sequence
$a_*\stackrel{\ph_*}\la B\stackrel{\psi_*}\la c_*\la\T a_*$.

By construction the map of sequences $\psi_*:B\la c_*$,
which is part of the Cauchy sequence of triangles
$B\stackrel{\psi_*}\la c_*\la \T a_*$, satisfies
$\T a_*\subset\T\cm_{m+1}\subset\cm_m$. Hence
the morphism $Y(B)=\cy F(B)\stackrel{Y(\psi)}\la\cy(C)$,
in the category
$\fl(\cs)\subset\MMod\cs$, is of type-$m$ with respect
to $(B,c_*)$. And we also know that $\cy(C)=\colim\,Y(c_*)$
belongs to
\[
\cy\Big[\big(\cl'_n\big)^\perp\Big]\eq
\cy\Big[\cy^{-1}\big(\cl_n^\perp\big)\Big]\sub\cl_n^\perp
\]
where the equality by applying $\cy$ to the
equality 
$\big(\cl'_n\big)^\perp=\cy^{-1}\big(\cl_n^\perp\big)$
of Lemma~\ref{L93.-3}, and the inclusion is
obvious.
Clearly $\cy(C)=\colim\,Y(c_*)$ also belongs to
$\fl(\cs)$, making it an object of
$\fl(\cs)\cap\cl_n^\perp=
\fs(\cs)\cap\cl_n^\perp$, where the equality is immediate
from Remark~\ref{R3.2??}.
Thus our morphism $\cy(\psi):Y(B)\la\cy(C)$
satisfies the requirements of Definition~\ref{D3.3}(iii).
\eprf

Proposition~\ref{P94.3} solves for us the recognition
issue, formulated in Problem~\ref{P94.1}(i). We remind
the reader. Suppose we are given a triangulated category
$\cs$ with a very good metric $\{\cm_i,\,i\in\nn\}$.
Assume that we are also given a good extension $F:\cs\la\ct$
with respect to the metric. Can the good extension
help us recognize that the metric isn't only very good,
but also satisfies the condition for excellence?

The short answer is Yes, provided that the the subcategories
$\cm_i$ are closed in $\cs$ under direct summands.

The time has now come to switch our attention
to how the theory, developed for excellent metrics
in the absence of a good extension, can be elaborated when
a good extension is present.
We begin with the following.

\obs{O94.5323}
Let $\cs$ be a triangulated category, and assume
$\{\cm_i,\,i\in\nn\}$ is an
excellent metric on $\cs$.
Let $F:\cs\la\ct$ be a good extension with respect to
the metric.
And, as we have done since Remark~\ref{R93.-11},
we commit the notational sin of using the equivalence
of categories $\cy$ to identify
$\fl'(\cs)\cap\cy^{-1}\big(\fc(\cs)\big)$, which is
a triangulated subcategory of $\ct$, with
$\fs(\cs)$, which is a subcategory of
$\MMod\cs$.
Now consider the diagram
\[\xymatrix@R-15pt@C+20pt{
\fl'(\cs)\op\ar[r]\ar[dd] &
\ct\op\ar[rd]^{\cy'}\ar[dd]_{\cy\op} &
 \\
& & \MMod{\fs(\cs)\op}\\
\fl(\cs)\op\ar[r] &
\big(\MMod\cs\big)\op\ar[ru]_-{\wh Y} &
}\]
The horizontal functors in the diagram are
inclusions. The functor labeled
$\cy:\ct\la\MMod\cs$ comes
from
\cite[Notation~\ref{N21.-100}]{Neeman18A}
applied to the fully faithful
functor $F:\cs\la\ct$, the functor
labeled $\cy'$  comes
from
\cite[Notation~\ref{N21.-100}]{Neeman18A}
applied to the fully faithful
functor $\wt F:\fs(\cs)\op\la\ct\op$,
and the functor labeled $\wh Y$
was introduced in
Discussion~\ref{D3.6.5}.
And finally: the vertical map on the left
is the map making the square commute,
that is the restriction to $\fl'(\cs)$
of the map $\cy:\ct\la\MMod\cs$.

Given that all the maps are either restricted
Yoneda maps or inclusions, the reader should
be warned that
\emph{the triangle on the right need not commute.}
It comes down to asking whether, for $E\in\ct$ and
all $X\in\fs(\cs)$, the natural map
$\Hom_\ct^{}(E,X)\la\Hom\big(\cy(E),\cy(X)\big)$
is an isomorphism---there is no reason to
expect this. But in
\cite[Lemma~\ref{L21.13}]{Neeman18A} we saw that,
if $E$ is restricted to belong to $\fl'(\cs)\subset\ct$.
then the map is an isomorphism. In other words:
the pentagon forming the perimeter of the
diagram does lead to equal composites.
\eobs

We sum this up in

\lem{L94.5325}
Let $\cs$ be a triangulated category, and assume
$\{\cm_i,\,i\in\nn\}$ is an
excellent metric on $\cs$.
Let $F:\cs\la\ct$ be a good extension with respect to
the metric.
Then the diagram below commutes
\[\xymatrix@C-10pt{
 &\fl'(\cs)\op\ar[ld]_{\cy\op}\ar[rd]
\ar@/^1pc/[rrrd]^{\cy'} & & \\
\fl(\cs)\op\ar[rr]_-{\wh Y} & &
\fl\big(\fs(\cs)\op\big)
\ar@{^{(}->}[rr] &&
\MMod{\fs(\cs)\op}
}\]
with $\wh Y$ the equivalence
$\wh Y:\fl(\cs)\op\la\fl\big(\fs(\cs)\op\big)$ of
Proposition~\ref{P3.94949}.
\elem

\prf
If we delete the unlabeled slanted
arrow in the diagram,
what is left is the perimeter of the
diagram of Observation~\ref{O94.5323}, which
commutes. And this slanted arrow exists
just because
the commutativity forces
$\cy':\fl'(\cs)\op\la\MMod{\fs(\cs)\op}$
to take $\fl'(\cs)\op$ into the essential image of
$\fl(\cs)\op$ under the map to $\MMod{\fs(\cs)\op}$.
\eprf

In future we will refer to the
unlabeled slanted arrow, in the
diagram of Lemma~\ref{L94.5325},
as $\cy':\fl'(\cs)\op\la\fl\big(\fs(\cs)\op\big)$.

\cor{C94.5327}
Let the hypotheses be as in Lemma~\ref{L94.5325}.
Consider the commutative triangle below
\[\xymatrix@C-10pt{
 &\fl'(\cs)\op\ar[ld]_{\cy\op}\ar[rd]^{\cy'} \\
\fl(\cs)\op\ar[rr]_-{\wh Y} & &
\fl\big(\fs(\cs)\op\big)
}\]
obtained from the diagram of 
Lemma~\ref{L94.5325} by deleting the $\MMod{\fs(\cs)\op}$
from the picture.

Then the functors $\cy$ and $\cy'$ are both essentially
surjective. Also: both are full, and the kernels
of the natural maps
\[
\begin{array}{ccc}
\Hom_{\fl'(\cs)}^{}(A,B)\la \Hom_{\fl(\cs)}^{}\big(\cy(A),\cy(B)\big)&
\qquad& \text{ and }\\
\Hom_{\fl'(\cs)}^{}(A,B)\la \Hom_{\fl\big(\fs(\cs)\op\big)\op}^{}\big(\cy'(A),\cy'(B)\big)& & 
\end{array}
\]
are the same.
\ecor

\prf
The commutativity of the triangle,
coupled with the fact that $\wh Y$ is an
equivalence, immediately forces the kernels
to be the same.
Also: since the functor $\wh Y$ is
an equivalence and the triangle commutes,
the other two statements can be checked on
the functor $\cy$; after all $\cy'$ is isomorphic
to $\cy\op$.

And now: the essential surjectivity of $\cy$
comes from \cite[Observation~\ref{O21.9}]{Neeman18A},
and the fullness
of the functor, meaning the surjectivity
of the map
\[\xymatrix@C+30pt{
\Hom_{\fl'(\cs)}^{}(A,B)\ar[r] & \Hom_{\fl(\cs)}^{}\big(\cy(A),\cy(B)\big)\ ,
}\]
is immediate from
\cite[Lemma~\ref{L21.11}]{Neeman18A}.
\eprf

\rmk{R94.101017}
The assertion in Corollary~\ref{C94.5327}, that
the kernels are the same, can be reformulated
more concretely as follows.
For a morphism $f:A\la B$ in the category
$\fl'(\cs)$, the following are equivalent:
\be
\item
For all objects $s\in\cs$, any composite
$F(s)\la A\stackrel f\la B$ vanishes.
\item
For all objects $t\in\fs(\cs)$, any composite
$A\stackrel f\la B\la t$ vanishes.
\ee
\ermk

\dis{D94.5}
It's high time
to come to Problem~\ref{P94.1}(ii). Recall: excellent
metrics were shown to be very special. The metric
$\{\cm_i,\,i\in\nn\}$ on $\cs$ induces a metric
$\{\cn\op_i,\,i\in\nn\}$
on $\fs(\cs)\op$---so far we are not appealing to
excellence, see Lemma~\ref{L3.1}.
But from Proposition~\ref{P97.1} we
learn that, if the metric $\{\cm_i,\,i\in\nn\}$
is excellent, then so is $\{\cn\op_i,\,i\in\nn\}$.
Moreover: by Proposition~\ref{P3.94949}
we have a natural equivalence
$\wh Y:\fl(\cs)\op\la\fl\big(\fs(\cs)\op\big)$,
which restricts to equivalences of the subcatgories
$\cl\op_i\subset\fl(\cs)\op$ with
$\wh\cl_i\subset\fl\big(\fs(\cs)\op\big)$.

And now the Problem~\ref{P94.1}(ii) can be made a
little more precise: given a good
extension $F:\cs\la\ct$, how much of this picture
does $\ct$ see? Is the inclusion $G:\fs(\cs)\op\la\ct\op$
a good extension? Is there a version of
Proposition~\ref{P3.94949} for the subcategories
$\fl'(\cs)\subset\ct$ and
$\fl'\big(\fs(\cs)\op\big)\subset\ct\op$?
\edis

\dsc{D94.5321}
Formulated in the generality of
Discussion~\ref{D94.5}, I do not know the answer to
Problem~\ref{P94.1}(ii). However: from
Corollary~\ref{C93.95} we learn that, if $F:\cs\la\ct$
is a good extension satisfying the hypotheses
of the Corollary, then the embedding
$\cs\la\fl'(\cs)$ is already a good extension.
And because the distance between $\cs$ and $\fl'(\cs)$
is substantially shorter, we will be able to
prove results about the embedding
$\fs(\cs)\op\la\fl'(\cs)\op$.
\edsc

We prepare the ground with the easy

\lem{L94.6}
Let $\cs$ be a triangulated category, and assume
$\{\cm_i,\,i\in\nn\}$ is an
excellent metric on $\cs$.
Let $F:\cs\la\ct$ be a good extension with respect to
the metric.
Then we have the containment
\[
F(\cs)\sub\bigcup_{i\in\nn}{^\perp\cl'_i}\ .
\]
\elem

\prf
In Remark~\ref{R3.?!?2} we saw the inclusion
$
Y(\cs)\subset\bigcup_{i\in\nn}{^\perp\cl_i}
$
which, by Lemma~\ref{L93.-1}, we can rewrite as
$
Y(\cs)\subset\bigcup_{i\in\nn}{^\perp\cy(\cl'_i)}
$.
Reformulating this yet again, we obtain that
$\Hom\big(Y(E),\cy(\cl'_i)\big)=0$ for any $E\in\cs$
and for $i=i(E)\gg0$. But then
\cite[Observation~\ref{O21.-1}]{Neeman18A}
tells us that
$\Hom\big(F(E),\cl'_i\big)=0$ for any $E\in\cs$
and for $i=i(E)\gg0$, which is a reformulation of
the inclusion in the Lemma.
\eprf

\lem{L94.657}
Let $\cs$ be a triangulated category, and assume
$\{\cm_i,\,i\in\nn\}$ is an
excellent metric on $\cs$.
Let $F:\cs\la\ct$ be a good extension with respect to
the metric. Let $\fs(\cs)$ have its induced metric
$\{\cn_i,\,i\in\nn\}$ of Lemma~\ref{L3.1}.
And, by the abuse of notation we adopted
in Remark~\ref{R93.-11}, we will view $\fs(\cs)$ as
a subcategory of $\fl'(\cs)\subset\ct$.

Suppose $a_*$ is a Cauchy sequence in $\cs$
with respect to the metric $\{\cm_i,\,i\in\nn\}$,
and let
$b_*$ be a Cauchy sequence in $\fs(\cs)\op$
with respect to the metric $\{\cn\op_i,\,i\in\nn\}$.
Then,
for any object $A\in\cs$ and for any object
$B\in\fs(\cs)$, we have the following:
\be
\item
The inverse sequence $\Hom\big(F(A),b_*\big)$ eventually stabilizes,
hence it follows that
$\climone\Hom\big(F(A),b_*\big)=0$.
\item
The inverse sequence $\Hom\big(F(a_*),B\big)$
eventually stabilizes,
hence it follows that
$\climone\Hom\big(F(a_*),B\big)$.
\ee
\elem

\prf
To prove (i), we first
appeal to Lemma~\ref{L94.6}
to choose an integer $n>0$ with $F(A)\in{^\perp\cl'_n}$.
And then, since the sequence $b_*$ is Cauchy, we can
choose $N$ large enough so that, for all $i\geq N$,
in the triangles $d_i\la b_{i+1}^{}\la b_i\la\T d_i$
we have that $d_i$ and $\T d_i$ belong to
$\cn_n$, and by
Lemma~\ref{L94.-9}
we have $\cn_n\subset\cl'_n$.
It follows that, for all $i\geq N$,
the map $\Hom\big(F(A),-\big)$ takes
the morphism $b_{i+1}^{}\la b_i$ to an isomorphism.

Now for the proof of (ii). We know that the object $B$
belongs to $\fs(\cs)\subset\fl'(\cs)$, and
Lemma~\ref{L93.-5} tells us that it must belong
to $(\cl'_n)^\perp$ for some $n>0$. Now choose
the integer $N>0$ so large that, for all $i\geq N$,
the triangle $d_i\la a_i\la a_{i+1}^{}\la\T d_i$ has
$d_i$ and $\T d_i$ both in $\cm_n$. Thus
$F(d_i)$ and $F(\T d_i)$ both belong
to $F(\cm_n)\subset\cl'_n$, and we conclude that
the functor $\Hom(-,B)$ must take
$F(a_i)\la F(a_{i+1}^{})$
to an isomorphism.
\eprf

\lem{L94.7}
Let $\cs$ be a triangulated category, and assume
$\{\cm_i,\,i\in\nn\}$ is an
excellent metric on $\cs$.
Let $F:\cs\la\ct$ be a good extension with respect to
the metric. Let $\fs(\cs)$ have its induced metric
$\{\cn_i,\,i\in\nn\}$ of Lemma~\ref{L3.1}.
Assume further that:
\be
\item
The category $\fs(\cs)$ has  
countable
products of the form $\prod_{n\in\nn}E$, for
any object $E\in\fs(\cs)$.
\item
For every integer $i>0$ there exists an integer 
$r=r(i)>0$ such that,
if $E$ belongs to $\cn_{i+r}$,
then $\prod_{n\in\nn}E$ belongs to $\cn_i$.
\setcounter{enumiv}{\value{enumi}}
\ee
Then for any Cauchy sequence $E_*$ in
the category $\fs(\cs)\op$,
with respect to the metric $\{\cn\op_i,\,i\in\nn\}$,
the product $\prod_{n\in\nn}E_i$
exists in the category $\fl'(\cs)$.
Moreover: the functor
$\cy:\fl'(\cs)\la\fl(\cs)$
respects these products.
\elem

\prf
Recall the commutative triangle of
Corollary~\ref{C94.5327}.
First of all: Corollary~\ref{C93.95}, applied
to the category $\fs(\cs)\op$ embedded
in $\fl\big(\fs(\cs)\op\big)$, tells us
that the coproduct of the objects
$\cy'(E_i)$ may be formed in
$\fl\big(\fs(\cs)\op\big)$. Therefore the
product of the objects $\cy(E_i)$ exists in
the category $\fl(\cs)\cong\fl\big(\fs(\cs)\op\big)\op$.
But
by
\cite[Observation~\ref{O21.9}]{Neeman18A}
the map
$\cy:\fl'(\cs)\la\fl(\cs)$ is essentially surjective,
and hence we may choose an object $E\in\fl'(\cs)$
such that $\cy(E)$ is the product
in $\fl(\cs)$ of the objects $\cy(E_i)$.

Now the morphisms $\cy(E)\la\cy(E_i)$ in the
category $\fl(\cs)$, which describe
$\cy(E)$ as the product of the objects $\cy(E_i)$,
are elements of $\Hom\big(\cy(E),\cy(E_i)\big)$
with $E_i\in\cy^{-1}\big(\fc(\cs)\big)$.
By
\cite[Lemma~\ref{L21.13}]{Neeman18A}
each of these morphisms lifts,
uniquely, to a morphism $\pi_i:E\la E_i$ in
the category $\fl'(\cs)$. What we need to prove
is that $E$, together with the projections
$\pi_i:E\la E_i$, satisfies in
$\fl'(\cs)$ the universal property of a product.

Let $a\in\cs$ be any object, and
observe the following commutative square
\[\xymatrix@C+40pt{
\Hom_{\ct}^{}\big(F(a),E\big)\ar[r]\ar[d]_\wr &
\ds\prod_{n=1}^\infty\Hom_{\ct}^{}\big(F(a),E_n\big)
\ar[d]^\wr \\
\Hom_{\MMod\cs}^{}\big(Y(a),\cy(E)\big)\ar[r]^-\sim &
\ds\prod_{n=1}^\infty\Hom_{\MMod\cs}^{}\big(Y(a),\cy(E_n)\big)
}\]
The vertical maps are isomorphisms by
\cite[Observation~\ref{O21.-1}]{Neeman18A},
and the bottom horizontal
map is an isomorphism because, in the category
$\fl(\cs)$, the object
$\cy(E)$ is the product of the $\cy(E_i)$.
From the commutativity of the
square we learn that the map in the top row
is an isomorphism.

Now let $a_*$ be a Cauchy sequence in $\cs$.
The functor $\Hom(-,E)$ takes the sequence
$F(a_*)$ to the inverse sequence
\[
\prod_{i=1}^\infty\Hom\big(F(a_*),E_i\big)
\]
which, by Lemma~\ref{L94.657}(ii), is a
product of sequences each of which stabilizes.
Since $\climone$ respect products
of abelian groups. We have that
\be
\setcounter{enumi}{\value{enumiv}}
\item
For every Cauchy sequence $a_*$ in $\cs$,
we have
$\climone\Hom\big(F(a_*),E\big)=0$.
\setcounter{enumiv}{\value{enumi}}
\ee
Even more easily, since $E_i\in\fs(\cs)$ and
$a_*$ is a Cauchy sequence in $\cs$,
Lemma~\ref{L94.657}(ii) gives that
\be
\setcounter{enumi}{\value{enumiv}}
\item
For every $i\in\nn$ and
every Cauchy sequence $a_*$ in $\cs$,
we have
$\climone\Hom\big(F(a_*),E_i\big)=0$.
\setcounter{enumiv}{\value{enumi}}
\ee
Let $A$ be an object of $\fl'(\cs)$,
and choose a Cauchy sequence
$a_*$ in the category $\cs$ with
$A=\hoco F(a_*)$.
And now \cite[Lemma~\ref{L21.11}]{Neeman18A} tells us 
that the vertical maps in the commutative square
below are isomorphisms
\[\xymatrix@C+40pt{
\Hom_{\ct}^{}\big(A,E\big)\ar[r]\ar[d]_\wr &
\ds\prod_{n=1}^\infty\Hom_{\ct}^{}\big(A,E_n\big)
\ar[d]^\wr \\
\Hom_{\MMod\cs}^{}\big(\cy(A),\cy(E)\big)\ar[r]^-\sim &
\ds\prod_{n=1}^\infty\Hom_{\MMod\cs}^{}\big(\cy(A),\cy(E_n)\big)
}\]
The bottom horizontal map is an isomorphism
because, in the category $\fl(\cs)$, the object $\cy(E)$
is the product of the objects $\cy(E_i)$. Hence
the horizontal map in the top row is an
isomorphism, and $E$ is the product of the $E_i$
in the category $\fl'(\cs)$.
\eprf

Since the next statement is a Proposition we
make the statement self-contained.

\pro{P94.9}
Let the hypotheses be as in Lemma~\ref{L94.7},
meaning that $\cs$ is a triangulated category, we assume
$\{\cm_i,\,i\in\nn\}$ is an
excellent metric on $\cs$,
and we let $F:\cs\la\ct$ be a good extension with respect to
the metric. Let $\fs(\cs)$ have its induced metric
$\{\cn_i,\,i\in\nn\}$

Assume in addition that:
\be
\item
The subcategories $\cn_i\subset\fs(\cs)$
are closed in $\fs(\cs)$ under direct
summands.
\item
The category $\fs(\cs)$ has countable
products of the form $\prod_{n\in\nn}E$, for
any object $E\in\fs(\cs)$.
\item
For every integer $i>0$ there exists an integer 
$r=r(i)>0$ such that,
if $E$ belongs to $\cn_{i+r}$,
then $\prod_{n\in\nn}E$ belongs to $\cn_i$.
\setcounter{enumiv}{\value{enumi}}
\ee
This completes the long list of hypotheses.

Then the inclusion $\fs(\cs)\op\la\fl'(\cs)\op$ is
a good extension with respect to the metric
$\{\cn\op_i,\,i\in\nn\}$.

Moreover: every object in $\fl'(\cs)\op$
is isomorphic to the homotopy colimit of
a Cauchy sequence in
$\fs(\cs)\op$ with respect to the metric
$\{\cn\op_i,\,i\in\nn\}$. In other words: for our particular
good extension $\fs(\cs)\op\la\fl'(\cs)\op$,
we have that the inclusion 
\[
\fl'\big(\fs(\cs)\op\big)\sub
\fl'(\cs)\op\ 
\]
is an equality.
\epro

\prf
First of all: by Lemma~\ref{L94.7}, if $b_*$ is
a Cauchy sequence in the category $\fs(\cs)\op$,
then the coproduct of the $b_i$ exists in $\fl'(\cs)\op$.
In other words: in
the category $\fl'(\cs)$ we can form $\prod_{i\in\nn}b_i$.
And by Proposition~\ref{P93.9} the subcategory
$\fl'(\cs)\subset\ct$ is triangulated,
meaning we can form in $\fl'(\cs)\op$ the homotopy
colimits of Cauchy sequences in $\fs(\cs)\op$;
in $\fl'(\cs)$ these are the homotopy limits of
the Cauchy inverse sequences in $\fs(\cs)$.
In any case:
\cite[Definition~\ref{D21.1}(i)]{Neeman18A}
is satisfied
for the inclusion $\fs(\cs)\op\la\fl'(\cs)\op$.
To show that the inclusion is a good extension with
respect to the metric it suffices to check
\cite[Definition~\ref{D21.1}(ii)]{Neeman18A}.
In the notation
of the commutative diagram of Lemma~\ref{L94.5325},
it remains to prove that the composite
$\cy':\fl'(\cs)\op\la\fl\big(\fs(\cs)\op\big)\hookrightarrow\MMod{\fs(\cs)\op}$
takes homotopy colimits of Cauchy sequences
to colimits.

Let $b_*$ be a Cauchy sequence on $\fs(\cs)\op$, and
form in $\fl'(\cs)$ the triangle
\[\xymatrix{
\holim b_*\ar[r]^-\ph &
\ds\prod_{i=1}^\infty b_i\ar[rr]^-{\id-\sh} &&
\ds\prod_{i=1}^\infty b_i\ar[r] &
\T\holim b_*
}\]
The next step is to apply the functor
$\cy:\cl'(\cs)\la\cl(\cs)$, which amounts to studying
what happens when you apply 
$\Hom\big(F(a),-\big)$ to the triangle, with $a\in\cs$.
But here Lemma~\ref{L94.657}(i) comes to our aid:
it tells us that, for $a\in\cs$ and $b_*$
a Cauchy sequence in $\fs(\cs)\op$, we have 
$\climone \Hom\big(F(a),b_*\big)=0$. Hence the map
\[\xymatrix{
\ds\Hom\left(F(a)\,,\,\prod_{i=1}^\infty b_i\right)\ar[rr]^-{\id-\sh} &&
\ds\Hom\left(F(a)\,,\,\prod_{i=1}^\infty b_i\right)
}\]
is an epimorphism. The long exact sequence, obtained
by applying the functor
$\Hom\big(F(a),-\big)$ to the triangle,
breaks up into short exact sequences. We deduce that,
when we  apply the functor $\Hom\big(F(a),-\big)$
to the map $\ph$ in the triangle, we obtain
an isomorphism
\[\xymatrix{
\Hom\big(F(a)\,,\,\holim b_*\big)
\ar[r] &
\clim\Hom\big(F(a),b_*\big)\ .
}\]
As this is true for every object
$a\in\cs$, we have proved that the functor
$\cy:\fl'(\cs)\la\MMod\cs$ satisfies
\[
\cy\big(\holim F(b_*)\big)\eq \clim b_*\ .
\]
Of course: because the essential image of $\cy$
is contained in $\fl(\cs)\subset\MMod\cs$,
the object $\cy\big(\holim F(b_*)\big)\cong\clim b_*$
must belong to
$\fl(\cs)$. As it satisfies in $\MMod\cs$ the
universal property of a limit, it must satisfy it
in the full subcategory $\fl(\cs)$. Passing to opposite
categories gives that the functor
$\cy\op:\fl'(\cs)\op\la\fl(\cs)\op$ takes homotopy
colimits,
of Cauchy sequences $b_*$
in $\fs(\cs)\op\subset\fl'(\cs)\op$, to colimits
in $\fl(\cs)\op$.

But now the commutative triangle of
Corollary~\ref{C94.5327} comes to our aid: it follows
immediately that the functor
$\cy':\fl'(\cs)\op\la\fl\big(\fs(\cs)\op\big)$
takes homotopy
colimits,
of Cauchy sequences $b_*$
in $\fs(\cs)\op\subset\fl'(\cs)\op$, to colimits
in $\fl\big(\fs(\cs)\op\big)$.
But by
\cite[Definition~\ref{D20.11}(i)]{Neeman18A},
the full subcategory
$\fl\big(\fs(\cs)\op\big)\subset\MMod{\fs(\cs)\op}$
has for objects all the colimits, in
$\MMod{\fs(\cs)\op}$, of Cauchy sequences $b_*$ contained
in $\fs(\cs)\op\subset\MMod{\fs(\cs)\op}$.
Therefore the colimits of such
Cauchy sequences exist in the subcategory,
and agree with the ones formed in 
$\MMod{\fs(\cs)\op}$.

This concludes the proof that the embedding
$\fs(\cs)\op\la\fl'(\cs)\op$ is a good extension with
respect to the metric $\{\cn\op_i,\,i\in\nn\}$.
It remains to show that every object in $\fl'(\cs)\op$
is isomorphic to the homotopy colimit
of a Cauchy sequence in $\fs(\cs)\op$.

Take therefore any object $X\in\fl'(\cs)$.
By
\cite[Lemma~\ref{L3.7}(iii)]{Neeman18A},
applied to the object
$\cy(X)\in\fl(\cs)$, there exists in $\fs(\cs)\op$
a Cauchy sequence $g_*$ with $\cy(X)=\clim g_*$.
But by the above we have that
we can form $\holim g_*$ in the
triangulated category $\fl'(\cs)$, and that
$\cy\big(\holim g_*\big)=\clim g_*$.
This gives an isomorphism
$\cy(X)\cong\cy\big(\holim g_*\big)$,
in the category $\fl(\cs)$, which,
by
\cite[Corollary~\ref{C21.908}]{Neeman18A},
can be
lifted to an isomorphism
$X\la\holim g_*$ in $\fl'(\cs)$.
\eprf

\rmk{R94.10}
Next we apply Definition~\ref{D93.11} to the
good extension $\fs(\cs)\op\la\fl'(\cs)\op$
of Proposition~\ref{P94.9}. 
We form subcategories $\wh\cl'_\ell\subset\fl'(\cs)\op$,
where the definition is that $\wh\cl'_\ell$
has for objects all the homotopy colimits
of Cauchy sequences in $\cn\op_\ell$.
\ermk

\cor{C94.11}
Let the hypotheses be as in Proposition~\ref{P94.9},
and let the subcategories
$\wh\cl'_\ell\subset\fl'(\cs)\op$ be as
in Remark~\ref{R94.10}. Then we have equalities
of subcategories of $\fl'(\cs)\op$
\[
(\cl'_\ell)\op=\wh\cl'_\ell\ .
\]
\ecor

\prf
First recall that the last
assertion of Proposition~\ref{P94.9}
was the equality
\[
\fl'\big(\fs(\cs)\op\big)\eq
\fl'(\cs)\op\ .
\]
Now \cite[Corollary~\ref{C21.908}]{Neeman18A}
tells us that
both $\cy:\fl'(\cs)\la\fl(\cs)$ and
$\cy':\fl'\big(\fs(\cs)\op\big)\la\fl\big(\fs(\cs)\op\big)$ induce
bijections between isomorphism classes of objects,
and Lemma~\ref{L93.-1} says that $\cy(\cl'_\ell)=\cl_\ell$
and $\cy'(\wh\cl'_\ell)=\wh\cl_\ell$. Thus
the current corollary follows from the commutative triangle
of Corollary~\ref{C94.5327}, coupled with
Proposition~\ref{P3.94949}, which tells us that the
equivalence of categories
$\wh Y:\fl(\cs)\op\la\fl\big(\fs(\cs)\op\big)$
restricts to equivalences $\wh Y:\cl\op_\ell\la\wh\cl_\ell$.
\eprf

\section{An abundant source of examples of very good metrics}
\label{S95}

This section will show that very good metrics are
plentiful. We will only seriously study one class
of examples, ones that come from bounded-above \tstr{s}.
Triangulated categories $\cs$ possessing bounded-above
\tstr{s} all have natural very good metrics.

Now suppose we are given a triangulated
category $\ct$ a with a nondegenerate \tstr,
and assume $\ct$ has countable coproducts.
Then $\ct^-$ has the very good metric
of the last paragraph.
We will show, under mild hypotheses, that the embedding
$\ct^-\subset\ct$ is a good extension with
respect to the metric, and furthermore
the inclusion is an equality $\fl'(\ct^-)=\ct$.

And, at the end of the section, we will show that
in some sense this is a general phenomenon.
Suppose
$\cs$ is a triangulated category with a bounded-above
\tstr, the very good metric is the one coming from the
\tstr, and $F:\cs\la\ct$ is any good extension with
respect to the metric. Then
the subcategory $\fl'(\cs)\subset\ct$
has a non-degenrate \tstr\ extending the
one on $\cs$. And, under mild hypotheses,
we find ourselves in
the situation in the paragraph above: we have that
$\cs=\fl'(\cs)^-$ and the fully faithful
functor $\cs=\fl'(\cs)^-\la\fl'(\cs)$ is a good extension.

Let us first make a definition
and 
prove a couple of technical little lemmas,
which will be helpful.

\dfn{D95.-5}
Let $\ca$ be a category. A sequence $A_*$ in the category
$\ca$, meaning a diagram
$A_1\la A_2\la A_3\la$, is said to \emph{eventually
stabilize} if there
exists an integer $N>0$ such that,
for all integers $j>i\geq N$,
the map $A_i\la A_j$ is an
isomorphism.
\edfn

\rmk{R111119576}
In a random category $\ca$, colimits do not normally
exist. But sequences that eventually stabilize are
an exception: if the sequence
$A_*$ eventually stabilizes, just choose your
integer $N$ large enough so that,
for all integers $j>i\geq N$,
the maps $A_i\la A_j$ are
isomorphisms. Then not only does $\colim\,A_*$ exist,
we further know that the map $A_N\la\colim\,A_*$ is an
isomorphism.
\ermk

\rmk{R95.-7}
If $\ct$ is a triangulated category with a \tstr\
$\tst\ct$, if $A_*$ is any sequence in $\ct$,
and if $n\in\zz$ is an integer, by the
symbol
$A_*^{\leq n}$ we mean
the sequence $A_1^{\leq n}\la A_2^{\leq n}\la A_3^{\leq n}\la \cdots$. We will care about sequences $A_*$ for which
the sequences $A_*^{\leq n}$
eventually stabilize.
\ermk

\lem{L95.11}
Assume $\ct$ is a triangulated category with
a \tstr\ $\tst\ct$. Assume
also that countable coproducts
exist in $\ct$.

Suppose further that there exists an 
an integer $r>0$ such that,
if $\{X_i,\,i\in\nn\}$ 
is a countable set of objects in $\ct^{\geq0}$,
then
$\coprod_{i=1}^\infty X_i$ belongs to $\ct^{\geq-r}$.

If $A_*$ is a sequence in $\ct$ such that,
for any integer $n\in\zz$, 
the sequence $A_*^{\leq n}$ eventually stablizes,
then the map
\[\xymatrix@C+30pt{
\colim\,A_*^{\leq n}\ar[r]^{\ph_n} &\Big(\hoco A_*\Big)^{\leq n}
}\]
is an isomorphism.
\elem

\prf
Pick a sequence $A_*$ in $\ct$ satisfying the hypothesis.
Choose any integer $n$, and we want to prove that
the map $\ph_n$ is an isomorphism.

With $r>0$ as in the
second paragraph of the
statement of the Lemma, consider the sequence of triangles
$A_*^{\leq n+r+1}\la A_*\la A_*^{\geq n+r+2}\la\T A_*^{\leq n+r+1}$.
The coproducts in 
the following diagram are countable
and hence exist in $\ct$,
the rows are
triangles, and the square commutes
\[\xymatrix@C+30pt{
\ds\coprod_{i=1}^\infty A_i^{\leq n+r+1} \ar[r]\ar[d]_{\id-\sh} &
\ds\coprod_{i=1}^\infty A_i \ar[r]\ar[d]^{\id-\sh} &
\ds\coprod_{i=1}^\infty A_i^{\geq n+r+2} \\
\ds\coprod_{i=1}^\infty A_i^{\leq n+r+1} \ar[r] &
\ds\coprod_{i=1}^\infty A_i \ar[r] &
\ds\coprod_{i=1}^\infty A_i^{\geq n+r+2} \\
}\]
We may complete the square in $\ct$
to a $3\times3$ diagram whose
rows and columns are triangles
\[\xymatrix@C+30pt{
\ds\coprod_{i=1}^\infty A_i^{\leq n+r+1} \ar[r]\ar[d]_{\id-\sh} &
\ds\coprod_{i=1}^\infty A_i \ar[r]\ar[d]^{\id-\sh} &
\ds\coprod_{i=1}^\infty A_i^{\geq n+r+2}\ar[d] \\
\ds\coprod_{i=1}^\infty A_i^{\leq n+r+1} \ar[r]\ar[d] &
\ds\coprod_{i=1}^\infty A_i \ar[r]\ar[d] &
\ds\coprod_{i=1}^\infty A_i^{\geq n+r+2}\ar[d] \\
\colim\,A_*^{\leq n+r+1} \ar[r] &
\hoco A_* \ar[r] & C
}\]
In the first column, because the sequence eventually
stabilizes, the homotopy colimit is a genuine colimit
and is equal to $A_i^{\leq n+r+1}$ for $i\gg0$. And the third
column gives a triangle
\[\xymatrix@C+30pt{
\ds\coprod_{i=1}^\infty A_i^{\geq n+r+2}\ar[r] &
C\ar[r] &\ds
\T\left(\coprod_{i=1}^\infty A_i^{\geq n+r+2}\right)
}\]
which, by the properties
of $r>0$ in second paragraph of the
statement of the Lemma,
exhibits $C$ as an objects in
$\ct^{\geq n+2}*\ct^{\geq n+1}\subset\ct^{\geq n+1}$.

Choose $i\gg0$ large enough so that
$A_i^{\leq n+r+1}\la\colim\,A_*^{\leq n+r+1}$ is an
isomorphism, and complete the composable morphisms
$A_i^{\leq n}\la A_i^{\leq n+r+1}\la \hoco A_*$ to
an octahedron
\[\xymatrix@C+30pt{
A_i^{\leq n}\ar@{=}[d]\ar[r] &
  A_i^{\leq n+r+1}\ar[r]\ar[d] &
  \Big(A_i^{\leq n+r+1}\Big)^{\geq n+1}\ar[d] \\ 
A_i^{\leq n}\ar[r] &  \hoco A_*\ar[r]\ar[d] &
      B\ar[d] \\
 & C\ar@{=}[r] & C
}\]
The third column exhibits $B$ as an object of
$\ct^{\geq n+1}*\ct^{\geq n+1}\subset\ct^{\geq n+1}$,
and hence the triangle in the second row must be
canonically isomorphic to
\[\xymatrix@C+30pt{
\Big(\hoco A_*\Big)^{\leq n} \ar[r] &      
\hoco A_* \ar[r] &
\Big(\hoco A_*\Big)^{\geq n+1}\ .
}\]
Thus, for our choice of $i\gg0$, the map
$A_i^{\leq n}\la \big(\hoco A_*\big)^{\leq n}$ is
indeed an isomorphism.
\eprf

\ntn{N95.-3}
Let 
$\ct$ be a triangulated category with
a \tstr\ $\tst\ct$, and let $Y\in\ct$ be an object.
The sequence
$Y^{\leq1}\la Y^{\leq2}\la Y^{\leq3}\la \cdots$,
will be denoted  $Y^{\leq(*)}$.
\entn

For any object $Y\in\ct$, letting $A_*$ be
the sequence $Y^{\leq(*)}$ of Notation~\ref{N95.-3}
gives a sequence satisfying the hypothesis
on $A_*$
in Lemma~\ref{L95.11}. We apply this in

\lem{L95.-1}
As in Lemma~\ref{L95.11},
assume $\ct$ is a triangulated category with
a \tstr\ $\tst\ct$,
and that countable coproducts
exist in $\ct$. Still as in Lemma~\ref{L95.11},
we further suppose that there exists an 
an integer $r>0$ such that,
if $\{X_i,\,i\in\nn\}$ 
is a countable set of objects in $\ct^{\geq0}$,
then
$\coprod_{i=1}^\infty X_i$ belongs to $\ct^{\geq-r}$.

And now we add the assumption that the \tstr\
$\tst\ct$ is nondegenerate.

Suppose $Y\in\ct$ is
an object, let $Y^{\leq(*)}$ be as in
Notation~\ref{N95.-3}, and let
$X=\hoco Y^{\leq(*)}$. Then
there is a (non-canonical) isomorphism $\psi:X\la Y$.
\elem

\prf
We have the following diagram
\[\xymatrix@C+5pt{
\ds\coprod_{i=1}^\infty Y^{\leq i} \ar[rrr]^-{\id-\sh} & & & 
\ds\coprod_{i=1}^\infty Y^{\leq i} \ar[rr]^-\gamma \ar[drr]_\rho & &
\hoco Y^{\leq(*)}\ar@{=}[r] &X \\
& & & && Y
}\]
where the top row is a triangle, the map
$\rho$ is the unique morphism from the coproduct
which restricts to
the natural map $Y^{\leq i}\la Y$ for every $i\in\nn$,
and the composite $\rho\circ(\id-\sh)$ vanishes.
Therefore there exists a factorization of $\rho$
as $\psi\circ\gamma$. This produces for us the desired map
$\psi:X\la Y$.
We need to prove it an isomorphism.

By construction, for every integer $n>0$ the
triangle below commutes
\[\xymatrix@R-20pt@C+10pt{
                & & &X\ar@{=}[r]\ar[dd]^\psi &\hoco Y^{\leq(*)}\\
Y^{\leq n} \ar[rrru]\ar[rrrd] & & & & \\
 & & & Y & \\
}\]
Now apply the functor $(-)^{\leq n}$,
to obtain the commutative triangle
\[\xymatrix@R-20pt@C+10pt{
                & & &X^{\leq n}\ar[dd]^{\psi^{\leq n}}\\
Y^{\leq n} \ar[rrru]^-{\ph_n}\ar[rrrd]_-\id & & & \\
 & & & Y^{\leq n}  \\
}\]
where $\ph_n$ is the isomorphism of Lemma~\ref{L95.11}.
The commutativity forces
$\psi^{\leq n}:X^{\leq n}\la Y^{\leq n}$ to be an isomorphism.

Now complete $\psi:X\la Y$ to a triangle
$X\stackrel\psi\la Y\la D\la\T X$.
By the above we deduce that $D^{\leq n}=0$ for all $n\in\zz$,
and hence $D\in\cap_{n\in\zz}\ct^{\geq n}$.
But the \tstr\ is assumed nondegenerate, and
we deduce that $D=0$. Hence $\psi:X\la Y$ nust be
an isomorphism.
\eprf

With the basic lemmas out of the way, we come to
the example we care about most---this is the one
for which we have immediate applications
in a forthcoming paper.

\exm{E95.1}
Let $\cs$ be a triangulated category with a
bounded-above
\tstr\ $\tst\cs$. Recall: being bounded-above means
that
\[
\cs\eq\bigcup_{i=0}^\infty \cs^{\leq i}\ .
\]
Define the metric $\{\cm_i,\,i\in\nn\}$
by the rule $\cm_i=\cs^{\geq i}$.

This metric is always very good.

After all: $^\perp\cm_i={^\perp\cs^{\geq i}}=\cs^{\leq i-1}$,
and the hypothesis $\cs=\cup_{i\in\nn}{^\perp\cm_i}$
of Definition~\ref{D93.1}(i)
comes from the \tstr\ being bounded-above,
as exhibited in
the displayed formula of the last paragraph.
And the hypothesis of Definition~\ref{D93.1}(ii)
is also easy: given any object $B\in\cs$,
there exists a triangle $B^{\leq i-1}\la B\la B^{\geq i}$,
with $B^{\leq i-1}\in\cs^{\leq i-1}={^\perp\cm_i}$ and with
$B^{\geq i}\in\cs^{\geq i}=\cm_i$.
\eexm

\lem{L95.2}
Let $\cs$ be a triangulated category with a bounded-above
\tstr\ $\tst\cs$, and let the metric $\{\cm_i,\,i\in\nn\}$
be as in Example~\ref{E95.1}.

A sequence $E_1\la E_2\la E_3\la\cdots$ in the category
$\cs$ is Cauchy if and only if, for every
unteger $n>0$, the sequence
$E_*^{\leq n}$, of Remark~\ref{R95.-7}, eventually
stabilizes as in Definition~\ref{D95.-5}.
\elem

\prf
Suppose the sequence is Cauchy and we are given an integer
$n>0$. By the Cauchyness, there must exist an integer
$N>0$ such that, if $j>i>N$ are integers, then in
the triangle $E_i\la E_j\la D_{i,j}$ we have
$D_{i,j}\in\cm_{n+1}=\cs^{\geq n+1}$. But then the
map $E_i^{\leq n}\la E_j^{\leq n}$ must be an isomorphism.

Conversely: if we are given a sequence
$E_1\la E_2\la E_3\la\cdots$, and we know that for every
integer $n>0$ the
sequence
$E_1^{\leq n}\la E_2^{\leq n}\la E_3^{\leq n}\la\cdots$ eventually
stabilizes, we need to show the Cauchyness.
Choose therefore an integer $n>0$, and then an integer
$N>0$ such that, for all $j>i>N$, the map
$E_i^{\leq n}\la E_j^{\leq n}$
is an isomorphism. Completing
$E_i^{\leq n}\la E_i\la E_j$ to an octahedron
\[\xymatrix@C+20pt{
E_i^{\leq n}\ar[r]\ar@{=}[d] & 
  E_i\ar[d]\ar[r] & E_i^{\geq n+1}\ar[d] \\
E_i^{\leq n}\ar[r] & E_j\ar[r]\ar[d] & E_j^{\geq n+1}\ar[d] \\
 & D_{i,j} \ar@{=}[r] & D_{i,j}
}\]
the triangle $E_j^{\geq n+1}\la D_{i,j}\la\T E_i^{\geq n+1}$
allows us to deduce that $D_{i,j}$ belongs
to $\cs^{\geq n}*\cs^{\geq n}\subset\cs^{\geq n}=\cm_n$.
\eprf

And now we come to what
(for us) will be
the main examples of good extensions
$F:\cs\la\ct$, for  which the category
$\fl'(\cs)$ of
\cite[Definition~\ref{D21.7}]{Neeman18A}
turns
out to satisfy $\ct=\fl'(\cs)$.

\exm{E95.13}
Let $\ct$ be a triangulated category with countable
coproducts, and let $\tst\ct$ be a
nondegenerate \tstr\ on $\ct$.
Assume further that there exists an
integer $r>0$ such that, for any
countable set of
objects $\{X_i,\,i\in\nn\}$ contained
in $\ct^{\geq0}$, the coproduct
$\coprod_{i\in\nn}X_i$ belongs to $\ct^{\geq-r}$.
Then
\be
\item
The category $\cs=\ct^-=\cup_{i\in\nn}\ct^{\leq i}$,
with the metric
$\cm_i=\ct^-\cap\ct^{\geq i}$, is a triangulated category
with a very good metric. And the embedding
$\cs=\ct^-\la\ct$ is a good extension with respect
to the metric, with
$\fl'(\cs)=\ct$.
\item
The category
$\cs=\ct^b$, with the metric
$\cm_i=\ct^b\cap\ct^{\geq i}$, is a triangulated category
with a very good metric. And the embedding
$\cs=\ct^b\la\ct^+$ is a good extension with respect
to the metric, with
$\fl'(\cs)=\ct^+$.  
\ee
\eexm

\prf
First
of all: both 
$\ct^b$ and $\ct^-$, with the \tstr\ restricted
from $\ct$, are triangulated
categories with a bounded-above \tstr{s}.
Example~\ref{E95.1} tells us that the metrics given,
respectively,
by $\cm_i=\ct^b\cap\ct^{\geq i}$
and by $\cm_i=\ct^-\cap\ct^{\geq i}$, are both very good.
And
by Lemma~\ref{L95.2} a sequence $A_*$,
in either $\ct^b$ or $\ct^-$, 
is Cauchy if and only if,
for every integer $n\in\zz$, the sequence
$A_*^{\leq n}$ eventually stabilizes.

Next: for each of the embeddings $\ct^b\la\ct^+$
and $\ct^-\la\ct$ we need to prove that
it is a good extension, meaning satisfies
\cite[Definition~\ref{D21.1}~(i) and (ii)]{Neeman18A}.
For the embedding $\ct^-\subset\ct$
we have that $\ct$ has all countable coproducts,
and hence
\cite[Definition~\ref{D21.1}(i)]{Neeman18A}
holds trivially.
For the embedding $\ct^b\la\ct^+$ we need to prove
something, we need to
show that, for any Cauchy sequence $A_*$ in $\ct^b$,
we can form in $\ct^+$ the coproduct $\coprod_{i\in\nn}A_i$.

The sequence $A_*$ is assumed Cauchy,
meaning that, 
for every integer $n\in\zz$, the sequence
$A_*^{\leq n}$ eventually stabilizes.
The special case where $n=0$ tells us that there must exist
an integer $N>0$ such that, for
$i\geq N$, the object $A_i^{\leq0}$
are all isomorphic. Because $A_N^{\leq0}$ belongs to
$\ct^b$, there is an integer $m>0$ with
$A_N^{\leq0}\in\ct^{\geq-m}$, and therefore
the isomorphic $A_i^{\leq0}$ will belong
to $\ct^{\geq-m}$ for all $i\geq N$.
It follows that $A_i\in\ct^{\geq-m}$ for all
$i\geq N$. And, since there are only
finitely many objects $A_i$ with $i<N$,
we may (after increasing $m$ if nessecary) choose
an integer $m>0$ such that $A_i\in\ct^{\geq-m}$ for all
$i\in\nn$. Therefore
the object $\coprod_{i\in\nn}A_i$, which we
may certainly form in the category $\ct$,
actually belongs
to $\ct^{\geq-m-r}\subset\ct^+$.

Next we prove that
\cite[Definition~\ref{D21.1}(ii)]{Neeman18A}
is satisfied, for both of
the inclusions $\ct^b\la\ct^+$ and
$\ct^-\la\ct$.
Recall: for a fully faithful functor
$F:\cs\la\ct$ to satisfy
\cite[Definition~\ref{D21.1}(ii)]{Neeman18A}
means that, for every object $X\in\cs$ and for every
Cauchy sequence $A_*$ in $\cs$, the natural map
\[\xymatrix@C+40pt{
\colim\,\Hom(X,A_*)\ar[r] &
\ds\Hom\left(F(X),\hoco F(A_*)\right)
}\]
is an isomorphism. By the previous paragraph
the embedding $\ct^+\la\ct$ respects the coproducts
of the form $\coprod_{i\in\nn}A_i$, where $A_*$
is a Cauchy sequence in $\ct^b$. Hence $\hoco A_*$,
formed in $\ct^+$ for a Cauchy sequence in $\ct^b$,
agrees with the $\hoco A_*$ formed in
$\ct$. Thus it suffices to prove
the assertion for $X\in\ct^-$ and $A_*$ a Cauchy
sequence in $\ct^-$; the assertion about $\ct^b$
is a special case.

But now: $X\in\ct^-$ implies that there exists an
integer $n>0$ with $X\in\ct^{\leq n}$.
And Lemma~\ref{L95.11} gives that the map
\[\xymatrix@C+30pt{
\colim\,A_*^{\leq n}\ar[r]^{\ph_n} &\Big(\hoco A_*\Big)^{\leq n}
}\]
is an isomorphism.  The fact that
\cite[Definition~\ref{D21.1}(ii)]{Neeman18A}
is satisfied follows
from the commutative square
\[\xymatrix@C+30pt{
\Hom\left(X\,,\,\colim\,A_*^{\leq n}\right)\ar[r]^-{\sim} 
\ar[d]_\wr &
\Hom\left(X\,,\,\Big(\hoco A_*\Big)^{\leq n}\right)
\ar[d]^\wr \\
\colim\,\Hom(X,A_*)\ar[r] &
\ds\Hom\left(X,\hoco A_*\right)
}\]
where the vertical isomorphisms come from $X$ belonging
to $\ct^{\leq n}$, and
the isomorphism in the top row comes
from applying $\Hom(X,-)$ to the isomorphism
$\ph_n$ of Lemma~\ref{L95.11}.

It remains to prove that, for the embedding
$\ct^b\la\ct^+$ we have $\fl'(\ct^b)=\ct^+$ and
for the embedding $\ct^-\la\ct$ we have
$\fl'(\ct^-)=\ct$. Choose therefore an
object $Y$, in either $\ct^+$ or $\ct$, and
we need to express it as the homotopy colimit
of a Cauchy sequence.

The Cauchy sequence we choose
is the sequence $Y^{\leq(*)}$
of Notation~\ref{N95.-3}. It is
definitely Cauchy, and if $Y\in\ct^+$ then
$Y^{\leq(*)}$ is contained in the category
$\ct^b$.
And Lemma~\ref{L95.-1} now tells us
that $Y\cong\hoco Y^{\leq(*)}$, and the right-hand-side
clearly belongs to
$\fl'(\cs)$. Thus we can express
every $Y$, in either $\ct^+$ or in $\ct$,
as the homotopy colimit
of some Cauchy sequence.

This completes the proof.
\eprf

\rmk{R95.-10}
In Example~\ref{E95.13} we started with a
triangulated category $\ct$, with a \tstr,
and observed
\be
\item
The category $\ct^-=\cup_{i\in\nn}\ct^{\leq i}$,
with the restricted \tstr,
is a triangulated category with
a bounded-above \tstr. Hence
it has the very good metric of 
Example~\ref{E95.1}.
\item
Under not-too-restrictive hypotheses,
the inclusion $\ct^-\la\ct$ is a
very good extension with respect
to the metric, and satisfies $\ct=\fl'(\cs)$.
\ee
We can wonder how general this phenomenon
is. In other words: suppose $\cs$ is a triangulated
category with a bounded-above \tstr, the very
good metric on $\cs$ is as in
Example~\ref{E95.1}, meaning $\cm_i=\ct^{\leq i}$,
and we are given some good extension $F:\cs\la\ct$.
Then does $\fl'(\cs)$ always have a \tstr\
extending the one on $\cs$?
Is the embedding $\cs\la\fl'(\cs)$ always
a good extension?
Is it always true that $\cs=\fl'(\cs)^-$?

Although some hypotheses need to be made
to guarantee this, they are remarkably mild.
The next few results will address this. 
\ermk

We start
with

\lem{L95.4}
Let $\cs$ be a triangulated category with a bounded-above
\tstr\ $\tst\cs$, and let the metric $\{\cm_i,\,i\in\nn\}$
be as in Example~\ref{E95.1}.  Assume $F:\cs\la\ct$
is a good extension with respect to the metric.
Then the category $\cl'_i\subset\fl'(\cs)$, of
Definition~\ref{D93.11}, is contained
in
$\fl'(\cs)\cap F\big(\cs^{\leq i-1}\big)^\perp$.
\elem

\prf
Suppose $C\in\cl'_i$ is an object.
Then there is a Cauchy sequence
$c_*$, contained in $\cm_i=\cs^{\geq i}$,
with $C=\hoco F(c_*)$.
By the conventions of Definition~\ref{D3.-3}(i)
\begin{itemize}
\item
It follows that $\cy(C)= \colim\,Y(c_*)$
must belong to $\cl_i\subset\fl(\cs)$.
\end{itemize}
Now we have the string of equalities
\[
\begin{array}{cclcl}
^\perp\cl_i\cap Y(\cs)
 &=&{^\perp Y(\cm_i)}\cap Y(\cs)
    &\quad&\text{By Lemma~\ref{L3.2.5}(i)}\\
  &=& Y\big({^\perp\cm_i}\big) & &
       \text{Because $Y$ is fully faithful and }\cm_i\subset\cs\\
&=& Y\big({^\perp\cs^{\geq i}}\big)
  & &\text{Because }\cm_i=\cs^{\geq i}\\
&=& Y\big({\cs^{\leq i-1}}\big)
  & &  \text{As }\cs^{\leq i-1}={^\perp\cs^{\geq i}}\ .  
\end{array}
\]
From this
we deduce the inclusion
$Y\big({\cs^{\leq i-1}}\big)\subset{^\perp\cl_i}$.
By $\bullet$ we know that $\cy(C)\in\cl_i$, and hence
in $\MMod\cs$ we have that
$\Hom\big[Y\big({\cs^{\leq i-1}}\big),\cy(C)\big]=0$.
But by
\cite[Notation~\ref{N21.-100}]{Neeman18A}
this rewrites as
\[
0\eq\Hom_{\MMod\cs}^{}\big[Y\big({\cs^{\leq i-1}}\big),\cy(C)\big]
\eq
\Hom_\ct^{}\big[F\big({\cs^{\leq i-1}}\big),C\big]\ ,
\]
and the vanishing tells us that
$C\in F\big({\cs^{\leq i-1}}\big)^\perp$. 
\eprf

\pro{P95.5}
Let $\cs$ be a triangulated category with a bounded-above
\tstr\ $\tst\cs$, and let the metric $\{\cm_i,\,i\in\nn\}$
be as in Example~\ref{E95.1}. Let $F:\cs\la\ct$ be
a good extension with respect to the metric,
and write $\wt\cs$ for the triangulated
subcategory $\wt\cs=\fl'(\cs)$ of the category $\ct$.

Then the triangulated category $\wt\cs$
has a \tstr\ $\tst{\wt\cs}$,
and the following holds:
\be
\item
The relation between the category $\wt\cs$ and
its triangulated subcategory $\cs$ is  
that the essential image of
$F$ is the category ${\wt\cs}^-$. We remind the reader,
this means:
\[
F(\cs)\eq{\wt\cs}^-\eq\bigcup_{i\in\nn}{\wt\cs}^{\leq i}\ .
\]
Moreover: the \tstr\ $\tst\cs$,
on the subcategory $\cs\cong {\wt\cs}^-$,
is just the inverse
image, under the fully faithful
functor $F$, of the restriction to $F(\cs)$ of the \tstr\
$\tst{\wt\cs}$ on the category $\wt\cs$.
\item
If $X_*$ is a 
Cauchy sequence in $\cs$, and 
$X\in\wt\cs$ is given by $X=\hoco F(X_*)$,
then the map 
$F\big(\colim\,X_*^{\leq n}\big)\la X^{\leq n}$ is an
isomorphism.
\setcounter{enumiv}{\value{enumi}}
\ee
\epro

\prf
To prove the existence of a \tstr\ $\tst{\wt\cs}$
as in (i),
it suffices to show 
that the embedding
$F\big(\cs^{\leq0}\big)\subset F(\cs)\subset\wt\cs$
makes $F\big(\cs^{\leq0}\big)$ into the aisle of a \tstr\
on $\wt\cs$. We clearly have that
$\T F\big(\cs^{\leq0}\big)\subset F\big(\cs^{\leq0}\big)$ and that
$F\big(\cs^{\leq0}\big)* F\big(\cs^{\leq0}\big)\subset F\big(\cs^{\leq0}\big)$.
What needs
showing is that, for any object
$X\in\wt\cs$, there exists a triangle
$A\la X\la C$ with $A\in F\big(\cs^{\leq0}\big)$ and with
$C\in F\big(\cs^{\leq0}\big)^\perp$.

Choose an object $X\in\wt\cs=\fl'(\cs)$, and express it as
$X=\hoco F(X_*)$ for some Cauchy sequence $X_*$ in
$\cs$. By Lemma~\ref{L95.2} we may, after replacing $X_*$ by
a subsequence, assume that the maps in the sequence
$X_1^{\leq0}\la X_2^{\leq0}\la X_3^{\leq0}\la\cdots$
are all isomorphisms. Now: by Lemma~\ref{L95.2}
the sequence $X_*^{\geq1}$ is also Cauchy,
and there is a canonical map of
Cauchy sequences $\ph_*:X_*\la X_*^{\geq1}$. Passing to
homotopy colimits of Cauchy sequences,
which exist in $\ct$ because
the embedding $F:\cs\la\ct$ is a
good extension, we are given
that
$X=\hoco F(X_*)$ and we set $C=\hoco F\big(X_*^{\geq1}\big)$.
And we choose in $\ct$
a (non-unique) morphism $\ph:X\la C$ rendering commutative
the square
\[\xymatrix@C+20pt{
F(X_*) \ar[d]^\beta\ar[r]^-{F(\ph_*)} & F(X_*^{\geq1})\ar[d]^\gamma \\
X\ar[r]^-\ph & C 
}\]
Complete the morphism $\ph:X\la C$ to a triangle
$A\la X\stackrel\ph\la C\la\T A$ 
in the category $\ct$. Proposition~\ref{P93.9}
asserts first that we may replace the sequence of
morphisms $\ph_*:X_*\la X_*^{\geq1}$ by a subsequence,
then extend in $\cs$ to a Cauchy sequence
of triangles
$A_*\la X_*\stackrel{\ph_*}\la X_*^{\geq1}\la\T A_*$.
And then Proposition~\ref{P93.9}(iv) asserts further that
we can do it in such a way that
the commutative square above 
extends, in $\ct$, to
a commutative diagram
\[\xymatrix@C+20pt{
F(A_*) \ar[d]^\alpha\ar[r] & F(X_*)\ar[r]^-{F(\ph_*)}\ar[d]^\beta &
 F(X_*^{\geq1})\ar[d]^\gamma \ar[r] &\T F(A_*)\ar[d]^{\T\alpha} \\
A\ar[r]^-\ph & X\ar[r] & C\ar[r] &\T A
}\]
And by Proposition~\ref{P93.9}(v)
\be
\setcounter{enumi}{\value{enumiv}}
\item
There is an isomorphism
$\hoco F(A_*)\la A$ in the category $\ct$,
which exhibits the pre-triangle
$\cy(A)\la\cy(X)\stackrel{\cy(\ph)}\la\cy(C)\la\T\cy(A)$
as naturally isomorphic to
the strong triangle
obtained by taking
the colimit in $\MMod\cs$ of the sequence 
$Y(A_*)\la Y(X_*)\stackrel{Y(\ph_*)}\la Y(X_*^{\geq1})\la\T Y(A_*)$.
\setcounter{enumiv}{\value{enumi}}
\ee

And the whole point here is that the extension of
$\ph_*:X_*\la X_*^{\geq1}$, to a Cauchy sequence of triangles
$A_*\la X_*\stackrel{\ph_*}\la X_*^{\geq1}\la\T A_*$
in the category $\cs$,
is canonically unique. It must be
isomorphic in $\cs$ to the sequence
$X_*^{\leq0}\la X_*\stackrel{\ph_*}\la X_*^{\geq1}\la\T X_*^{\leq0}$,
and in the beginning of the second paragraph of
the proof we made sure, by passing to a subsequence,
that the morphisms in the sequence $X_*^{\leq0}$ were all
isomorphisms. Hence $A=\hoco F\big(X_*^{\leq0}\big)$
is isomorphic to
$F\big(X_1^{\leq0}\big)\in F\big(\cs^{\leq0}\big)$.
Thus $A$ belongs to $F\big(\cs^{\leq0}\big)$.
And Lemma~\ref{L95.4} guarantees 
that $C\in\cl'_1\subset\fl'(\cs)$
belongs to $F\big(\cs^{\leq0}\big)^\perp$. This
concludes most of the proof: we have
shown both that
$F\big(\cs^{\leq0}\big)$ is the aisle of a
\tstr\ on $\wt\cs$, and that, with respect to this
\tstr, we have that
$X^{\leq 0}=\colim\,F\big(X_*^{\leq0}\big)$.

Now: because $\wt\cs^{\leq0}=F\big({\cs^{\leq0}}\big)$,
by shifting we have that
$\wt\cs^{\leq n}=F\big({\cs^{\leq n}}\big)$.
Hence
\[
\wt\cs^-\eq\bigcup_{n\in\nn}\wt\cs^{\leq n}\eq
\bigcup_{n\in\nn}F\big({\cs^{\leq n}}\big) \eq
F(\cs)\ .
\]
This completes the proof of (i).

Next for (ii): by shifting we may assume
$n=0$. Now note: in the course of the proof of
(i) we exhibited in $\wt\cs$ a triangle
$A\la X\la C$, with $A\in F\big(\cs^{\leq0}\big)=\wt\cs^{\leq0}$
and with
$C\in F\big(\cs^{\leq0}\big)^\perp=\wt\cs^{\geq1}$.
And the $A$ we exhibited was
$\colim\,F\big(X_*^{\leq0}\big)=F\big(\colim\,X_*^{\leq 0}\big)$. 
\eprf

Applying Corollary~\ref{C93.95} to
the metrics of this section
yields:

\lem{L95.-15}
Let $\cs$ be a triangulated category
with a bounded-above \tstr\
$\tst\cs$.
Assume in addition that
\be
\item
The category
$\cs$ has the property that, for any object
$E\in\cs$, the coproduct of $\coprod_{n\in\nn}E$ exists
in $\cs$. By this we mean: the coproduct of
countably many copies of the object
$E$ exists.
\item
There exists an integer
$r>0$ such that,
for every object $E\in\cs^{\geq r}$,
the coproduct $\coprod_{n\in\nn}E$ belongs
to $\cs^{\geq0}$.
\item
We are given a good extension $F:\cs\la\ct$,
such that $F$ respects
the coproducts in $\cs$
of the form
$\coprod_{n\in\nn}E$,
whose existence is guaranteed in (i).
\ee
Then the embedding $\cs\la\wt\cs=\fl'(\cs)$
is a good extension.
\elem

\prf
This just comes by applying Corollary~\ref{C93.95},
which is stated for general
very good metrics $\{\cm_i,\,i\in\nn\}$,
to the special case where $\cm_i=\cs^{\leq i}$ for
a bounded-above \tstr\ on $\cs$;
see Example~\ref{E95.1}.
\eprf

And now for corollaries of Proposition~\ref{P95.5}.

\cor{C95.6}
Let the notation and hypotheses be
as in Proposition~\ref{P95.5}.
If $A_*$ is a Cauchy sequence in $\cs$ and
$X=\hoco F(A_*)$, then
\be
\item
The inverse image under the fully
faithful functor $F$, of the sequence
$X^{\leq1}\la X^{\leq2}\la X^{\leq3}\la\cdots$,
is a Cauchy sequence in $\cs$ which we will
criminally denote
$X^{\leq(*)}$. The notational crime
is that we are identifying $\cs$ with
its essential image under $F$.
\item
The Cauchy sequences $A_*$ and $X^{\leq(*)}$ are
ind-isomorphic.   

To spell this out: there exists a Cauchy sequence
$B_*$ such that the odd subsequence $B_{2(*)-1}$ is
a subsequence of $A_*$ and the even subsequence
$F(B_{2(*)})$ is a subsequence of $X^{\leq(*)}$.
\ee
\ecor

\prf
First of all:
the objects $X^{\leq i}$ belong to
$\wt\cs^{\leq i}\subset F(\cs)$,
and hence the sequence
$X^{\leq(*)}$ is in the essential image of the
fully faithful functor $F$.
And the Cauchyness of the sequence is
by Lemma~\ref{L95.2}, after all the functor
$(-)^{\leq n}$ takes the sequence to one that eventually
stabilizes.

Now for (ii).
Let $B_1=A_1$. We construct the sequence $B_*$
by induction.

Suppose we are given $B_{2n-1}$ and want to construct
$B_{2n}$. By induction $B_{2n-1}=A_m\in\cs$ for some $m$,
and it belongs to $\cs^{\leq N}$ for some $N>m+2n-1$. Hence
the map $F(B_{2n-1})=F(A_m)\la\hoco F(A_*)=X$ must factor
uniquely through $X^{\leq N}$. We choose a $B_{2n}$ with
$F(B_{2n})=X^{\leq N}$, and the map
$F(B_{2n-1})\la F(B_{2n})=X^{\leq N}$
has a unique lifting to a morphism $B_{2n-1}\la B_{2n}$
in the category $\cs$.

Now assume $B_{2n}$ has been defined, with
$F(B_{2n})=X^{\leq N}$. By the Cauchyness of $A_*$ the
sequence $A_*^{\leq N}$ must stabilize, and by
Proposition~\ref{P95.5}(ii) there
exists some integer $M>2n+N+1$ at which
we reach stability and $F(A_M)^{\leq N}=X^{\leq N}$.
This gives us a map $F(B_{2n})=F(A_M)^{\leq N}\la F(A_M)$.
We set $B_{2n+1}=A_M$ and let $B_{2n}\la B_{2n+1}$ be the
unique preimage in $\cs$ of the morphism
$F(B_{2n})\la F(A_M)$ above.

This defines the sequence $B_*$, with the odd subsequence
and the even subsequence as specified. And the fact that
the sequence $B_*$ is Cauchy is again by Lemma~\ref{L95.2}:
the functor $(-)^{\leq n}$ takes the sequence $B_*$ that
we have constructed to one that eventually
stabilizes.
\eprf

\cor{C95.6.5}
Let the notation and hypotheses
be as in Proposition~\ref{P95.5}.
For any integer $i>0$
the subcategory $\cl'_i\subset\fl'(\cs)=\wt\cs$, of
Definition~\ref{D93.11}, turns out
to be equal to $\wt\cs^{\geq i}\subset\wt\cs$.
\ecor

\prf
The inclusion $\cl'_i\subset\wt\cs^{\geq i}$
was proved back in Lemma~\ref{L95.4}; we need to prove
the reverse inclusion. Therefore let $X$ be an
object in $\wt\cs^{\geq i}$.

Because $X$ belongs to $\wt\cs=\fl'(\cs)$ there
exists a Cauchy sequence $A_*$ in $\cs$ with
$X=\hoco F(A_*)$. But by Corollary~\ref{C95.6}
the Cauchy sequences $F(A_*)$ and $X^{\leq(*)}$
are ind-isomorphic, and hence $X=\hoco F(A_*)=\hoco X^{\leq(*)}$.
And because $X\in\wt\cs^{\geq i}$ we have that
$X^{\leq n}\in\cs^{\geq i}$ for all $n\in\nn$. This
expresses $X$ as the homotopy colimit of
a Cauchy sequence $X^{\leq(*)}$ contained in
$\cm_i=\cs^{\geq i}$, and therefore $X\in\cl'_i$.
\eprf

\section{The situation in which the metrics of
  Section~\protect{\ref{S95}} are excellent}
\label{S96}

Let $\cs$ be a triangulated category with
a very good metric $\{\cm_i,\,i\in\nn\}$,
and let $F:\cs\la\ct$ be good extension with
respect to the metric.
Proposition~\ref{P94.3} gave us a necessary and sufficient
criterion for the metric to be excellent.
Specializing this to the very good metrics
of Example~\ref{E95.1}, and then narrowing
attention to the ones possessing good extensions
$F:\cs\la\ct$ as in Example~\ref{E95.13}, we have:

\pro{P96.1}
Suppose $\ct$ is a triangulated category, with
countable coproducts and a 
nondegenerate \tstr\ $\tst\ct$.
Assume further that there exists an integer
$r>0$ such that, for any countable set
of objects $\{X_i,i\in\nn\}$
contained in $\ct^{\geq r}$, we have
that $\coprod_{i\in\nn}X_i$ is contained in
$\ct^{\geq0}$.
Then 
the very good metric $\{\ct^-\cap\ct^{\geq i},\,i\in\nn\}$,
on the triangulated category $\ct^-$,
is excellent 
if and only if there exists an integer $n>1$ such that
every object $B\in\ct^-$ admits,
in the category $\ct$, a distinguished
triangle $A\la B\la C$ with $A\in\ct^{\geq1}$ and with
$C\in\big(\ct^{\geq n}\big)^\perp$.
\epro

\prf
The proof is a direct application
Proposition~\ref{P94.3}.
We first observe that $\cm_i=\ct^-\cap\ct^{\geq i}$
is closed in $\ct^-$ under direct summands, and
hence the hypotheses of
Proposition~\ref{P94.3} are fulfilled. The
necessity of the condition
is by setting the integer
$m>0$ to be $m=1$, letting
$n>m$ be the corresponding $n$
of Proposition~\ref{P94.3}.
Now Example~\ref{E95.13} tells
us that $\cl'(\cs)=\ct$,
and Corollary~\ref{C95.6.5} computes
for us that $\cl'_i=\ct^{\geq i}$.
And the existence
of the triangle $A\la B\la C$ is exactly as in
Proposition~\ref{P94.3}.

The simplifying feature, in dealing with metrics
coming from \tstr{s}, is that the case $m=1$
suffices by shifting. For any $m>0$, apply to
the object $\T^{m-1}B$ the condition for $m=1$. We
deduce a triangle $A\la\T^{m-1}B\la C$ with
$A\in\ct^{\geq1}$ and with
$C\in\cap\big(\ct^{\geq n}\big)^\perp$.
And from this we deduce a triangle
$\T^{1-m}A\la B\la \T^{1-m}C$
with $\T^{1-m}A\in\ct^{\geq m}$ and
with $\T^{1-m}C\in\cap\big(\ct^{\geq m+n-1}\big)^\perp$.
This proves that, for any $m>0$, the integer
$m+n-1$ works and satisfies the condition
of Proposition~\ref{P94.3}.
\eprf

For $\ct^b\subset\ct$ the statement simplifies a little.

\pro{P96.3}
Suppose $\ct$ is a triangulated category, with
countable coproducts and a 
nondegenerate \tstr\ $\tst\ct$.
Assume further that there exists an integer
$r>0$ such that, for any countable set
of objects $\{X_i,i\in\nn\}$
contained in $\ct^{\geq r}$, we have
that $\coprod_{i\in\nn}X_i$ is contained in
$\ct^{\geq0}$.

Then the
very good metric $\{\ct^b\cap\ct^{\geq i},\,i\in\nn\}$,
on the triangulated category $\ct^b$,
is excellent 
if and only if there exists an integer $n>1$ such that
every object $B\in\ct^b\cap\ct^{\geq0}$ admits,
in the category $\ct$, a distinguished
triangle $A\la B\la C$ with $A\in\ct^{\geq1}$ and with
$C\in\ct^+\cap\big(\ct^{\geq n}\big)^\perp$.
\epro

\prf
We first of all assert:
\be
\item
The integer $n>1$ works for any $B\in\ct^b$, not
only for $B\in\ct^b\cap\ct^{\geq0}$. That
is: with $n>1$ as in the hypotheses of
the Proposition and $B\in\ct^b$ arbitrary, there
exists a triangle $A\la B\la C$ in $\ct^+$
with $A\in\ct^{\geq1}$ and with
$C\in\ct^+\cap\big(\ct^{\geq n}\big)^\perp$.
\ee
Let us begin by proving (i).

Because $B$ is assumed to belong to $\ct^b\subset\ct^+$,
it must be in $\ct^{\geq-\ell}$ for some $\ell>0$.
And now we proceed by induction to reduce $\ell$.
Put $B=A_0$, and assume that, for some $0\leq i\leq\ell$,
we have
constructed
a sequence $A_i\la A_{i-1}\la\cdots\la A_1\la A_0=B$,
with each
$A_j\in\ct^b\cap\ct^{\geq-\ell+j}$. Then we apply
the hypothesis of the current Proposition,
to $\T^{-\ell+i}A_i\in\ct^b\cap\ct^{\geq0}$, and
produce a triangle $A_{i+1}\la A_i\la C_{i+1}$
with $A_i\in\ct^{\geq-\ell+i+1}$ and with
$C_{i+1}\in\ct^+\cap\big(\ct^{\geq n-\ell-i}\big)^\perp\subset\ct^+\cap\big(\ct^{\geq n}\big)^\perp$.
We have that $A_{\ell+1}$ belongs to $\ct^b\cap\ct^{\geq1}$,
while $C_{i+1}\in\ct^+\cap\big(\ct^{\geq n}\big)^\perp$
for all $0\leq i\leq\ell$. And a
repeated application of the octahedral axiom
gives that, in the triangle $A_{\ell+1}\la B\la C$,
we have $C\in\ct^+\cap\big(\ct^{\geq n}\big)^\perp$.

OK: now that we have proved (i) we are ready to apply
Proposition~\ref{P94.3}. First of all note that,
because $\cm_i=\ct^b\cap\ct^{\geq i}$ is closed
in $\ct^b$ under direct summands, the hypotheses
of Proposition~\ref{P94.3} are all satisfied.

The
necessity of the condition
in the Proposition  
is again by setting the integer
$m>0$ to be $m=1$, letting
$n>m$ be the corresponding $n$
of Proposition~\ref{P94.3}.
As stated in the current Proposition,
the existence
of the triangle $A\la B\la C$ as in
Proposition~\ref{P94.3} is only assumed
for $B\in\ct^b\cap\ct^{\geq0}$, but in
(i) we proved that it holds for all $B\in\ct^b$.

And the remainder of the proof is as in
Proposition~\ref{P96.1}:
once again
the simplifying feature, in dealing with metrics
coming from \tstr{s}, is that the case $m=1$
suffices by shifting. For any $m>0$, apply to
the object $\T^{m-1}B$ the condition for $m=1$. We
deduce a triangle $A\la\T^{m-1}B\la C$ with
$A\in\ct^{\geq1}$ and with
$C\in\ct^+\cap\big(\ct^{\geq n}\big)^\perp$.
And from this we deduce a triangle
$\T^{1-m}A\la B\la \T^{1-m}C$
with $\T^{1-m}A\in\ct^{\geq m}$ and
with $\T^{1-m}C\in\ct^+\cap\big(\ct^{\geq m+n-1}\big)^\perp$.
This proves that, for any $m>0$, the integer
$m+n-1$ works and satisfies the condition
of Proposition~\ref{P94.3}.
\eprf

\exm{E96.5}
Let $\ct$ be a weakly approximable triangulated category,
and let $\tst\ct$ be any compactly generated \tstr\
in the preferred equivalence class.
The category $\ct$
has products, and
\cite[Lemma~4.2.1]{Canonaco-Neeman-Stellari24}
establishes that there exists an integer $r>0$ such
that, given any set $\{X_\lambda,\,\lambda\in\Lambda\}$ of
objects in $\ct^{\leq0}$, the product
$\prod_{X_\lambda,\,\lambda\in\Lambda}X_\lambda$
belongs to $\ct^{\leq r}$.
Applying Example~\ref{E95.13}(ii) to the category
$\ct\op$ with the opposite
\tstr, we obtain on 
$\big(\ct^b\big)\op$ the very good metric
$\{\cm_i,\,i\in\nn\}$, where 
$\cm\op_i=\big(\ct^b\big)\op\cap\big(\ct^{\leq-i}\big)\op$.
The inclusion $\cs=\big(\ct^b\big)\op\la\big(\ct^-\big)\op$
is a good extension with respect to the metric,
with $\fl'(\cs)=\big(\ct^-\big)\op$.

This puts us in the situation where the hypotheses
of Proposition~\ref{P96.3} hold. The Proposition
provides necessary and
sufficient conditions for the metric
$\cm\op_i=\big(\ct^b\big)\op\cap\big(\ct^{\leq-i}\big)\op$
on the category $\big(\ct^b\big)\op$ to be excellent.

Let $G$ be a compact generator for $\ct$, and as
in \cite[Proposition~2.6]{Neeman24}
choose an integer $n>0$ such that
\be
\item
$\Hom_\ct^{}\big(G,\ct^{\leq-n}\big)=0$ and $G\in\ct^{\leq A}$.
\item
For any object $B\in\ct^{\leq 0}$, there exists
in $\ct$ a triangle $C\la B\la A$ with
$C\in\ogenu G{}{-n,n}$ and with $A\in\ct^{\leq-1}$
\ee
This triangle does the job, proving that this
metric on $\big(\ct^b\big)\op$
is excellent. After all: in the
triangle of (ii) the object
$A$ belongs to
$\ct^{\leq-1}$, and by (i) the
object $C$ belongs to
$^\perp\big(\ct^{\leq-2n}\big)\cap\ct^{\leq 2n}$. 
\eexm

\exm{E96.7}
Let $\ct$ be a weakly approximable triangulated category. In
Example~\ref{E96.5} we learned that
the triangulated category $\cs=\big(\ct^b\big)\op$
has an excellent metric
$\cm\op_i=\big(\ct^b\big)\op\cap\big(\ct^{\leq-i}\big)\op$.
Moreover: we know from
Example~\ref{E95.13}(ii) that
$\fl'(\cs)=\big(\ct^-\big)\op$, while 
Corollary~\ref{C95.6.5} teaches us that
the subcategory $\cl'_i\subset\fl'(\cs)$
turn out to be $\cl'_i=\big(\ct^{\leq-i}\big)\op$.

Proposition~\ref{P97.1} teaches us that the
triangulated category
$\fs\big((\ct^b)\op\big)\op$,
with the metric $\cn\op_i$, is also a triangulated
category with an excellent metric.
But Lemmas~\ref{L93.-5} and \ref{L94.-9} give us
explicit formulas to compute $\fs\big((\ct^b)\op\big)\op$
and its metric. The formula for the
category $\fs\big((\ct^b)\op\big)\op$, as a subcategory
of $\ct=\big(\ct\op\big)\op$, comes down to 
\begin{eqnarray*}
\fs\big((\ct^b)\op\big)\op&=&
\ct^-\cap\bigcup_{i\in\nn}\left[\Big(\big(\ct^{\leq-i}\big)\op\Big)^\perp\right]\op\\
&=& \ct^-\cap\bigcup_{i\in\nn}{^\perp\big(\ct^{\leq-i}\big)}\\
&=&\bigcup_{i,j\in\nn}\left(\ct^{\leq j}\cap{^\perp\big(\ct^{\leq-i}\big)}\right)
\end{eqnarray*}
while the formula for the excellent
metric $\{\cn\op_\ell,\,\ell\in\nn\}$
is 
\begin{eqnarray*}
\cn\op_\ell&=&
\ct^{\leq-\ell}\cap\bigcup_{i\in\nn}\left[\Big(\big(\ct^{\leq-i}\big)\op\Big)^\perp\right]\op\\
&=& \ct^{\leq-\ell}\cap\bigcup_{i\in\nn}{^\perp\big(\ct^{\leq-i}\big)}\\
&=&\bigcup_{i\in\nn}\left(\ct^{\leq-\ell}\cap{^\perp\big(\ct^{\leq-i}\big)}\right)\ .
\end{eqnarray*}
Next we ask ourselves if these expressions can be simplified.

Suppose $B$ is an object in
$\ct^{\leq j}\cap{^\perp\big(\ct^{\leq-i}\big)}$,
and let the integer $n>0$ be as in
Example~\ref{E96.5}~(i) and (ii).
By \cite[Corollary~2.2.1]{Neeman24} there
exists a triangle $A\la B\la C$ in $\ct^-$,
with $A\in\ogenu G{}{-i-n+1,j+n}$
and with $C\in\ct^{\leq -i}$. But then the
map from $B\in{^\perp\big(\ct^{\leq-i}\big)}$ to
$C\in\ct^{\leq-i}$ must vanish, and $B$ is
a direct summand of $A\in\ogenu G{}{-i-n+1,j+n}$.
Therefore $B\in\ogenu G{}{-i-n+1,j+n}$.

By the above $\fs\big((\ct^b)\op\big)\op$ is contained
in the union $\bigcup_{m\in\nn}\ogenu G{}{-m,m}$.
But,
by the inclusions of
Example~\ref{E96.5}~(i), we also have
$\ogenu G{}{-m,m}\subset\ct^{\leq m+n}\subset\ct^-$,
and also
$\ogenu G{}{-m,m}\subset{^\perp\big(\ct^{\leq-m-n}\big)}$.
We deduce the equality
\[
\fs\big((\ct^b)\op\big)\op\eq\bigcup_{m\in\nn}\ogenu G{}{-m,m}\ .
\]
\eexm

Let us give this category a name:

\dfn{D96.9034}
Let $\ct$ be a weakly approximable triangulated category.
Choose a compact generator, and define
\[
\tsb\eq\bigcup_{m\in\nn}\ogenu G{}{-m,m}\ .
\]
This category possesses an excellent metric
$\{\cn\op_\ell,\,\ell\in\nn\}$, whose formula comes to
\[
\cn\op_\ell=\tsb\cap\ct^{\leq-\ell}\ .
\]
\edfn

Of course: the way we arrived at this category
is indirect, it showed up
$\tsb=\fs\big((\ct^b)\op\big)\op$ with the induced
excellent metric.

\bibliographystyle{amsplain}
\bibliography{stan}

\end{document}